\documentclass[10pt,leqno,oneside,letterpaper]{amsart}
\usepackage[letterpaper,left=1.75truein, top=1truein]{geometry}
\usepackage{amsmath,amsfonts,amssymb,amscd,amsthm,amsbsy,epsf,graphics}
\usepackage{xy}
\xyoption{all}

\usepackage{color}

\usepackage[colorlinks,linkcolor=blue,citecolor=blue]{hyperref}
\usepackage{epsfig}
\usepackage{graphicx}
\usepackage{ifpdf}
\usepackage{epstopdf}
\allowdisplaybreaks
\usepackage{longtable}
\makeatletter
\def\thm@space@setup{\thm@preskip=0pt \thm@postskip=0pt}
\makeatother
\usepackage[nodisplayskipstretch]{setspace}
\newtheoremstyle{mystyle}
  {3pt} 
  {3pt} 
  {\itshape} 
  {} 
  {\bfseries} 
  {.} 
  {.5em} 
  {} 
\theoremstyle{mystyle}
\newtheorem{thm}{Theorem}[section]

\newtheorem{lemma}[thm]{Lemma}

\theoremstyle{definition}

\def\qed{{\hspace{2mm}{\small $\square$}}}

\theoremstyle{cases}

\numberwithin{subcase}{section}
\numberwithin{equation}{section}

\def\sfrac#1#2{\kern.1em\raise.5ex\hbox{$#1$}
    \kern-.1em/\kern-.05em\lower.25ex\hbox{$#2$}}
\def\A{{\mathcal A}}

\def\T{{\mathcal T}}
\def\wT{\widetilde{\T}}

\def\wA{\widetilde{\A}}

\def\sfrac#1#2{\kern.1em\raise.5ex\hbox{$#1$}
        \kern-.1em/\kern-.05em\lower.25ex\hbox{$#2$}}

 \def\d{{\delta}}

 \def\a{{\alpha}}
 
 \def\p{{\partial}}
 
 \def\ra{{\rightarrow}}

 \def\lra{{\longrightarrow}}

 \def\z{{\mathbb Z}}
 \def\2{{\mathbb Z_2}}
 \def\q{{\mathbb Q}}

 \def\sl2{{SL(2,\mathbb C)}}

 \def\pf{{\noindent{\bf Proof.\hspace{2mm}}}}
 
 \def\sm{{{\mbox{\footnotesize M}}}}
  \def\sl{{{\mbox{\footnotesize L}}}}
 
\def\pq{{{\mbox{\tiny T(p,q)}}}}

 \def\ab{{{\mbox{\tiny T(a,b)}}}}
\def\aa{{{\mbox{\tiny T(a,2)}}}}

\def\sk{{{\mbox{\tiny K}}}}
\def\su{{{\mbox{\tiny U}}}}

\begin{document}

\title{The AJ conjecture and connected sums of torus knots}

\author{Xingru Zhang}
\address{Department of Mathematics, University at Buffalo, Buffalo, NY, 14214-3093, USA.}
\email{xinzhang@buffalo.edu}

\maketitle

\begin{abstract}The set of isotopy classes of
 nontrivial torus knots $T(p,q)$ in $S^3$  is in bijection with the set  of
coprime integer pairs $(p,q)$  satisfying $|p|>q\geq 2$.
We verify the AJ conjecture for the connected sums
$T(p,q)\# T(a,b)$ when $p$ and $a$ have the same sign.
Notably, in cases where $pq=ab$ but $p\ne a$,
the recurrence polynomial
$\a(t,\sm,\sl)$ of $T(p,q)\#T(a,b)$ has repeated factors involving
the  variable  $\sl$  after evaluation at $t=-1$. These appear to be
 the first  examples of  knots exhibiting this phenomenon.
Therefore,  the AJ conjecture
requires a slight  modification to accommodate this possibility.
  \end{abstract}

\section{Introduction}

The AJ conjecture, formulated by Garoufalidis  \cite[Conjecture 2]{G}, predicts a fundamental connection between two important knot invariants for knots in $S^3$: the A-polynomial and the colored Jones polynomial.

 For a knot $K$ in $S^3$,
its (normalized) $A$-polynomial $A_\sk(\sm,\sl)$ is a two-variable polynomial
in $\z[\sm,\sl]$ with no repeated factors and with relatively prime integer coefficients,
which is uniquely associated with $K$ up to a sign (see \cite{CCGLS}).
Note that $A_\sk(\sm,\sl)$ always contains the factor $\sl-1$.

For a knot $K$ in $S^3$ and each positive integer $n$,  the  $n$-colored Jones polynomial  $J_{\sk}(n)$ is a one-variable Laurent polynomial in $\z[t^{\pm1}]$ with the zero framing (see \cite{Tu}). In this paper, $J_{\sk}(n)$ is normalized  such that for the unknot $U$,
\begin{equation}\label{[n]}
J_\su(n) =[n]:= \frac{t^{2n}-t^{-2n}}{t^2-t^{-2}}.
\end{equation}
By defining $J_{\sk}(-n) := -J_{\sk}(n)$ and $J_\sk(0):=0$, one may  treat $J_{\sk}(n)$ as
 a discrete function  $J_{\sk}(\cdot): \z \to \z[t^{\pm 1}]$. The quantum torus $\T$, a non-commutative ring with presentation
\[ \T = \q[t^{\pm 1}] \left<\sl^{\pm 1}, \sm^{\pm 1} \right> / (\sl \sm-t^2\sm \sl),\]
acts on   the set of discrete functions $f: \z \to \q[t^{\pm 1}]$ by
\begin{equation}\label{equ:operator}
(\sm f)(n) := t^{2n}f(n), \quad  (\sl f)(n) := f(n+1).
\end{equation}
 The set of
  annihilators of $J_{\sk}(n)$ in $\T$,  denoted by
  $${\mathcal A}_\sk := \{P \in \T \mid P J_{\sk}(n) = 0\},$$
  is   a left ideal of $\T$ called the \textit{recurrence ideal} of $K$ (as each operator $P\in {\mathcal A}_\sk$ defines
 a recurrence relation for $J_\sk(n)$ and vice versa).
  A fundamental result proven in  \cite{GaLe}  states that for every knot $K$ in $S^3$, $J_\sk(n)$ has a nontrivial recurrence relation, or equivalently, ${\mathcal A}_\sk$
  is not the zero ideal.

The ring $\T$ can be extended to a principal left ideal domain $\wT$ by adjoining inverses of polynomials in $t$ and $\sm$; that is, $\wT$ is the set of Laurent polynomials in $\sl$ with coefficients in $\q(t,\sm)$ (which is the set of
rational functions in $t$ and $\sm$ with coefficients in $\q$) under the product defined by
\begin{equation}\label{equ:product}
f(t,\sm) \sl^j \cdot g(t,\sm) \sl^k = f(t,\sm) g(t,t^{2j}\sm)\sl^{j+k}.
\end{equation}
The extended left ideal $\wA_\sk=\wT {\mathcal A_\sk}$ is then generated by a single non-zero element in $\wT$.
Such  a generator can be normalized to lie in ${\mathcal A}_\sk$ and takes  the form
\[\alpha_\sk(t,\sm,\sl) = \sum_{i=0}^d D_i \sl^i,\]
with $d>0$ minimum (among all non-zero annihilators of $J_\sk(n)$) and with  $D_0,\dots,D_d \in \z[t,\sm]$ being coprime in $\z[t,\sm]$. This polynomial $\alpha_\sk(t,\sm,\sl)$ is uniquely determined up to a sign and is called the (\textit{normalized}) \textit{recurrence polynomial} of $K$.

The original AJ conjecture   states that for every knot $K$ in $S^3$, $\alpha_\sk(-1,\sm,\sl)$ is equal to
$A_\sk(\sm,\sl)$  up to a non-zero factor depending on $\sm$ only.
Based on examples found in this paper, for some knots $K$, $\a_\sk(-1,\sm,\sl)$
may  contain repeated factors involving the variable $\sl$.
Consequently, the $AJ$ conjecture requires  a slight modification:

{\bf The AJ conjecture:} For every knot $K$ in $S^3$,
 $\a_\sk(-1,\sm,\sl)$, after removing all repeated factors, is equal to $A_\sk(\sm,\sl)$
 up to a non-zero factor depending on $\sm$ only.

The conjecture remains largely open, having been proven true only for some simple classes of knots, such as some $2$-bridge knots \cite{G, Ta, Le, GK, LT, LZ},  torus knots \cite{Hikami, Tran},  some
 pretzel knots \cite{LT}, and  some cabled  knots  over certain knots \cite{RZ, R,Tran2,Tran3,Tran4,D}.
 In this paper, we consider the AJ conjecture for connected sums of torus knots.

Recall that the set of isotopy classes of nontrivial torus knots $T(p,q)$ in $S^3$
is in one-to-one correspondence with the set of coprime integer pairs
$(p,q)$ satisfying  $|p|>q\geq 2$.
We show:
\begin{thm}\label{main result}
The $AJ$ conjecture holds for the connected sums $T(p,q)\#T(a,b)$
when $p$ and $ a$ have the same sign.
\end{thm}

As we will see,  for a knot $K=T(p,q)\#T(a,b)$  given as in Theorem \ref{main result},
its recurrence polynomial $\a_\sk(t,\sm,\sl)$, after evaluation at $t=-1$, has repeated factors involving the variable
$\sl$  exactly when $pq=ab$ but $p\ne a$.

 In Section 2, we list the formulas to be applied.
Theorem 1.1 is then proved in Sections 3-9, which  correspond to seven cases.

For  convenience, in this paper we often represent the $A$-polynomial
$A_\sk(\sm,\sl)$ as an element in $\z[\sm^{\pm 1}, \sl]$ that  differs from the normalized
$A$-polynomial only by a  factor  of an integer power of $\sm$.
Similarly, in each case of $K=T(p,q)\# T(a,b)$ that we consider,  we shall obtain
a recurrence polynomial $\a_\sk(t, \sm,\sl)$ up to normalization, i.e., up to multiplying
by a  non-zero  rational function in $t$ and $\sm$.

Our proof  procedure consists of the  following three steps:
\begin{enumerate}
\item  Produce a candidate minimum $\sl$-degree non-zero annihilator $\a_\sk(t, \sm, \sl)$
in $\wA_\sk$ of the form
$\displaystyle\a_\sk(t, \sm,\sl)=\sum_{i=0}^d D_i\sl^i $ with $D_i\in \q(t,\sm)$
 such that $\a_\sk(-1, \sm, \sl)$ is well-defined and has the same $\sl$-degree as $\a_\sk(t,\sm,\sl)$;
\item Verify that $\a_\sk(-1,\sm,\sl)$, after removal of all repeated factors, is equal to
the $A$-polynomial $A_\sk(\sm,\sl)$ up to multiplying  by  a non-zero  rational function in $\sm$;
\item  Prove that the $\sl$-degree $d$ of our candidate annihilator  $\a_\sk(t,\sm,\sl)$
is indeed  the minimum among
all non-zero annihilators of $J_\sk(n)$ in $\wA_\sk$.
\end{enumerate}
This  procedure is sufficient  to  verify the  AJ conjecture for the given knot.

\section{Some formulas}\label{sec:formulas}

The $A$-polynomial of  $T(p,q)$ is given by (see e.g.,\cite[Example 4.1]{Z})
\begin{equation}\label{Apoly of torus knots}
A_{\pq}(\sm, \sl)= \begin{cases}(\sl-1)(\sm^{2pq}\sl^2-1)&  \text{if $|p|>q>2$,}\\
(\sl-1)(\sm^{2p}\sl+1)&  \text{if $|p|>q=2$.}\end{cases}\end{equation}

For knots $K_1$ and $K_2$ in $S^3$, if $\displaystyle A_{\sk_1}(\sm,\sl)=\prod_{i\in I}(\sl \sm^{r_i}-\d_i)$
and $\displaystyle A_{\sk_2}(\sm,\sl)=\prod_{j\in J}(\sl \sm^{s_i}-\d_j)$, where $I$ and $J$ are some finite index sets,
$r_i$ and $s_j$ are integers, and $\d_i, \d_j\in \{1, -1\}$,
then the $A$-polynomial of  the connected sum $K_1\# K_2$ is given by (see \cite[Theorem 1.1]{S})
\begin{equation}\label{Apoly of sum}A_{\sk_1\#\sk_2}(\sm,\sl)=\text{Red}\Big[\prod_{(i,j)\in I\times J}(\sl\sm^{r_i+s_j}-\d_i\d_j)\Big]
\end{equation} where $\text{Red}[\cdot]$ denotes the  removal of  all
repeated factors. Plugging  (\ref{Apoly of torus knots})  into (\ref{Apoly of sum}),
we obtain the following lemma:
\begin{lemma}\label{lem:A-poly of sum}The $A$-polynomial of $T(p,q)\#T(a,b)$ is given by:
\begin{enumerate}
\item[(i)]   If  $|p|>q>2$, $|a|>b>2$, $pq\ne ab$, and $pq\ne -ab$,  then
\[ A_{\pq\#\ab}(\sm, \sl)=
 (\sl-1)(\sm^{2pq}\sl^2-1)(\sm^{2ab}\sl^2-1)(\sm^{2(pq+ab)}\sl^2-1).\]
\item[(ii)] If $|p|>q>2$, $|a|>b>2$,  and $pq=ab$, then
\[ A_{\pq\#\ab}(\sm, \sl)=
(\sl-1)(\sm^{2ab}\sl^2-1)(\sm^{4ab}\sl^2-1).\]
\item[(iii)] If $|p|>q>2$, $|a|>b>2$,  and $pq=-ab$, then
\[ A_{\pq\#\ab}(\sm, \sl)=
 (\sl^2-1)(\sm^{-2ab}\sl^2-1)(\sm^{2ab}\sl^2-1).\]
\item[(iv)]   If  $|p|>q>2$, $|a|>b=2$, $pq\ne 2a$, and $pq\ne -2a$,  then
\[A_{\pq\#\ab}(\sm, \sl)=
 (\sl-1)(\sm^{2pq}\sl^2-1)(\sm^{2a}\sl+1)(\sm^{2(pq+2a)}\sl^2-1).\]
\item[(v)]   If  $|p|>q>2$, $|a|>b=2$, and $pq=2a$,  then
\[ A_{\pq\#\ab}(\sm, \sl)=
 (\sl-1)(\sm^{4a}\sl^2-1)(\sm^{8a}\sl^2-1).\]
\item[(vi)]   If  $|p|>q>2$, $|a|>b=2$, and $pq=-2a$,  then
\[ A_{\pq\#\ab}(\sm, \sl)=
 (\sl^2-1)(\sm^{-4a}\sl^2-1)(\sm^{2a}\sl+1).\]
\item[(vii)]   If  $|p|>q=2$, $|a|>b=2$, $p\ne a$, and $p\ne -a$,  then
\[A_{\pq\#\ab}(\sm, \sl)=
 (\sl-1)(\sm^{2p}\sl+1)(\sm^{2a}\sl+1)(\sm^{2p+2a}\sl-1).\]
\item[(viii)]   If  $|p|>q=2$, $|a|>b=2$, $p= a$,  then
\[ A_{\pq\#\ab}(\sm, \sl)=
 (\sl-1)(\sm^{2p}\sl+1)(\sm^{4p}\sl-1).\]
\item[(ix)]   If  $|p|>q=2$, $|a|>b=2$, $p= -a$,  then
\[ A_{\pq\#\ab}(\sm, \sl)=
 (\sl-1)(\sm^{2p}\sl+1)(\sm^{-2p}\sl+1).\]
\end{enumerate}
\end{lemma}

Under the normalization (\ref{[n]}), the connected sum formula for the colored Jones polynomial  is given by
\begin{equation}\label{connected sum for CJ}
 [n]J_{\sk_1\#\sk_2}(n)=J_{\sk_1}(n)J_{\sk_2}(n).
\end{equation}

For a Laurent polynomial $f(t)\in \q[t^{\pm1}]$, let $\ell(f)$ and $\hbar(f)$
 denote the lowest and highest degrees of $f$ in $t$, respectively.
Clearly,  for any two non-zero $f(t)$ and $g(t)$ in $\q[t^{\pm1}]$,
  $\ell{fg}=\ell(f)+\ell(g)$ and $\hbar(fg)=\hbar(f)+\hbar(g)$.

\begin{lemma}\label{degree of J_T} $(1)$ If $p>q$, then
\begin{align*}&
\ell(J_{\pq}(n))=-pqn^2+pq+\frac{1}{2}(1-(-1)^{n-1})(p-2)(q-2),\\
&\hbar(J_{\pq}(n))=2(p+q-pq)|n|+2(pq-p-q).
\end{align*}
$(2)$
If $p<-q$, then \begin{align*}
&\ell(J_{\pq}(n))=2(p-q-pq)|n|+2(pq-p+q),\\
&\hbar(J_{\pq}(n))=-pqn^2+pq+\frac{1}{2}(1-(-1)^{n-1})(p+2)(q-2).\end{align*}
\end{lemma}

 The formula for $\ell(J_\pq(n))$ in part (1) of Lemma \ref{degree of J_T} is from \cite[Lemma 2.4]{Tran}, and the rest  of
Lemma \ref{degree of J_T}
is from \cite[Lemma 2.1]{RZ}.

Combining Lemma \ref{degree of J_T} with formula (\ref{connected sum for CJ}), we obtain the following
lemma.

\begin{lemma}\label{degree of connected sum}
\begin{enumerate}
\item If $p>q$ and  $a>b$, then
\begin{align*}
\ell([n]J_{\pq\#\ab}(n))&=-pqn^2+pq+\frac{1}{2}(1-(-1)^{n-1})(p-2)(q-2)
\\&\quad -abn^2+ab+\frac{1}{2}(1-(-1)^{n-1})(a-2)(b-2),\\
 \hbar ([n]J_{\pq\#\ab}(n))&=2(p+q-pq)|n|+2(pq-p-q)\\
& \quad +2(a+b-ab)|n|+2(ab-a-b).\end{align*}
\item If $p>q$ and  $a<-b$, then
\begin{align*}
\ell([n]J_{\pq\#\ab}(n))&=-pqn^2+pq+\frac{1}{2}(1-(-1)^{n-1})(p-2)(q-2)
\\&\quad+2(a-b-ab)|n|+2(ab-a+b),\\
\hbar([n]J_{\pq\#\ab}(n))&=2(p+q-pq)|n|+2(pq-p-q)\\
&\quad-abn^2+ab+\frac{1}{2}(1-(-1)^{n-1})(a+2)(b-2).\end{align*}
\item If $p<-q$ and  $a<-b$, then
\begin{align*}
\ell([n]J_{\pq\#\ab}(n))&
=2(p-q-pq)|n|+2(pq-p+q)\\
&\quad+2(a-b-ab)|n|+2(ab-a+b),\\
\hbar([n]J_{\pq\#\ab}(n))&
=-pqn^2+pq+\frac{1}{2}(1-(-1)^{n-1})(p+2)(q-2)\\
&\quad-abn^2+ab+\frac{1}{2}(1-(-1)^{n-1})(a+2)(b-2).\end{align*}
\end{enumerate}
\end{lemma}

By \cite[Lemma 2.5]{Tran}, we have
\begin{equation}\label{(p,2)}
J_{T(p,2)}(n+1)=-t^{-4pn-2p}J_{T(p,2)}(n)+t^{-2pn}[2n+1].
\end{equation}
Note that although only $p>2$ was considered  in \cite{Tran}, (\ref{(p,2)}) is valid for $p<-2$ as well.

Since $J_{T(p,2)}(n+1)=\sl J_{T(p,2)}(n)$ and $\sm=t^{2n}$,
 it follows  that
\begin{equation}\label{(p,2,n+1)}
(\sl+\sm^{-2p}t^{-2p})J_{T(p,2)}(n) =\sm^{-p}[2n+1].
\end{equation}
Applying (\ref{(p,2)}) twice,  we obtain
\begin{equation}\label{(p,2,n+2)}
\begin{aligned}
J_{T(p,2)}(n+2)&=-t^{-4pn-6p}J_{T(p,2)}(n+1)+t^{-2pn-2p}[2n+3]\\
&=-t^{-4pn-6p}\left(-t^{-4pn-2p}J_{T(p,2)}(n)+t^{-2pn}[2n+1]\right)
+t^{-2pn-2p}[2n+3]\\
&=t^{-8p(n+1)}J_{T(p,2)}(n)-t^{-6p(n+1)}[2n+1]+t^{-2p(n+1)}[2n+3].
\end{aligned}
\end{equation}
By   \cite[Lemma 2.1]{Tran},  when $|p|>q>2$, we have
\begin{equation}\label{relation of J_T(n+2)}J_{\pq}(n+2) = t^{-4pq(n+1)} J_{\pq}(n)+t^{-2pq(n+1)} \delta(p,q,n),
\end{equation}
where
\begin{equation}\label{equ:delta}
\begin{aligned}
\delta(p,q,n)&:= \frac{t^{2(p+q)(n+1)+2}+t^{-2(p+q)(n+1)+2} -t^{2(q-p)(n+1)-2}-t^{-2(q-p)(n+1)-2}}{t^2-t^{-2}}\\
&= \frac{\sm^{p+q}t^{2(p+q+1)}+\sm^{-p-q}t^{-2(p+q-1)} -\sm^{q-p}t^{2(q-p-1)}
-\sm^{p-q}t^{-2(q-p+1)}}{t^2-t^{-2}}.
\end{aligned}\end{equation}
Again, although only $p>q$ was considered in \cite{Tran}, (\ref{relation of J_T(n+2)}) and (\ref{equ:delta})
are  valid for $p<-q$ as well.

Since $J_{\pq}(n+2)=\sl^2J_{\pq}(n)$ and $\sm=t^{2n}$,  we have from (\ref{relation of J_T(n+2)})
that
\begin{equation}\label{an of J_T(n+2)}
(\sl^2-\sm^{-2pq}t^{-4pq})J_{\pq}(n) =\sm^{-pq}t^{-2pq} \delta(p,q, n).
\end{equation}

\section{Case $|p|>q>2$, $|a|>b>2$, $pq\ne ab$, $p$ and $a$ have the same sign}\label{sec:pqab}

 As described in the introduction, the first step of the proof procedure in each of the following sections
-- which correspond to seven cases of knots $K=T(p,q)\# T(a,b)$--is
 to produce a candidate annihilator $\a_\sk(t,\sm,\sl)$  in $\wA_\sk$.
 Applying (\ref{connected sum for CJ}) and (\ref{relation of J_T(n+2)}), we have
\begin{equation}\label{relation for pqab}
\begin{aligned}&
[n+2]J_{\pq\#\ab}(n+2)=
J_{\pq}(n+2)J_{\ab}(n+2)
\\&
=\left (t^{-4pq(n+1)} J_{\pq}(n)+t^{-2pq(n+1)} \delta(p,q,n)\right)\\&\;\;\;\;\left( t^{-4ab(n+1)} J_{\ab}(n)+
t^{-2ab(n+1)} \delta(a,b,n)\right)\\
&=t^{-4(pq+ab)(n+1)}J_{\pq}(n)J_{\ab}(n)+t^{-(4pq+2ab)(n+1)}
\delta(a,b,n)J_{\pq}(n)\\
&\;\;\;\;+t^{-(2pq+4ab)(n+1)}\delta(p,q,n)J_{\ab}(n)+t^{-2(pq+ab)(n+1)}
\delta(p,q,n)\delta(a,b,n)\\
&=  t^{-4(pq+ab)(n+1)}[n]J_{\pq\#\ab}(n)+t^{-(4pq+2ab)(n+1)}\delta(a,b,n)J_{\pq}(n)
\\&
\;\;\;\;+t^{-(2pq+4ab)(n+1)}\delta(p,q,n)J_{\ab}(n)+t^{-2(pq+ab)(n+1)}\delta(p,q,n)\delta(a,b,n).
\end{aligned}
\end{equation}
 Since $J_{\pq\#\ab}(n+2)= \sl^2J_{\pq\#\ab}(n)$ and
$t^{2n}=\sm$, we can transform (\ref{relation for pqab}) into
\begin{align*}
&\left([n+2]\sl^2- \sm^{-2pq-2ab}t^{-4pq-4ab}[n]\right)
J_{\pq\#\ab}(n)\\
&=\sm^{-2pq-ab}t^{-4pq-2ab}\delta(a,b,n)J_{\pq}(n)
+\sm^{-pq-2ab}t^{-2pq-4ab}\delta(p,q,n)J_{\ab}(n)\\
&\;\;\;\;+
\sm^{-pq-ab}t^{-2pq-2ab} \delta(p,q,n)\delta(a,b,n).
\end{align*}
Dividing from the left by the non-zero function
$\sm^{-2pq-ab}t^{-4pq-2ab}\delta(a,b,n)$,
 which is the coefficient of $J_{\pq}(n)$,  this equality becomes
$$\begin{array}{l}
\frac{\sm^{2pq+ab}t^{4pq+2ab}}{\d(a,b,n)}
\left([n+2]\sl^2- \sm^{-2pq-2ab}t^{-4pq-4ab}[n]\right)
J_{\pq\#\ab}(n)\\
=J_{\pq}(n)+\frac{\sm^{pq-ab}t^{2pq-2ab}\d(p,q,n)}{\d(a,b,n)}
J_{\ab}(n)+\sm^{pq}t^{2pq} \delta(p,q,n).
\end{array}$$
Applying the operator $(\sl^2-\sm^{-2pq}t^{-4pq})$ from the left
to both sides of the last equality, and then applying (\ref{an of J_T(n+2)}),
 (\ref{equ:operator}) (treating $\sm$ as $t^{2n}$),
 and (\ref{relation of J_T(n+2)}),  we obtain
\begin{small}\begin{equation}\label{3.2}
\begin{aligned}
&\left(\sl^2-\sm^{-2pq}t^{-4pq}\right)\frac{\sm^{2pq+ab}t^{4pq+2ab}}{\d(a,b,n)}
\left([n+2]\sl^2- \sm^{-2pq-2ab}t^{-4pq-4ab}[n]\right)
J_{\pq\#\ab}(n)
\\&
=\sm^{-pq}t^{-2pq} \delta(p,q,n)+\left(\sl^2-\sm^{-2pq}t^{-4pq}\right)\left(\frac{\sm^{pq-ab}t^{2pq-2ab}\d(p,q,n)}{\d(a,b,n)}
J_{\ab}(n)+\sm^{pq}t^{2pq} \delta(p,q,n)\right)
\\&
=\sm^{-pq}t^{-2pq} \delta(p,q,n)+
\frac{\sm^{pq-ab}t^{6pq-6ab}\d(p,q,n+2)}{\d(a,b,n+2)}
J_{\ab}(n+2)+\sm^{pq}t^{6pq} \delta(p,q,n+2)
\\&
\;\;\;\;-\frac{\sm^{-pq-ab}t^{-2pq-2ab}\d(p,q,n)}{\d(a,b,n)}
J_{\ab}(n)-\sm^{-pq}t^{-2pq} \delta(p,q,n)
\\&
=\sm^{-pq}t^{-2pq} \delta(p,q,n)+
\frac{\sm^{pq-ab}t^{6pq-6ab}\d(p,q,n+2)}{\d(a,b,n+2)}
\left(t^{-4ab(n+1)}J_{\ab}(n)+t^{-2ab(n+1)}\d(a,b,n)\right)
\\&
\;\;\;\;+\sm^{pq}t^{6pq} \delta(p,q,n+2)
-\frac{\sm^{-pq-ab}t^{-2pq-2ab}\d(p,q,n)}{\d(a,b,n)}
J_{\ab}(n)-\sm^{-pq}t^{-2pq} \delta(p,q,n)
\\&
=\left(\frac{\sm^{pq-3ab}t^{6pq-10ab}\d(p,q,n+2)}{\d(a,b,n+2)}
-\frac{\sm^{-pq-ab}t^{-2pq-2ab}\d(p,q,n)}{\d(a,b,n)}\right)J_{\ab}(n)
\\&
\;\;\;\;+\frac{\sm^{pq-2ab}t^{6pq-8ab}\d(p,q,n+2)\d(a,b,n)}{\d(a,b,n+2)}
+\sm^{pq}t^{6pq} \delta(p,q,n+2).
\end{aligned}
\end{equation}
\end{small}
Let $\displaystyle f(t,\sm)=\frac{\sm^{pq-3ab}t^{6pq-10ab}\d(p,q,n+2)}{\d(a,b,n+2)}
-\frac{\sm^{-pq-ab}t^{-2pq-2ab}\d(p,q,n)}{\d(a,b,n)}$,
which is the coefficient of $J_{\ab}(n)$ on the right-hand side  of  the last equality.
Then $f(t,\sm)\ne 0$.
This is easily verified
by specializing the function at $t=-1$ (using (\ref{equ:delta})):
\begin{equation}\label{equ:f(-1,m)}
 f(-1,\sm)
=\frac{(\sm^{pq-3ab}-\sm^{-pq-ab})
(\sm^{p+q}+\sm^{-p-q}-\sm^{q-p}-\sm^{p-q})}
{\sm^{a+b}+\sm^{-a-b}-\sm^{b-a}-\sm^{a-b}}
\end{equation}
which is non-zero
since $pq\ne ab$.

Applying the operator $(\sl^2-\sm^{-2ab}t^{-4ab})f^{-1}(t,\sm)$ from the left to both sides of (\ref{3.2}), and then using (\ref{an of J_T(n+2)}) and (\ref{equ:operator}), we obtain
\begin{small}
\begin{align*}
&\left(\sl^2-\sm^{-2ab}t^{-4ab}\right)f^{-1}(t,\sm) \left(\sl^2-\sm^{-2pq}t^{-4pq}\right)\frac{\sm^{2pq+ab}t^{4pq+2ab}}{\d(a,b,n)} \left([n+2]\sl^2- \sm^{-2pq-2ab}t^{-4pq-4ab}[n]\right) J_{\pq\#\ab}(n) \\
&\quad =\sm^{-ab}t^{-2ab}\d(a,b,n) + \left(\sl^2-\sm^{-2ab}t^{-4ab}\right)f^{-1}(t,\sm)\left(\frac{\sm^{pq-2ab}t^{6pq-8ab}\d(p,q,n+2)\d(a,b,n)}{\d(a,b,n+2)} + \sm^{pq}t^{6pq} \delta(p,q,n+2)\right) \\
&\quad = \sm^{-ab}t^{-2ab}\d(a,b,n) + f^{-1}(t,t^4\sm)\left(\frac{\sm^{pq-2ab}t^{10pq-16ab}\d(p,q,n+4)\d(a,b,n+2)}{\d(a,b,n+4)} + \sm^{pq}t^{10pq} \delta(p,q,n+4)\right) \\
& \qquad - \sm^{-2ab}t^{-4ab}f^{-1}(t,\sm)\left(\frac{\sm^{pq-2ab}t^{6pq-8ab}\d(p,q,n+2)\d(a,b,n)}{\d(a,b,n+2)} + \sm^{pq}t^{6pq} \delta(p,q,n+2)\right).
\end{align*}
\end{small}
Let $g(t,\sm)$ be the function on the right-hand side of  the last  equality.
Then $g(t, \sm)\ne 0$. Indeed, applying (\ref{equ:delta}) and
(\ref{equ:f(-1,m)}), we have (noting that
$\displaystyle\lim_{t\ra -1}(t^2-t^{-2})\d(p,q,n+j)=
\lim_{t\ra -1}(t^2-t^{-2})\d(p,q,n)$ for any fixed positive  integer $j$)
\begin{equation}\label{equ:g(-1,m)}
\lim_{t\ra -1}(t^2-t^{-2})g(t, \sm)=
\frac{\left(\sm^{a+b}+\sm^{-a-b}-\sm^{b-a}-\sm^{a-b}\right)
(\sm^{pq}-\sm^{-pq-2ab})}{\sm^{pq-3ab}-\sm^{-pq-ab}}\ne 0
\end{equation}
Therefore,
\begin{align*}
\a(t,\sm,\sl)=&(\sl-1)\frac{1}{(t^2-t^{-2})g(t, \sm)}
\left(\sl^2-\sm^{-2ab}t^{-4ab}\right)f^{-1}(t,\sm)
\left(\sl^2-\sm^{-2pq}t^{-4pq}\right)\\&\frac{\sm^{2pq+ab}t^{4pq+2ab}}{\d(a,b,n)}
\left([n+2]\sl^2- \sm^{-2pq-2ab}t^{-4pq-4ab}[n]\right)
\end{align*}
is an annihilator of $J_{\pq\#\ab}(n)$.
Also, using (\ref{equ:g(-1,m)}), (\ref{equ:f(-1,m)}), (\ref{equ:delta}), (\ref{[n]}), and Lemma \ref{lem:A-poly of sum} (i)  we have
\[\a(-1, \sm,\sl)=a(\sm)
(\sl-1)(\sl^2-\sm^{-2ab})(\sl^2-\sm^{-2pq})
(\sl^2- \sm^{-2pq-2ab})
=a(\sm)\sm^{-4pq-4ab}A_{\pq\#\ab}(\sm,\sl)\]
where
\[ a(\sm)=
\frac{\sm-\sm^{-1}}
{\left(\sm^{-pq-ab}-\sm^{-3pq-3ab}\right)\left
(\sm^{p+q}+\sm^{-p-q}-\sm^{q-p}-\sm^{p-q}\right)
\left(\sm^{a+b}+\sm^{-a-b}-\sm^{b-a}-\sm^{a-b}\right)}\]
is a well-defined, non-zero rational function
of $\sm$.

We now proceed to show that $\a(t,\sm,\sl)$ is the recurrence polynomial of $J_{\pq\#\ab}(n)$ up to a normalization.
Since $\a(t, \sm,\sl)$ has  $\sl$-degree $7$, we just need to show that $7$ is the minimum $\sl$-degree
among all non-zero annihilators of  $J_{\pq\#\ab}(n)$.
It is enough to show that if  $P = D_6\sl^6+D_5\sl^5+D_4 \sl^4 + D_3 \sl^3 + D_2 \sl^2 + D_1 \sl + D_0$
with $D_0,...,D_6\in \q(t,\sm)$ is an annihilator of $J_{\pq\#\ab}(n)$, then $P = 0$.

Suppose $P J_{\pq\#\ab}(n) = 0$, that is,
\begin{align*}
&D_6J_{\pq\#\ab}(n+6)+D_5J_{\pq\#\ab}(n+5)+D_4 J_{\pq\#\ab}(n+4)
 + D_3J_{\pq\#\ab}(n+3)\\& + D_2 J_{\pq\#\ab}(n+2)
+ D_1 J_{\pq\#\ab}(n+1)
 + D_0 J_{\pq\#\ab}(n) = 0.
\end{align*}
We wish to show that $D_i = 0$ for all $i = 0,1,...,6$.

Applying (\ref{relation for pqab}) and (\ref{relation of J_T(n+2)})
repeatedly, we have
\begin{align*}&D_6J_{\pq\#\ab}(n+6)+D_5J_{\pq\#\ab}(n+5)+D_4 J_{\pq\#\ab}(n+4)
 + D_3J_{\pq\#\ab}(n+3)\\ &+ D_2 J_{\pq\#\ab}(n+2) + D_1 J_{\pq\#\ab}(n+1)
 + D_0 J_{\pq\#\ab}(n) \\
&=\frac{D_6}{[n+6]}\left( t^{-4(pq+ab)(n+5)}[n+4]
J_{\pq\#\ab}(n+4)+t^{-(4pq+2ab)(n+5)}\delta(a,b,n+4)J_{\pq}(n+4)\right.&\\
&\;\;\;\;\left.+t^{-(2pq+4ab)(n+5)}\delta(p,q,n+4)J_{\ab}(n+4)+t^{-2(pq+ab)(n+5)}\delta(p,q,n+4)\delta(a,b,n+4)\right)\\&\;\;
+\frac{D_5}{[n+5]}
\left( t^{-4(pq+ab)(n+4)}[n+3]J_{\pq\#\ab}(n+3)+t^{-(4pq+2ab)(n+4)}\delta(a,b,n+3)J_{\pq}(n+3)\right.&\\
&\;\;\;\;\left.+t^{-(2pq+4ab)(n+4)}\delta(p,q,n+3)
J_{\ab}(n+3)+t^{-2(pq+ab)(n+4)}\delta(p,q,n+3)\delta(a,b,n+3)\right)
\\&\;\;+D_4 J_{\pq\#\ab}(n+4)
 + D_3J_{\pq\#\ab}(n+3)+ D_2 J_{\pq\#\ab}(n+2) \\ &\;\;+ D_1 J_{\pq\#\ab}(n+1)
 + D_0 J_{\pq\#\ab}(n)
\\
&
=\left(\frac{D_6}{[n+6]} t^{-4(pq+ab)(n+5)}[n+4]+D_4\right)J_{\pq\#\ab}(n+4)
\\
&\;\;+\left(\frac{D_5}{[n+5]} t^{-4(pq+ab)(n+4)}[n+3]+D_3\right)J_{\pq\#\ab}(n+3)
\\&\;\;
+\frac{D_6}{[n+6]}t^{-(4pq+2ab)(n+5)}\delta(a,b,n+4)J_{\pq}(n+4)
+\frac{D_6}{[n+6]}t^{-(2pq+4ab)(n+5)}\delta(p,q,n+4)J_{\ab}(n+4)
\\&\;\;
+\frac{D_6}{[n+6]}t^{-2(pq+ab)(n+5)}\delta(p,q,n+4)\delta(a,b,n+4)+\frac{D_5}{[n+5]}t^{-(4pq+2ab)(n+4)}\delta(a,b,n+3)J_{\pq}(n+3)
\\&\;\;+\frac{D_5}{[n+5]}t^{-(2pq+4ab)(n+4)}\delta(p,q,n+3)J_{\ab}(n+3)
+\frac{D_5}{[n+5]}t^{-2(pq+ab)(n+4)}\delta(p,q,n+3)\delta(a,b,n+3)
\\ &\;\;+ D_2 J_{\pq\#\ab}(n+2) + D_1 J_{\pq\#\ab}(n+1)
 + D_0 J_{\pq\#\ab}(n)
\\&
=\left(\frac{D_6}{[n+6]} t^{-4(pq+ab)(n+5)}+\frac{D_4}{[n+4]}\right)
\left( t^{-4(pq+ab)(n+3)}[n+2]J_{\pq\#\ab}(n+2)\right.\\&\;\;\;\;+t^{-(4pq+2ab)(n+3)}\delta(a,b,n+2)J_{\pq}(n+2)+t^{-(2pq+4ab)(n+3)}\delta(p,q,n+2)J_{\ab}(n+2)\\&\;\;\;\;\left.+t^{-2(pq+ab)(n+3)}\delta(p,q,n+2)\delta(a,b,n+2)\right)
\\&\;\;+\left(\frac{D_5}{[n+5]} t^{-4(pq+ab)(n+4)}+\frac{D_3}{[n+3]}\right)\left( t^{-4(pq+ab)(n+2)}[n+1]J_{\pq\#\ab}(n+1)\right.\\&\;\;\;\;+t^{-(4pq+2ab)(n+2)}\delta(a,b,n+1)J_{\pq}(n+1)+t^{-(2pq+4ab)(n+2)}\delta(p,q,n+1)J_{\ab}(n+1)\\&\;\;\;\;\left.+t^{-2(pq+ab)(n+2)}\delta(p,q,n+1)\delta(a,b,n+1)\right)
\\&\;\;
+\frac{D_6}{[n+6]}t^{-(4pq+2ab)(n+5)}\delta(a,b,n+4)
\left( t^{-4pq(n+3)} J_{\pq}(n+2)+t^{-2pq(n+3)} \delta(p,q,n+2)\right)
\\&\;\;
+\frac{D_6}{[n+6]}t^{-(2pq+4ab)(n+5)}\delta(p,q,n+4)
\left( t^{-4ab(n+3)} J_{\ab}(n+2)+t^{-2ab(n+3)} \delta(a,b,n+2)\right)
\\&\;\;
+\frac{D_6}{[n+6]}t^{-2(pq+ab)(n+5)}\delta(p,q,n+4)\delta(a,b,n+4)\\&\;\;
+\frac{D_5}{[n+5]}t^{-(4pq+2ab)(n+4)}\delta(a,b,n+3)
\left( t^{-4pq(n+2)} J_{\pq}(n+1)+t^{-2pq(n+2)} \delta(p,q,n+1)\right)\\&\;\;
+\frac{D_5}{[n+5]}t^{-(2pq+4ab)(n+4)}\delta(p,q,n+3)
\left( t^{-4ab(n+2)} J_{\ab}(n+1)+t^{-2ab(n+2)} \delta(a,b,n+1)\right)\\&\;\;
+\frac{D_5}{[n+5]}t^{-2(pq+ab)(n+4)}\delta(p,q,n+3)\delta(a,b,n+3)
\\&\;\;+ D_2 J_{\pq\#\ab}(n+2) + D_1 J_{\pq\#\ab}(n+1)
 + D_0 J_{\pq\#\ab}(n)
\\&=\left[\left(\frac{D_6}{[n+6]} t^{-4(pq+ab)(n+5)}+\frac{D_4}{[n+4]}\right)
 t^{-4(pq+ab)(n+3)}[n+2]+D_2\right]J_{\pq\#\ab}(n+2)
\\&\;\;
+\left[\left(\frac{D_6}{[n+6]} t^{-4(pq+ab)(n+5)}+\frac{D_4}{[n+4]}\right)
 t^{-(4pq+2ab)(n+3)}\delta(a,b,n+2)\right.\\&\;\;
\;\;\;\;\;\;+\left.\frac{D_6}{[n+6]}t^{-(4pq+2ab)(n+5)}\delta(a,b,n+4)
 t^{-4pq(n+3)} \right]J_{\pq}(n+2)
\\&\;\;
+\left[\left(\frac{D_6}{[n+6]} t^{-4(pq+ab)(n+5)}+\frac{D_4}{[n+4]}\right)\right.t^{-(2pq+4ab)(n+3)}\delta(p,q,n+2)
\\&\;\;\;\;\;\;\;\;
+\left.\frac{D_6}{[n+6]}t^{-(2pq+4ab)(n+5)}\delta(p,q,n+4)
t^{-4ab(n+3)} \right]J_{\ab}(n+2)
\\&\;\;+\left(\frac{D_6}{[n+6]} t^{-4(pq+ab)(n+5)}+\frac{D_4}{[n+4]}\right)t^{-2(pq+ab)(n+3)}\delta(p,q,n+2)\delta(a,b,n+2)
\\&\;\;
+\left[\left(\frac{D_5}{[n+5]} t^{-4(pq+ab)(n+4)}+\frac{D_3}{[n+3]}\right)
 t^{-4(pq+ab)(n+2)}[n+1]+D_1\right]J_{\pq\#\ab}(n+1)
\\&\;\;
+\left[\left(\frac{D_5}{[n+5]} t^{-4(pq+ab)(n+4)}+\frac{D_3}{[n+3]}\right)
t^{-(4pq+2ab)(n+2)}\delta(a,b,n+1)\right.
\\&\;\;\;\;\;\;\;\;
\left.+\frac{D_5}{[n+5]}t^{-(4pq+2ab)(n+4)}\delta(a,b,n+3)
 t^{-4pq(n+2)}\right]J_{\pq}(n+1)
\\&\;\;
+\left[\left(\frac{D_5}{[n+5]} t^{-4(pq+ab)(n+4)}+\frac{D_3}{[n+3]}\right)
t^{-(2pq+4ab)(n+2)}\delta(p,q,n+1)\right.
\\&\;\;\;\;\;\;\;\;
\left.+\frac{D_5}{[n+5]}t^{-(2pq+4ab)(n+4)}\delta(p,q,n+3)
 t^{-4ab(n+2)}\right]J_{\ab}(n+1)
\\&\;\;
+\left(\frac{D_5}{[n+5]} t^{-4(pq+ab)(n+4)}+\frac{D_3}{[n+3]}\right)
t^{-2(pq+ab)(n+2)}\delta(p,q,n+1)\delta(a,b,n+1)
\\&\;\;
+\frac{D_6}{[n+6]}t^{-(4pq+2ab)(n+5)}\delta(a,b,n+4)
 t^{-2pq(n+3)} \delta(p,q,n+2)
\\&\;\;
+\frac{D_6}{[n+6]}t^{-(2pq+4ab)(n+5)}\delta(p,q,n+4)
 t^{-2ab(n+3)} \delta(a,b,n+2)
\\&\;\;
+\frac{D_6}{[n+6]}t^{-2(pq+ab)(n+5)}\delta(p,q,n+4)\delta(a,b,n+4)\\&\;\;
+\frac{D_5}{[n+5]}t^{-(4pq+2ab)(n+4)}\delta(a,b,n+3)
t^{-2pq(n+2)} \delta(p,q,n+1)
\\&\;\;
+\frac{D_5}{[n+5]}t^{-(2pq+4ab)(n+4)}\delta(p,q,n+3)
t^{-2ab(n+2)} \delta(a,b,n+1)
\\&\;\;+\frac{D_5}{[n+5]}t^{-2(pq+ab)(n+4)}\delta(p,q,n+3)\delta(a,b,n+3)
 + D_0 J_{\pq\#\ab}(n)
\\
&
=\left[\left(\frac{D_6}{[n+6]} t^{-4(pq+ab)(n+5)}+\frac{D_4}{[n+4]}\right)
 t^{-4(pq+ab)(n+3)}+\frac{D_2}{[n+2]}\right]
\\&\;\;\;\;
 \left(t^{-4(pq+ab)(n+1)}[n]J_{\pq\#\ab}(n)+t^{-(4pq+2ab)(n+1)}\delta(a,b,n)J_{\pq}(n)\right.\\&
\;\;\;\;\;\;\;\left.+t^{-(2pq+4ab)(n+1)}\delta(p,q,n)J_{\ab}(n)
+t^{-2(pq+ab)(n+1)}\delta(p,q,n)\delta(a,b,n)\right)
\\&\;\;
+\left[\left(\frac{D_6}{[n+6]} t^{-4(pq+ab)(n+5)}+\frac{D_4}{[n+4]}\right)
 t^{-(4pq+2ab)(n+3)}\delta(a,b,n+2)\right.\\&\;\;
\;\;\;\;\;\;+\left.\frac{D_6}{[n+6]}t^{-(4pq+2ab)(n+5)}\delta(a,b,n+4)
 t^{-4pq(n+3)} \right]
\left(t^{-4pq(n+1)}J_{\pq}(n)+t^{-2pq(n+1)}\d(p,q,n)\right)
\\&\;\;
+\left[\left(\frac{D_6}{[n+6]} t^{-4(pq+ab)(n+5)}+\frac{D_4}{[n+4]}\right)\right.t^{-(2pq+4ab)(n+3)}\delta(p,q,n+2)
\\&\;\;\;\;\;\;\;\;
+\left.\frac{D_6}{[n+6]}t^{-(2pq+4ab)(n+5)}\delta(p,q,n+4)
t^{-4ab(n+3)} \right]\left(t^{-4ab(n+1)}J_{\ab}(n)
+t^{-2ab(n+1)}\d(a,b,n)\right)\\
&\;\;+\left(\frac{D_6}{[n+6]} t^{-4(pq+ab)(n+5)}+\frac{D_4}{[n+4]}\right)t^{-2(pq+ab)(n+3)}\delta(p,q,n+2)\delta(a,b,n+2)
\\&\;\;
+\left[\left(\frac{D_5}{[n+5]} t^{-4(pq+ab)(n+4)}+\frac{D_3}{[n+3]}\right)
 t^{-4(pq+ab)(n+2)}[n+1]+D_1\right]J_{\pq\#\ab}(n+1)
\\&\;\;
+\left[\left(\frac{D_5}{[n+5]} t^{-4(pq+ab)(n+4)}+\frac{D_3}{[n+3]}\right)
t^{-(4pq+2ab)(n+2)}\delta(a,b,n+1)\right.
\\&\;\;\;\;\;\;\;\;
\left.+\frac{D_5}{[n+5]}t^{-(4pq+2ab)(n+4)}\delta(a,b,n+3)
 t^{-4pq(n+2)}\right]J_{\pq}(n+1)
\\&\;\;
+\left[\left(\frac{D_5}{[n+5]} t^{-4(pq+ab)(n+4)}+\frac{D_3}{[n+3]}\right)
t^{-(2pq+4ab)(n+2)}\delta(p,q,n+1)\right.
\\&\;\;\;\;\;\;\;\;
\left.+\frac{D_5}{[n+5]}t^{-(2pq+4ab)(n+4)}\delta(p,q,n+3)
 t^{-4ab(n+2)}\right]J_{\ab}(n+1)
\\&\;\;
+\left(\frac{D_5}{[n+5]} t^{-4(pq+ab)(n+4)}+\frac{D_3}{[n+3]}\right)
t^{-2(pq+ab)(n+2)}\delta(p,q,n+1)\delta(a,b,n+1)
\\&\;\;
+\frac{D_6}{[n+6]}t^{-(4pq+2ab)(n+5)}\delta(a,b,n+4)
 t^{-2pq(n+3)} \delta(p,q,n+2)
\\&\;\;
+\frac{D_6}{[n+6]}t^{-(2pq+4ab)(n+5)}\delta(p,q,n+4)
 t^{-2ab(n+3)} \delta(a,b,n+2)
\\&\;\;
+\frac{D_6}{[n+6]}t^{-2(pq+ab)(n+5)}\delta(p,q,n+4)\delta(a,b,n+4)\\&\;\;
+\frac{D_5}{[n+5]}t^{-(4pq+2ab)(n+4)}\delta(a,b,n+3)
t^{-2pq(n+2)} \delta(p,q,n+1)
\\&\;\;
+\frac{D_5}{[n+5]}t^{-(2pq+4ab)(n+4)}\delta(p,q,n+3)
t^{-2ab(n+2)} \delta(a,b,n+1)
\\&\;\;+\frac{D_5}{[n+5]}t^{-2(pq+ab)(n+4)}\delta(p,q,n+3)\delta(a,b,n+3)
 + D_0 J_{\pq\#\ab}(n)
\\
&=\left[\left(\frac{D_5}{[n+5]} t^{-4(pq+ab)(n+4)}+\frac{D_3}{[n+3]}\right)
 t^{-4(pq+ab)(n+2)}[n+1]+D_1\right]J_{\pq\#\ab}(n+1)
\\&\;\;
+\left\{\left[\left(\frac{D_6}{[n+6]} t^{-4(pq+ab)(n+5)}+\frac{D_4}{[n+4]}\right)
 t^{-4(pq+ab)(n+3)}+\frac{D_2}{[n+2]}\right]t^{-4(pq+ab)(n+1)}[n]+D_0\right\}\\&\;\;\;\;J_{\pq\#\ab}(n)
+\left[\left(\frac{D_5}{[n+5]} t^{-4(pq+ab)(n+4)}+\frac{D_3}{[n+3]}\right)
t^{-(4pq+2ab)(n+2)}\delta(a,b,n+1)\right.
\\&\;\;\;\;\;\;\;\;
\left.+\frac{D_5}{[n+5]}t^{-(4pq+2ab)(n+4)}\delta(a,b,n+3)
 t^{-4pq(n+2)}\right]J_{\pq}(n+1)
\\&\;\;
+\left\{\left[\left(\frac{D_6}{[n+6]} t^{-4(pq+ab)(n+5)}+\frac{D_4}{[n+4]}\right)
 t^{-4(pq+ab)(n+3)}+\frac{D_2}{[n+2]}\right]\right.
\\&\;\;\;\;
t^{-(4pq+2ab)(n+1)}\delta(a,b,n)
+\left[\left(\frac{D_6}{[n+6]} t^{-4(pq+ab)(n+5)}+\frac{D_4}{[n+4]}\right)
 t^{-(4pq+2ab)(n+3)}\delta(a,b,n+2)\right.\\&\;\;
\;\;\;\;\;\;+\left.
\left.\frac{D_6}{[n+6]}t^{-(4pq+2ab)(n+5)}\delta(a,b,n+4)
 t^{-4pq(n+3)} \right]
t^{-4pq(n+1)}\right\}J_{\pq}(n)
\\&\;\;
+\left[\left(\frac{D_5}{[n+5]} t^{-4(pq+ab)(n+4)}+\frac{D_3}{[n+3]}\right)
t^{-(2pq+4ab)(n+2)}\delta(p,q,n+1)\right.
\\&\;\;\;\;\;\;\;\;
\left.+\frac{D_5}{[n+5]}t^{-(2pq+4ab)(n+4)}\delta(p,q,n+3)
 t^{-4ab(n+2)}\right]J_{\ab}(n+1)
\\&\;\;
+\left\{\left[\left(\frac{D_6}{[n+6]} t^{-4(pq+ab)(n+5)}+\frac{D_4}{[n+4]}\right)
 t^{-4(pq+ab)(n+3)}+\frac{D_2}{[n+2]}\right]\right.
\\&\;\;\;\;
 t^{-(2pq+4ab)(n+1)}\delta(p,q,n)
+\left[\left(\frac{D_6}{[n+6]} t^{-4(pq+ab)(n+5)}+\frac{D_4}{[n+4]}\right)\right.t^{-(2pq+4ab)(n+3)}\delta(p,q,n+2)
\\&\;\;\;\;\;\;\;\;
+\left.\left.\frac{D_6}{[n+6]}t^{-(2pq+4ab)(n+5)}\delta(p,q,n+4)
t^{-4ab(n+3)} \right]t^{-4ab(n+1)}\right\}J_{\ab}(n)
\\&\;\;
+\left[\left(\frac{D_6}{[n+6]} t^{-4(pq+ab)(n+5)}+\frac{D_4}{[n+4]}\right)
 t^{-4(pq+ab)(n+3)}+\frac{D_2}{[n+2]}\right]
t^{-2(pq+ab)(n+1)}\delta(p,q,n)\delta(a,b,n)
\\&\;\;
+\left[\left(\frac{D_6}{[n+6]} t^{-4(pq+ab)(n+5)}+\frac{D_4}{[n+4]}\right)
 t^{-(4pq+2ab)(n+3)}\delta(a,b,n+2)\right.\\&\;\;
\;\;\;\;\;\;+\left.\frac{D_6}{[n+6]}t^{-(4pq+2ab)(n+5)}\delta(a,b,n+4)
 t^{-4pq(n+3)} \right]
t^{-2pq(n+1)}\d(p,q,n)
\\&\;\;
+\left[\left(\frac{D_6}{[n+6]} t^{-4(pq+ab)(n+5)}+\frac{D_4}{[n+4]}\right)\right.t^{-(2pq+4ab)(n+3)}\delta(p,q,n+2)
\\&\;\;\;\;\;\;\;\;
+\left.\frac{D_6}{[n+6]}t^{-(2pq+4ab)(n+5)}\delta(p,q,n+4)
t^{-4ab(n+3)} \right]
t^{-2ab(n+1)}\d(a,b,n)\\
&\;\;+\left(\frac{D_6}{[n+6]} t^{-4(pq+ab)(n+5)}+\frac{D_4}{[n+4]}\right)t^{-2(pq+ab)(n+3)}\delta(p,q,n+2)\delta(a,b,n+2)
\\&\;\;
+\left(\frac{D_5}{[n+5]} t^{-4(pq+ab)(n+4)}+\frac{D_3}{[n+3]}\right)
t^{-2(pq+ab)(n+2)}\delta(p,q,n+1)\delta(a,b,n+1)
\\&\;\;
+\frac{D_6}{[n+6]}t^{-(4pq+2ab)(n+5)}\delta(a,b,n+4)
 t^{-2pq(n+3)} \delta(p,q,n+2)
\\&\;\;
+\frac{D_6}{[n+6]}t^{-(2pq+4ab)(n+5)}\delta(p,q,n+4)
 t^{-2ab(n+3)} \delta(a,b,n+2)
\\&\;\;
+\frac{D_6}{[n+6]}t^{-2(pq+ab)(n+5)}\delta(p,q,n+4)\delta(a,b,n+4)\\&\;\;
+\frac{D_5}{[n+5]}t^{-(4pq+2ab)(n+4)}\delta(a,b,n+3)
t^{-2pq(n+2)} \delta(p,q,n+1)
\\&\;\;
+\frac{D_5}{[n+5]}t^{-(2pq+4ab)(n+4)}\delta(p,q,n+3)
t^{-2ab(n+2)} \delta(a,b,n+1)
\\&\;\;+\frac{D_5}{[n+5]}t^{-2(pq+ab)(n+4)}\delta(p,q,n+3)
\delta(a,b,n+3)
\\&=E_6J_{\pq\#\ab}(n+1)+E_5J_{\pq\#\ab}(n)
+E_4J_{T(p,q)}(n+1)+E_3J_{T(p,q)}(n)\\&\;\;+E_2J_{T(a,b)}(n+1)
+E_1J_{T(a,b)}(n)+E_0.\end{align*}

\begin{lemma}\label{lem:Eizero}For any $E_i\in\q(t,\sm)$, $i=0,1, \dots,6$, the equation
\begin{equation}\label{Ei are zero}
\begin{aligned}&E_6J_{\pq\#\ab}(n+1)+E_5J_{\pq\#\ab}(n)
+E_4J_{T(p,q)}(n+1)+E_3J_{T(p,q)}(n)\\&+E_2J_{T(a,b)}(n+1)
+E_1J_{T(a,b)}(n)+E_0=0\end{aligned}
\end{equation}implies
$E_i=0$, for all $i=0,1, \dots, 6$.
\end{lemma}

\pf  Up to multiplying equation (\ref{Ei are zero}) by a non-zero rational function
in $\q(t,\sm)$ (which clears the denominators of all $E_i$), we may and shall assume that
each $E_i$ is in $\z[t,\sm]$.

We first consider the subcase that  both $p$ and $a$ are positive.
From equation (\ref{Ei are zero}), we see that
at least two of the summands have the same lowest degree in $t$
(here we treat $\sm$ as $t^{2n}$). That is, at least one of
the following equalities holds:
\begin{align*}
&\ell\left(E_6J_{\pq\#\ab}(n+1)\right)=\ell\left(E_5J_{\pq\#\ab}(n)\right), \quad
\ell\left(E_6J_{\pq\#\ab}(n+1)\right)=\ell\left(E_4J_{\pq}(n+1)\right),\\
&\ell\left(E_6J_{\pq\#\ab}(n+1)\right)=\ell\left(E_3J_{\pq}(n)\right),
\quad\ell\left(E_6J_{\pq\#\ab}(n+1)\right)=\ell\left(E_2J_{\ab}(n+1)\right), \\
&\ell\left(E_6J_{\pq\#\ab}(n+1)\right)=\ell\left(E_1J_{\ab}(n)\right),
\quad\ell\left(E_6J_{\pq\#\ab}(n+1)\right)=\ell(E_0),\\
&\ell\left(E_5J_{\pq\#\ab}(n)\right)=\ell\left(E_4J_{\pq}(n+1)\right),
\quad\ell\left(E_5J_{\pq\#\ab}(n)\right)=\ell\left(E_3J_{\pq}(n)\right),\\&
\ell\left(E_5J_{\pq\#\ab}(n)\right)=\ell\left(E_2J_{\ab}(n+1)\right),
\quad\ell\left(E_5J_{\pq\#\ab}(n)\right)=\ell\left(E_1J_{\ab}(n)\right),\\&
\ell\left(E_5J_{\pq\#\ab}(n)\right)=\ell(E_0),
\quad\ell\left(E_4J_{\pq}(n+1)\right)=\ell\left(E_3J_{\pq}(n)\right),\\&
\ell\left(E_4J_{\pq}(n+1)\right)=\ell\left(E_2J_{\ab}(n+1)\right),
\quad\ell\left(E_4J_{\pq}(n+1)\right)=\ell\left(E_1J_{\ab}(n)\right),\\&
\ell\left(E_4J_{\pq}(n+1)\right)=\ell(E_0),
\quad\ell\left(E_3J_{\pq}(n)\right)=\ell\left(E_2J_{\ab}(n+1)\right),\\&
\ell\left(E_3J_{\pq}(n)\right)=\ell\left(E_1J_{\ab}(n)\right),
\quad\ell\left(E_3J_{\pq}(n)\right)=\ell(E_0),
\quad\ell\left(E_2J_{\ab}(n+1)\right)=\ell\left(E_1J_{\ab}(n)\right), \\&
\ell\left(E_2J_{\ab}(n+1)\right)=\ell(E_0),
\quad\ell\left(E_1J_{\ab}(n)\right)=\ell(E_0).
\end{align*}

If $\displaystyle\ell\left(E_6J_{\pq\#\ab}(n+1)\right)=\ell\left(E_5J_{\pq\#\ab}(n)\right)$,
then either both $E_6$ and $E_5$ are zero or both are non-zero. In the latter case, we have
$\displaystyle\ell(E_6)+\ell\left(J_{\pq\#\ab}(n+1)\right)=\ell(E_5)+\ell\left(J_{\pq\#\ab}(n)\right)$.
Applying  Lemma \ref{degree of connected sum} (1),
 we see that
\begin{align*}
\ell(E_6)-\ell(E_5)&=\ell\left(J_{\pq\#\ab}(n)\right)-\ell\left(J_{\pq\#\ab}(n+1)\right)\\
&=\left(-pqn^2+pq+\frac12(1-(-1)^{n-1})(p-2)(q-2)
-abn^2+ab\right.\\&\;\;\;\;\;\;\left.+\frac12(1-(-1)^{n-1})(a-2)(b-2)-2n-2\right)\\&\;\;
-\left(-pq(n+1)^2+pq+\frac12(1-(-1)^{n})(p-2)(q-2)
-ab(n+1)^2+ab\right.\\&\;\;\;\;\;\;\;\;
\left.+\frac12(1-(-1)^{n})(a-2)(b-2)-2(n+1)-2\right)
\\&=
2pqn+(-1)^n[(p-2)(q-2)+(a-2)(b-2)]+pq+2abn+ab+2.\end{align*}
However,  the far left-hand side of these connecting equations  is a linear function of $n$, while the far right-hand side is not a linear function of $n$ (as $(p-2)(q-2)+(a-2)(b-2)>0$), which is a contradiction.

If $\displaystyle\ell\left(E_6J_{\pq\#\ab}(n+1)\right)=\ell\left(E_4J_{T(p,q)}(n+1)\right)$,
then either both $E_6$ and $E_4$ are zero or both are non-zero. In the latter case,  by applying Lemmas \ref{degree of J_T} (1) and   \ref{degree of connected sum} (1) we have
\begin{align*}
\ell(E_6)-\ell(E_4)&=\ell\left(J_{T(p,q)}(n+1)\right)-\ell\left(J_{\pq\#\ab}(n+1)\right)\\
&=\left(-pq(n+1)^2+pq+\frac12(1-(-1)^{n})(p-2)(q-2)
\right)\\&\;\;
-\left(-pq(n+1)^2+pq+\frac12(1-(-1)^{n})(p-2)(q-2)
-ab(n+1)^2+ab\right.\\&\;\;\left.+\frac12(1-(-1)^{n})(a-2)(b-2)-2(n+1)-2\right)
\\&=ab(n+1)^2-ab-\frac12(1-(-1)^{n})(a-2)(b-2)+2(n+1)+2.
\end{align*}
Again, the far left-hand side of these connecting equations is a linear function of $n$ while the far right-hand side is not, yielding a contradiction.

Similarly, whenever two summands of (\ref{Ei are zero}) have the same lowest degree in $t$,
we can get a contradiction using the degree formulas given in Section \ref{sec:formulas},
unless both  involved $E_i$ and $E_j$
are zero.

It follows that $E_i=0$ for all $i=0,1, \dots,6$.

The subcase where both $p$ and $a$ are negative can be treated similarly, but using
the highest degree in $t$ instead of the lowest.
\qed

By Lemma \ref{lem:Eizero} and the calculation preceding it, we
obtain a system of seven  homogeneous  linear equations, $E_i=0$ for
$i \in \{0,1,\dots,6\}$,
in the  variables $D_0, \dots, D_6$ with coefficients in $\q(t,\sm)$.
We now demonstrate that this system  admits only the trivial solution:  $D_i=0$ for all $i$.

Our approach involves expressing all coefficients
as functions of $t$ and $\sm$ and  evaluating them at $t=-1$,
 followed by simplifications such as dividing equations
by non-zero functions of $\sm$.
Notably, under this reduction, the terms
 $[n+i]$ and $[n+j]$  coincide, as do  $\d(p,q,n+i)$ and $\d(p,q,n+j)$.
We then verify  that the resulting simplified linear system
has only the trivial solution.
Below, we use $\lra$ to denote
the reduction of the equation $E_i=0$.

\begin{align*}& E_6=\left(\frac{D_5}{[n+5]} t^{-4(pq+ab)(n+4)}+\frac{D_3}{[n+3]}\right)
 t^{-4(pq+ab)(n+2)}[n+1]+D_1=0
\\&\lra\;\; \sm^{-4(pq+ab)}D_5+\sm^{-2(pq+ab)}D_3+D_1=0,
\\&E_5=\left[\left(\frac{D_6}{[n+6]} t^{-4(pq+ab)(n+5)}+\frac{D_4}{[n+4]}\right)
 t^{-4(pq+ab)(n+3)}+\frac{D_2}{[n+2]}\right]t^{-4(pq+ab)(n+1)}[n]+D_0=0
\\& \lra\;\;
\sm^{-6(pq+ab)}D_6+\sm^{-4(pq+ab)}D_4+\sm^{-2(pq+ab)}D_2+D_0=0,
\\& E_4=\left(\frac{D_5}{[n+5]} t^{-4(pq+ab)(n+4)}+\frac{D_3}{[n+3]}\right)
t^{-(4pq+2ab)(n+2)}\delta(a,b,n+1)
\\&\;\;\;\;\;\;\;\;\;\; +\frac{D_5}{[n+5]}t^{-(4pq+2ab)(n+4)}\delta(a,b,n+3)
 t^{-4pq(n+2)}=0\\
& \lra\;\;
(\sm^{-2(pq+ab)}+\sm^{-2pq})D_5+D_3=0,
\\& E_3=\left[\left(\frac{D_6}{[n+6]} t^{-4(pq+ab)(n+5)}+\frac{D_4}{[n+4]}\right)
 t^{-4(pq+ab)(n+3)}+\frac{D_2}{[n+2]}\right]
t^{-(4pq+2ab)(n+1)}\delta(a,b,n)\\
& \;\;\;\;\;\;\;\;\;\; +\left[\left(\frac{D_6}{[n+6]} t^{-4(pq+ab)(n+5)}+\frac{D_4}{[n+4]}\right)
 t^{-(4pq+2ab)(n+3)}\delta(a,b,n+2)\right.\\
&\;\;\;\;\;\;\;\;\;\;\;\;\;\;  +
\left.\frac{D_6}{[n+6]}t^{-(4pq+2ab)(n+5)}\delta(a,b,n+4)
 t^{-4pq(n+3)} \right]
t^{-4pq(n+1)}=0\\& \lra\;\;
(\sm^{-4(pq+ab)}+\sm^{-4pq-2ab}+\sm^{-4pq})D_6
+(\sm^{-2(pq+ab)}+\sm^{-2pq})D_4+D_2=0,
\\&
E_2=\left(\frac{D_5}{[n+5]} t^{-4(pq+ab)(n+4)}+\frac{D_3}{[n+3]}\right)
t^{-(2pq+4ab)(n+2)}\delta(p,q,n+1)
\\& \;\;\;\;\;\;\;\;\;\; +\frac{D_5}{[n+5]}t^{-(2pq+4ab)(n+4)}\delta(p,q,n+3)
 t^{-4ab(n+2)}=0\\
& \lra\;\;
(\sm^{-2(pq+ab)}+\sm^{-2ab})D_5+D_3=0,
\\&E_1=\left[\left(\frac{D_6}{[n+6]} t^{-4(pq+ab)(n+5)}+\frac{D_4}{[n+4]}\right)
 t^{-4(pq+ab)(n+3)}+\frac{D_2}{[n+2]}\right]
 t^{-(2pq+4ab)(n+1)}\delta(p,q,n)\\
& \;\;\;\;\;\;\;\;\;\;  +\left[\left(\frac{D_6}{[n+6]} t^{-4(pq+ab)(n+5)}+\frac{D_4}{[n+4]}\right)\right.t^{-(2pq+4ab)(n+3)}\delta(p,q,n+2)
\\
& \;\;\;\;\;\;\;\;\;\;\;\; \;\;+\left.\frac{D_6}{[n+6]}t^{-(2pq+4ab)(n+5)}\delta(p,q,n+4)
t^{-4ab(n+3)} \right]t^{-4ab(n+1)}=0
\\& \lra\;\;
(\sm^{-4(pq+ab)}+\sm^{-2pq-4ab}+\sm^{-4ab})D_6
+(\sm^{-2(pq+ab)}+\sm^{-2ab})D_4+D_2=0,
\\
& E_0=
\left[\left(\frac{D_6}{[n+6]} t^{-4(pq+ab)(n+5)}+\frac{D_4}{[n+4]}\right)
 t^{-4(pq+ab)(n+3)}+\frac{D_2}{[n+2]}\right]
t^{-2(pq+ab)(n+1)}\delta(p,q,n)\delta(a,b,n)
\\&\;\;\;\;\;\;\;\;\;\;
\displaystyle+\left[\left(\frac{D_6}{[n+6]} t^{-4(pq+ab)(n+5)}+\frac{D_4}{[n+4]}\right)
 t^{-(4pq+2ab)(n+3)}\delta(a,b,n+2)\right.\\&\;\;
\;\;\;\;\;\;\;\;\;\;\;\;\displaystyle+\left.\frac{D_6}{[n+6]}t^{-(4pq+2ab)(n+5)}\delta(a,b,n+4)
 t^{-4pq(n+3)} \right]
t^{-2pq(n+1)}\d(p,q,n)
\\&\;\;\;\;\;\;\;\;\;\;
\displaystyle+\left[\left(\frac{D_6}{[n+6]} t^{-4(pq+ab)(n+5)}+\frac{D_4}{[n+4]}\right)\right.t^{-(2pq+4ab)(n+3)}\delta(p,q,n+2)
\\&\;\;\;\;\;\;\;\;\;\;\;\;\;\;
\displaystyle+\left.\frac{D_6}{[n+6]}t^{-(2pq+4ab)(n+5)}\delta(p,q,n+4)
t^{-4ab(n+3)} \right]
t^{-2ab(n+1)}\d(a,b,n)\\
&\;\;\;\;\;\;\;\;\;\;\displaystyle+\left(\frac{D_6}{[n+6]} t^{-4(pq+ab)(n+5)}+\frac{D_4}{[n+4]}\right)t^{-2(pq+ab)(n+3)}\delta(p,q,n+2)\delta(a,b,n+2)
\\&\;\;\;\;\;\;\;\;\;\;
\displaystyle+\left(\frac{D_5}{[n+5]} t^{-4(pq+ab)(n+4)}+\frac{D_3}{[n+3]}\right)
t^{-2(pq+ab)(n+2)}\delta(p,q,n+1)\delta(a,b,n+1)
\\&\;\;\;\;\;\;\;\;\;\;
\displaystyle+\frac{D_6}{[n+6]}t^{-(4pq+2ab)(n+5)}\delta(a,b,n+4)
 t^{-2pq(n+3)} \delta(p,q,n+2)
\\&\;\;\;\;\;\;\;\;\;\;
\displaystyle+\frac{D_6}{[n+6]}t^{-(2pq+4ab)(n+5)}\delta(p,q,n+4)
 t^{-2ab(n+3)} \delta(a,b,n+2)
\\&\;\;\;\;\;\;\;\;\;\;
\displaystyle+\frac{D_6}{[n+6]}t^{-2(pq+ab)(n+5)}\delta(p,q,n+4)\delta(a,b,n+4)\\&\;\;\;\;\;\;\;\;\;\;
\displaystyle+\frac{D_5}{[n+5]}t^{-(4pq+2ab)(n+4)}\delta(a,b,n+3)
t^{-2pq(n+2)} \delta(p,q,n+1)
\\&\;\;\;\;\;\;\;\;\;\;
\displaystyle+\frac{D_5}{[n+5]}t^{-(2pq+4ab)(n+4)}\delta(p,q,n+3)
t^{-2ab(n+2)} \delta(a,b,n+1)
\\&\;\;\;\;\;\;\;\;\;\;\displaystyle+\frac{D_5}{[n+5]}t^{-2(pq+ab)(n+4)}\delta(p,q,n+3)\delta(a,b,n+3)=0\\
&\lra\;\;
\displaystyle(\sm^{-4(pq+ab)}+\sm^{-4pq-2ab}+\sm^{-4pq}
+\sm^{-2pq-4ab}+\sm^{-4ab}
+\sm^{-2(pq+ab)}+\sm^{-2pq}
+\sm^{-2ab}+1)D_6
\\&\;\;\;\;\;\;\;\;\displaystyle
+(\sm^{-2(pq+ab)}+\sm^{-2pq}+\sm^{-2ab}
+1)D_5+(\sm^{-2(pq+ab)}+\sm^{-2pq}
+\sm^{-2ab}+1)D_4+D_3+D_2=0.
\end{align*}

That is, the given  system
has been reduced to the following one (with coefficients in  $\z[\sm^{\pm1}]$):
\begin{align*}& \sm^{-4(pq+ab)}D_5+\sm^{-2(pq+ab)}D_3+D_1=0,
\\&
\sm^{-6(pq+ab)}D_6+\sm^{-4(pq+ab)}D_4+\sm^{-2(pq+ab)}D_2+D_0=0,
\\&
(\sm^{-2(pq+ab)}+\sm^{-2pq})D_5+D_3=0,
\\&
(\sm^{-4(pq+ab)}+\sm^{-4pq-2ab}+\sm^{-4pq})D_6
+(\sm^{-2(pq+ab)}+\sm^{-2pq})D_4+D_2=0,
\\&
(\sm^{-2(pq+ab)}+\sm^{-2ab})D_5+D_3=0,
\\&
(\sm^{-4(pq+ab)}+\sm^{-2pq-4ab}+\sm^{-4ab})D_6
+(\sm^{-2(pq+ab)}+\sm^{-2ab})D_4+D_2=0,
\\&
(\sm^{-4(pq+ab)}+\sm^{-4pq-2ab}+\sm^{-4pq}
+\sm^{-2pq-4ab}+\sm^{-4ab}
+\sm^{-2(pq+ab)}+\sm^{-2pq}
+\sm^{-2ab}+1)D_6
\\&
+(\sm^{-2(pq+ab)}+\sm^{-2pq}+\sm^{-2ab}
+1)D_5+(\sm^{-2(pq+ab)}+\sm^{-2pq}
+\sm^{-2ab}+1)D_4+D_3+D_2=0.
\end{align*}

The third and fifth equations imply $D_5=D_3=0$ since $pq\ne ab$.
Then the first equation implies $D_1=0$.
The remaining  four equations become:
\begin{align*}&
\sm^{-6(pq+ab)}D_6+\sm^{-4(pq+ab)}D_4+\sm^{-2(pq+ab)}D_2+D_0=0,\\&
(\sm^{-4(pq+ab)}+\sm^{-4pq-2ab}+\sm^{-4pq})D_6
+(\sm^{-2(pq+ab)}+\sm^{-2pq})D_4+D_2=0,\\&
(\sm^{-4(pq+ab)}+\sm^{-2pq-4ab}+\sm^{-4ab})D_6
+(\sm^{-2(pq+ab)}+\sm^{-2ab})D_4+D_2=0,\\&
(\sm^{-4(pq+ab)}+\sm^{-4pq-2ab}+\sm^{-4pq}
+\sm^{-2pq-4ab}+\sm^{-4ab}
+\sm^{-2(pq+ab)}+\sm^{-2pq}
+\sm^{-2ab}+1)D_6,\\&
+(\sm^{-2(pq+ab)}+\sm^{-2pq}
+\sm^{-2ab}+1)D_4+D_2=0.
\end{align*}
The determinant of the coefficient matrix of the  four equations
is   $(\sm^{-2pq}-\sm^{-2ab})(\sm^{-2pq-2ab}-1)$
which is non-zero since $pq\ne ab$ and $pq\ne -ab$.
Hence $D_0=D_2=D_4=D_6=0$.

\section{Case $|p|>q>2$, $(a,b)=(p,q)$}\label{sec:(a,b)=(p,q)}

That is,  in this section, we verify the $AJ$-conjecture for $T(p,q)\# T(p,q)$ when $|p|>q>2$, which
 quickly follows from what has been done in Section \ref{sec:pqab}.
Equation  (\ref{3.2}) in Section \ref{sec:pqab} is still valid when $(p,q)=(a,b)$.
So, plugging $a=p$ and $b=q$ into the equation, we get
\begin{align*}&
\left(\sl^2-\sm^{-2pq}t^{-4pq}\right)\frac{\sm^{3pq}t^{6pq}}{\d(p,q,n)}\left([n+2]\sl^2-\sm^{-4pq}t^{-8pq}[n]\right)J_{\pq\#\pq}(n)
\\&
=\sm^{-pq}t^{-2pq}\d(p,q,n)+\sm^{pq}t^{6pq}\d(p,q,n+2).\end{align*}
The  function $f(t, \sm)=\sm^{-pq}t^{-2pq}\d(p,q,n)+\sm^{pq}t^{6pq}\d(p,q,n+2)$ is non-zero
since \[\lim_{t\ra -1}(t^2-t^{-2})f(t, \sm)=(\sm^{-pq}+\sm^{pq})(\sm^{p+q}+\sm^{-p-q}-\sm^{q-p}-\sm^{p-q})\ne 0.\]
Therefore,
\[\a(t, \sm,\sl)=(\sl-1)\frac{1}{(t^2-t^{-2})f(t,\sm)}(\sl^2-\sm^{-2pq}t^{-4pq})\frac{\sm^{3pq}t^{6pq}}{\d(p,q,n)}\left([n+2]\sl^2-\sm^{-4pq}t^{-8pq}[n]\right)\]
is an annihilator of  $J_{\pq\#\pq}(n)$ of $\sl$-degree $5$.
Also, \[\a(-1,\sm,\sl)=a(\sm)(\sl-1)(\sl^2-\sm^{-2pq})(\sl^2-\sm^{-4pq})=a(\sm)\sm^{-6pq}A_{T(p,q)\#T(p,q)}(\sm,\sl)\]
where \[a(\sm)=\frac{\sm^{3pq}(\sm-\sm^{-1})}{(\sm^{-pq}+\sm^{pq})(\sm^{p+q}+\sm^{-p-q}-\sm^{q-p}-\sm^{p-q})^2}\] is a well-defined, non-zero rational function of $\sm$.

We now show that $5$ is the minimum $\sl$-degree among all non-zero annihilators of
$J_{\pq\#\pq}(n)$. So suppose $D_4 \sl^4 + D_3 \sl^3 + D_2 \sl^2 + D_1 \sl + D_0$
with $D_0,\dots,D_4\in \q(t,\sm)$ is an annihilator of $J_{T(p,q)\#T(p,q)}(n)$,
that is,
\begin{align*}&D_4 J_{\pq\#\pq}(n+4)
 + D_3J_{\pq\#\pq}(n+3)+ D_2 J_{\pq\#\pq}(n+2)
\\&+ D_1 J_{\pq\#\pq}(n+1)
 + D_0 J_{\pq\#\pq}(n) = 0,\end{align*}
we want to show that $D_i = 0$ for all $i$.
Plugging $D_6=D_5=0$,  $a=p$, $b=q$ into the corresponding calculation carried
out in Section \ref{sec:pqab}, we get
\begin{align*}&D_4 J_{\pq\#\pq}(n+4)
 + D_3J_{\pq\#\pq}(n+3)+ D_2 J_{\pq\#\pq}(n+2) \\ &\;\;\;\;\;\;+ D_1 J_{\pq\#\pq}(n+1)
 + D_0 J_{\pq\#\pq}(n) \\&
=\left(\frac{D_3}{[n+3]}
 t^{-8pq(n+2)}[n+1]+D_1\right)J_{\pq\#\pq}(n+1)
\\&\;\;
+\left[\left(\frac{D_4}{[n+4]}
 t^{-8pq(n+3)}+\frac{D_2}{[n+2]}\right)t^{-8pq(n+1)}[n]+D_0\right]J_{\pq\#\pq}(n)
\\&\;\;+\frac{D_3}{[n+3]}2t^{-6pq(n+2)}\delta(p,q,n+1)J_{\pq}(n+1)
\\&\;\;
+\left[\left(\frac{D_4}{[n+4]}
 t^{-8pq(n+3)}+\frac{D_2}{[n+2]}\right)2t^{-6pq(n+1)}\delta(p,q,n)
+\frac{D_4}{[n+4]} 2t^{-6pq(n+3)}\delta(p,q,n+2)t^{-4pq(n+1)}\right]J_{\pq}(n)
\\&\;\;
+\left(\frac{D_4}{[n+4]} t^{-8pq(n+3)}+\frac{D_2}{[n+2]}\right)
t^{-4pq(n+1)}\d^2(p,q,n)
\\&\;\;
+\frac{D_4}{[n+4]}2 t^{-6pq(n+3)}\delta(p,q,n+2)  t^{-2pq(n+1)}\d(p,q,n)
\\&\;\;
+\frac{D_4}{[n+4]}t^{-4pq(n+3)}\d^2(p,q,n+2)
+\frac{D_3}{[n+3]}t^{-4pq(n+2)}\d^2(p,q,n+1)
\\&=E_4J_{\pq\#\pq}(n+1)+E_3J_{\pq\#\pq}(n)
+E_2J_{T(p,q)}(n+1)+E_1J_{T(p,q)}(n)+E_0\end{align*}

\begin{lemma}
$E_4J_{\pq\#\pq}(n+1)+E_3J_{\pq\#\pq}(n)
+E_2J_{T(p,q)}(n+1)+E_1J_{T(p,q)}(n)+E_0=0$ implies
$E_i=0$ for all $i=0,1,\dots,4$.
\end{lemma}

This lemma can be treated as a subcase of Lemma \ref{lem:Eizero}.

 Now we show that the system of  five  homogeneous  linear equations $E_i=0$, $i=0,1\dots,4$,
in five variables $D_0, \dots, D_4$ has only the trivial solution: $D_i=0$ for all $i=0,1,\dots,4$.
Again plugging $D_6=D_5=0$,  $a=p$, $b=q$ into what has been done in
Section \ref{sec:pqab}, we get
 the following simplified linear system:
\begin{align*}&
\sm^{-4pq}D_3+D_1=0,
\\&
\sm^{-8pq}D_4+\sm^{-4pq}D_2+D_0=0,
\\&
D_3=0,
\\&
(\sm^{-4pq}+\sm^{-2pq})D_4+D_2=0,
\\&
(\sm^{-4pq}+2\sm^{-2pq}+1)D_4+D_3+D_2=0.\end{align*}
It is now easy to see that the system has only the trivial solution
$D_i=0$ for all $i$.

\section{Case $|p|>q>2$, $|a|>b>2$, $pq=ab$, $p\ne a$}\label{sec:pq=ab}

The calculations carried out in Section \ref{sec:pqab} are largely applicable
to the present section, but there are a few  places where we
can no longer simplify by specializing  the involved functions of
$(t,\sm)$ at $t=-1$, as this results in the zero function.
Additionally, we will  see that for any knot
$K=T(p,q)\#T(a,b)$   in the present  case, its recurrence polynomial $\a_\sk(t,\sm,\sl)$ has repeated factors
involving the variable $\sl$ after evaluation at $t=-1$.

Formula (\ref{3.2}) in Section \ref{sec:pqab} is still valid when $pq=ab$; that is, we have
\begin{equation}\label{equ:pq=ab}
\begin{aligned}
&\left(\sl^2-\sm^{-2ab}t^{-4ab}\right)\frac{\sm^{3ab}t^{6ab}}{\d(a,b,n)}
\left([n+2]\sl^2- \sm^{-4ab}t^{-8ab}[n]\right)
J_{\pq\#\ab}(n)\\
&=\left(\frac{\sm^{-2ab}t^{-4ab}\d(p,q,n+2)}{\d(a,b,n+2)}
-\frac{\sm^{-2ab}t^{-4ab}\d(p,q,n)}{\d(a,b,n)}\right)
J_{\ab}(n)\\
&\;\;\;\;+\frac{\sm^{-ab}t^{-2ab}\d(p,q,n+2)\d(a,b,n)}{\d(a,b,n+2)}
+\sm^{ab}t^{6ab} \delta(p,q,n+2).
\end{aligned}
\end{equation}
We want to show  that the function
\begin{align*} f(t,\sm)&=\frac{\sm^{-2ab}t^{-4ab}\d(p,q,n+2)}{\d(a,b,n+2)}
-\frac{\sm^{-2ab}t^{-4ab}\d(p,q,n)}{\d(a,b,n)}\\&=\frac{\left(\d(p,q,n+2)\d(a,b,n)-\d(p,q,n)\d(a,b,n+2)\right)}
{\sm^{2ab}t^{4ab}\d(a,b,n+2)\d(a,b,n)}
\end{align*}
is not the zero function, even though $f(-1, \sm)=0$.
We can express $f(t,\sm)$ as
\begin{equation}\label{f in fraction}
f(t,\sm)=\frac{\phi(t,\sm)}{\psi(t,\sm)}
\end{equation}
where \begin{align*}
\phi(t,\sm)&=(t^2-t^{-2})^2\left(\d(p,q,n+2)\d(a,b,n)-\d(p,q,n)\d(a,b,n+2)\right),
\\
\psi(t,\sm)&=(t^2-t^{-2})^2\sm^{2ab}t^{4ab}\d(a,b,n+2)\d(a,b,n).
\end{align*}
Note that $\displaystyle\lim_{t\to -1}\phi(t,\sm)=0$ and
$\displaystyle\lim_{t\to -1}\psi(t,\sm)\ne 0$. To show that $f(t,\sm)$ is  non-zero, it suffices to show:
\begin{lemma}\label{lem:f non-zero}
$\phi(t,\sm)$ is not the zero function.
\end{lemma}
\pf By (\ref{equ:delta}), we have
\begin{align*}
\phi(t,\sm)&=(t^2-t^{-2})^2\left(\d(p,q,n+2)\d(a,b,n)-\d(p,q,n)\d(a,b,n+2)\right)
\\&=(\sm^{p+q}t^{6(p+q)+2}+\sm^{-p-q}t^{-6(p+q)+2}-\sm^{q-p}t^{6(q-p)-2}-\sm^{p-q}t^{6(p-q)-2})
\\&\;\;\;\;\;
(\sm^{a+b}t^{2(a+b)+2}+\sm^{-a-b}t^{-2(a+b)+2}-\sm^{b-a}t^{2(b-a)-2}-\sm^{a-b}t^{2(a-b)-2})
\\&
-(\sm^{p+q}t^{2(p+q)+2}+\sm^{-p-q}t^{-2(p+q)+2}-\sm^{q-p}t^{2(q-p)-2}-\sm^{p-q}t^{2(p-q)-2})
\\&\;\;\;\;\;(\sm^{a+b}t^{6(a+b)+2}+\sm^{-a-b}t^{-6(a+b)+2}-\sm^{b-a}t^{6(b-a)-2}-\sm^{a-b}t^{6(a-b)-2})
\end{align*}
We just need to show
that the  Laurent polynomial
\begin{align*}
\phi(t, t^{2n})&=(t^{2(p+q)(n+3)+2}+t^{-2(p+q)(n+3)+2}-t^{2(q-p)(n+3)-2}-t^{-2(q-p)(n+3)-2})\\
&\;\;\;\;\;(t^{2(a+b)(n+1)+2}+t^{-2(a+b)(n+1)+2}-t^{2(a-b)(n+1)-2}-t^{-2(a-b)(n+1)-2})\\
&-(t^{2(p+q)(n+1)+2}+t^{-2(p+q)(n+1)+2}-t^{2(q-p)(n+1)-2}-t^{-2(q-p)(n+1)-2})\\
&\;\;\;\;\;(t^{2(a+b)(n+3)+2}+t^{-2(a+b)(n+3)+2}-t^{2(a-b)(n+3)-2}-t^{-2(a-b)(n+3)-2})
\end{align*}
is not the zero function. This can be seen by considering the degree of
the polynomial.

In fact, first  suppose  $p$ is positive, and thus so is $a$,  and we may assume $p>a$.
In the above expression of $\phi(t, t^{2n})$, the first product  has a lowest degree of $-2(p+q)(n+3)+2-2(a+b)(n+1)+2$,
and the second product has a lowest degree of $-2(p+q)(n+1)+2-2(a+b)(n+3)+2$.
Their difference is
\begin{equation}\label{p>a}
4(a+b)-4(p+q)=4(a+b-p-\frac{ab}{p})=\frac{4(pa+pb-p^2-ab)}{p}
=\frac{4(p-a)(b-p)}{p}<0
\end{equation}
since $p>a>b>0$.
So the lowest degree  of $\phi(t, t^{2n})$ is
$-2(p+q)(n+3)+2-2(a+b)(n+1)+2<0$.

Suppose now  $p<0$, and thus so is $a$, and we may assume $p<a$.
Then the  first product  has a highest degree of $2(q-p)(n+3)-2+2(b-a)(n+1)-2$
and the second product has a highest degree of $2(q-p)(n+1)-2+4(b-a)(n+3)-2$.
Their difference is
\begin{equation}\label{p<a}
4(q-p)-4(b-a)=4(\frac{ab}{p}-p-b+a)=\frac{4(ab-p^2-pb+pa)}{p}=\frac{4(b+p)(a-p)}{p}>0
\end{equation}
 since $p<a<-b<0$.
So the highest degree  of $\phi(t, t^{2n})$ is
$2(q-p)(n+3)-2+2(b-a)(n+1)-2>0$. \qed

Now that $f(t,\sm)\ne 0$,  we may apply the operator $(\sl^2-\sm^{-2ab}t^{-4ab})f^{-1}(t,\sm)$ from the left to (\ref{equ:pq=ab}),
and we have
\begin{align*}
&(\sl^2-\sm^{-2ab}t^{-4ab})f^{-1}(t,\sm)\left(\sl^2-\sm^{-2ab}t^{-4ab}\right)\frac{\sm^{3ab}t^{6ab}}{\d(a,b,n)}
\left([n+2]\sl^2- \sm^{-4ab}t^{-8ab}[n]\right)
J_{\pq\#\ab}(n)\\&
=(\sl^2-\sm^{-2ab}t^{-4ab})
J_{\ab}(n)+(\sl^2-\sm^{-2ab}t^{-4ab})f^{-1}(t,\sm)\frac{\sm^{-ab}t^{-2ab}\d(p,q,n+2)\d(a,b,n)}{\d(a,b,n+2)}
+\sm^{ab}t^{6ab} \delta(p,q,n+2)\\&
=\sm^{ab}t^{-2ab}\d(a,b,n)+f^{-1}(t,t^4\sm)\left(\frac{\sm^{-ab}t^{-6ab}\d(p,q,n+4)\d(a,b,n+2)}{\d(a,b,n+4)}+\sm^{ab}t^{10ab}\d(p,q,n+4)\right)\\&
\;\;-\sm^{-2ab}t^{-4ab}f^{-1}(t,\sm)\left(\frac{\sm^{-ab}t^{-2ab}\d(p,q,n+2)\d(a,b,n)}{\d(a,b,n+2)}
+\sm^{ab}t^{6ab}\d(p,q,n+2)\right).
\end{align*}
Let $g(t,\sm)$ be the function on the right-hand side of the last equality.
Then $g(t,\sm)$ is non-zero since
$\displaystyle\lim_{t\ra-1}(t^2-t^{-2})g(t, \sm)=\infty$ (see (\ref{equ:g(-1,m)})).
Therefore,
\begin{align*}
&\a(t,\sm,\sl)\\&
=\left(\sl-1\right)\frac{1}{\left(t^2-t^{-2}\right)g(t,\sm)}\left(\sl^2-\sm^{-2ab}t^{-4ab}\right)f^{-1}(t,\sm)\left(\sl^2-\sm^{-2ab}t^{-4ab}\right)\frac{\sm^{3ab}t^{6ab}}{\d(a,b,n)}
\left([n+2]\sl^2- \sm^{-4ab}t^{-8ab}[n]\right)\\
&=\left(\sl-1\right)\frac{f^{-1}(t, t^4\sm)}{\left(t^2-t^{-2}\right)g(t,\sm)}\left(\sl^2-\frac{f^{-1}(t,\sm)}{f^{-1}(t,t^4\sm)}\sm^{-2ab}t^{-4ab}\right)
\left(\sl^2-\sm^{-2ab}t^{-4ab}\right)\frac{\sm^{3ab}t^{6ab}}{\d(a,b,n)}
\left([n+2]\sl^2- \sm^{-4ab}t^{-8ab}[n]\right)
\end{align*}
is an annihilator of $J_{\pq\#\ab}(n)$, where  the second equality follows from the product formula (\ref{equ:product}).
To eveluate $\a(t,\sm,\sl)$ at $t=-1$, we need the following lemma:
\begin{lemma}\label{lem:limit=1}
$\displaystyle \lim_{t\ra -1}\frac{f^{-1}((t,\sm)}{f^{-1}(t,t^4\sm)}=1.$
\end{lemma}

\pf  By (\ref{f in fraction}) and the note following it, we have  \[ \lim_{t\ra -1}\frac{f^{-1}(t,\sm)}{f^{-1}(t,t^4\sm)}
=\lim_{t\ra-1}\frac{f(t,t^4\sm)}{f(t,\sm)}=\lim_{t\ra -1}\frac{\phi(t, t^4\sm)\psi(t,\sm)}{\psi(t, t^4\sm)\phi(t,\sm)}
=\lim_{t\ra -1}\frac{\phi(t, t^4\sm)}{\phi(t,\sm)}.\]
Applying   L'H\^opital's rule (treating $\sm$ as an independent variable), we have
 \[ \lim_{t\ra -1}\frac{\phi(t, t^4\sm)}{\phi(t,\sm)}
= \lim_{t\ra -1}\frac{\dfrac{\p\phi(t, t^4\sm)}{\p t}}{\dfrac{\p \phi(t,\sm)}{\p t}}.\]
Now \begin{align*}
& \lim_{t\ra -1}\frac{\p\phi(t, \sm)}{\p t}\\
&=\lim_{t\ra -1}\frac{\p}{\p t}\left[\left(\sm^{p+q}t^{6(p+q)+2}+\sm^{-p-q}t^{-6(p+q)+2}-\sm^{q-p}t^{6(q-p)-2}-\sm^{p-q}t^{6(p-q)-2}\right)\right.\\&\;\;\;\;\;\;\;\;\;\;\;\;\;\;\;\;\;\;\;\;
\left(\sm^{a+b}t^{2(a+b)+2}+\sm^{-a-b}t^{-2(a+b)+2}-\sm^{b-a}t^{2(b-a)-2}-\sm^{a-b}t^{2(a-b)-2}\right)\\
&\;\;\;\;\;\;\;\;\;\;\;\;\;\;\;\;-\left(\sm^{p+q}t^{2(p+q)+2}+\sm^{-p-q}t^{-2(p+q)+2}-\sm^{q-p}t^{2(q-p)-2}-\sm^{p-q}t^{2(p-q)-2}\right)\\&\;\;\;\;\;\;\;\;\;\;\;\;\;\;\;\;\;\;\;\;\left.\left(\sm^{a+b}t^{6(a+b)+2}+\sm^{-a-b}t^{-6(a+b)+2}-\sm^{b-a}t^{6(b-a)-2}-\sm^{a-b}t^{6(a-b)-2}\right)\right]\\
&= \lim_{t\ra -1}
\left\{\left[(6p+6q+2)\sm^{p+q}t^{6(p+q)+1}
+(-6p-6q+2)\sm^{-p-q}t^{-6(p+q)+1}\right.\right.\\
&\;\;\;\;\;\;\;\;\;\;\;\;\;\;\;\;\;\left.-(6q-6p-2)\sm^{q-p}t^{6(q-p)-3}-(6p-6q-2)\sm^{p-q}t^{6(p-q)-3}\right]
\\&\;\;\;\;\;\;\;\;\;\;\;\;\;\;\;\;
\left(\sm^{a+b}t^{2(a+b)+2}+\sm^{-a-b}t^{-2(a+b)+2}-\sm^{b-a}t^{2(b-a)-2}-\sm^{a-b}t^{2(a-b)-2}\right)
\\&\;\;\;\;\;\;\;\;\;\;\;\;\;\;+\left(\sm^{p+q}t^{6(p+q)+2}+\sm^{-p-q}t^{-6(p+q)+2}-\sm^{q-p}t^{6(q-p)-2}-\sm^{p-q}t^{6(p-q)-2}\right)
\\&\;\;\;\;\;\;\;\;\;\;\;\;\;\;\;\;\;\;
\left[(2a+2b+2)\sm^{a+b}t^{2(a+b)+1}+(-2a-2b+2)\sm^{-a-b}t^{-2(a+b)+1}\right.\\&\;\;\;\;\;\;\;\;\;\;\;\;\;\;\;\;
\;\;\;\;\left.-
(2b-2a-2)\sm^{b-a}t^{2(b-a)-3}-(2a-2b-2)\sm^{a-b}t^{2(a-b)-3}\right]
\\&\;\;\;\;\;\;\;\;\;\;\;\;\;\;\;
-\left[(2p+2q+2)\sm^{p+q}t^{2(p+q)+1}+(-2p-2q+2)\sm^{-p-q}t^{-2(p+q)+1}\right.\\&\;\;\;\;\;\;\;\;\;\;\;\;\;\;\;\;\;\;\;\;\;\left.-(2q-2p-2)\sm^{q-p}t^{2(q-p)-3}-(2p-2q-2)\sm^{p-q}t^{2(p-q)-3}\right]
\\&\;\;\;\;\;\;\;\;\;\;\;\;\;\;\;\;\;\;\;\;\left(\sm^{a+b}t^{6(a+b)+2}+\sm^{-a-b}t^{-6(a+b)+2}-\sm^{b-a}t^{6(b-a)-2}-\sm^{a-b}t^{6(a-b)-2}\right)
\\&\;\;\;\;\;\;\;\;\;\;\;\;\;\;\;\;-\left(\sm^{p+q}t^{2(p+q)+2}+\sm^{-p-q}t^{-2(p+q)+2}-\sm^{q-p}t^{2(q-p)-2}-\sm^{p-q}t^{2(p-q)-2}\right)\\&\;\;\;\;\;\;\;\;\;\;\;\;\;\;\;\;\;\;\;\;\left[(6a+6b+2)\sm^{a+b}t^{6(a+b)+1}
+(-6a-6b+2)\sm^{-a-b}t^{-6(a+b)+1}\right.\\
&\;\;\;\;\;\;\;\;\;\;\;\;\;\;\;\;\;\;\;\;\;\left.\left.-(6b-6a-2)\sm^{b-a}t^{6(b-a)-3}-(6a-6b-2)\sm^{a-b}t^{6(a-b)-3}\right]\right\}
\\&=
\left[-(6p+6q+2)\sm^{p+q}
-(-6p-6q+2)\sm^{-p-q}+(6q-6p-2)\sm^{q-p}+(6p-6q-2)\sm^{p-q}\right]
\\&\;\;\;\;\;
\left(\sm^{a+b}+\sm^{-a-b}-\sm^{b-a}-\sm^{a-b}\right)
+\left(\sm^{p+q}+\sm^{-p-q}-\sm^{q-p}-\sm^{p-q}\right)
\\&\;\;\;\;\;
\left[-(2a+2b+2)\sm^{a+b}-(-2a-2b+2)\sm^{-a-b}+(2b-2a-2)\sm^{b-a}+(2a-2b-2)\sm^{a-b}\right]
\\&\;
-\left[-(2p+2q+2)\sm^{p+q}-(-2p-2q+2)\sm^{-p-q}+(2q-2p-2)\sm^{q-p}+(2p-2q-2)\sm^{p-q}\right]
\\&\;\;\;\;\;\;\left(\sm^{a+b}+\sm^{-a-b}-\sm^{b-a}-\sm^{a-b}\right)
-\left(\sm^{p+q}+\sm^{-p-q}-\sm^{q-p}-\sm^{p-q}\right)\\&\;\;\;\;\;\;\left[-(6a+6b+2)\sm^{a+b}
-(-6a-6b+2)\sm^{-a-b}+(6b-6a-2)\sm^{b-a}+(6a-6b-2)\sm^{a-b}\right]
\\&
=\left(\sm^{a+b}+\sm^{-a-b}-\sm^{b-a}-\sm^{a-b}\right)
\\&\;\;\;\;\left[-(4p+4q)\sm^{p+q}+(4p+4q)\sm^{-p-q}+
(4q-4p)\sm^{q-p}+(4p-4q)\sm^{p-q}\right]
\\&\;+\left(\sm^{p+q}+\sm^{-p-q}-\sm^{q-p}-\sm^{p-q}\right)
\\&\;\;\;\;\;\left[(4a+4b)\sm^{a+b}+(-4a-4b)\sm^{-a-b}-(4b-4a)\sm^{b-a}-(4a-4b)\sm^{a-b}\right]
\end{align*}
which is non-zero. In fact, when $p>0$ and thus $a>0$, we may assume $p>a$. The highest degree term of this  polynomial is
$-4(p+q-a-b)\sm^{p+q+a+b}$, which is non-zero since $p+q\ne a+b$ by (\ref{p>a}). Similarly,
 when $p<0$ and thus $a<0$, we may assume $p<a$; the highest degree term of this polynomial
 is $-4(q-p-b+a)\sm^{q-p+ b-a}$, which is non-zero  since $q-p\ne b-a$ by (\ref{p<a}).

On the other hand,
\begin{align*}
&
 \lim_{t\ra -1}\frac{\p\phi(t, t^4\sm)}{\p t}\\
&=\lim_{t\ra -1}\frac{\p}{\p t}\left[\left(\sm^{p+q}t^{10(p+q)+2}+\sm^{-p-q}t^{-10(p+q)+2}-\sm^{q-p}t^{10(q-p)-2}-\sm^{p-q}t^{10(p-q)-2}\right)\right.\\&\;\;\;\;\;\;\;\;\;\;\;\;\;\;\;\;\;\;
\left(\sm^{a+b}t^{6(a+b)+2}+\sm^{-a-b}t^{-6(a+b)+2}-\sm^{b-a}t^{6(b-a)-2}-\sm^{a-b}t^{6(a-b)-2}\right)\\
&\;\;\;\;\;\;\;\;\;\;\;\;\;\;\;-\left(\sm^{p+q}t^{6(p+q)+2}+\sm^{-p-q}t^{-6(p+q)+2}-\sm^{q-p}t^{6(q-p)-2}-\sm^{p-q}t^{6(p-q)-2}\right)\\&\;\;\;\;\;\;\;\;\;\;\;\;\;\;\;\;\;\;\;\left.\left(\sm^{a+b}t^{10(a+b)+2}+\sm^{-a-b}t^{-10(a+b)+2}-\sm^{b-a}t^{10(b-a)-2}-\sm^{a-b}t^{10(a-b)-2}\right)\right]\\&
= \lim_{t\ra -1}
\left\{\left[(10p+10q+2)\sm^{p+q}t^{10(p+q)+1}
+(-10p-10q+2)\sm^{-p-q}t^{-10(p+q)+1}\right.\right.\\
&\;\;\;\;\;\;\;\;\;\;\;\;\;\;\;\;\left.-(10q-10p-2)\sm^{q-p}t^{10(q-p)-3}-(10p-10q-2)\sm^{p-q}t^{10(p-q)-3}\right]
\\&\;\;\;\;\;\;\;\;\;\;\;\;\;\;\;\;
\left(\sm^{a+b}t^{6(a+b)+2}+\sm^{-a-b}t^{-6(a+b)+2}-\sm^{b-a}t^{6(b-a)-2}-\sm^{a-b}t^{6(a-b)-2}\right)
\\&\;\;\;\;\;\;\;\;\;\;\;\;\;+\left(\sm^{p+q}t^{10(p+q)+2}+\sm^{-p-q}t^{-10(p+q)+2}-\sm^{q-p}t^{10(q-p)-2}-\sm^{p-q}t^{10(p-q)-2}\right)
\\&\;\;\;\;\;\;\;\;\;\;\;\;\;\;\;\;\;\;
\left[(6a+6b+2)\sm^{a+b}t^{6(a+b)+1}+(-6a-6b+2)\sm^{-a-b}t^{-6(a+b)+1}\right.\\&\;\;\;\;\;\;\;\;\;\;\;\;\;\;\;\;\;\;
\left.-(6b-6a-2)\sm^{b-a}t^{6(b-a)-3}-(6a-6b-2)\sm^{a-b}t^{6(a-b)-3}\right]
\\&\;\;\;\;\;\;\;\;\;\;\;\;\;\;
-\left[(6p+6q+2)\sm^{p+q}t^{6(p+q)+1}+(-6p-6q+2)\sm^{-p-q}t^{-6(p+q)+1}\right.\\&\;\;\;\;\;\;\;\;\;\;\;\;\;\;\;\;\;\;\;\;
\left.-(6q-6p-2)\sm^{q-p}t^{6(q-p)-3}-(6p-6q-2)\sm^{p-q}t^{6(p-q)-3}\right]
\\&\;\;\;\;\;\;\;\;\;\;\;\;\;\;\;\;\;\;\;\;\left(\sm^{a+b}t^{10(a+b)+2}+\sm^{-a-b}t^{-10(a+b)+2}-\sm^{b-a}t^{10(b-a)-2}-\sm^{a-b}t^{10(a-b)-2}\right)
\\&\;\;\;\;\;\;\;\;\;\;\;\;\;\;\;\;-\left(\sm^{p+q}t^{6(p+q)+2}+\sm^{-p-q}t^{-6(p+q)+2}-\sm^{q-p}t^{6(q-p)-2}-\sm^{p-q}t^{6(p-q)-2}\right)\\&\;\;\;\;\;\;\;\;\;\;\;\;\;\;\;\;\;\;\;\;\left[(10a+10b+2)\sm^{a+b}t^{10(a+b)+1}
+(-10a-10b+2)\sm^{-a-b}t^{-10(a+b)+1}\right.\\&
\;\;\;\;\;\;\;\;\;\;\;\;\;\;\;\;\;\;\;\;\left.\left.-(10b-10a-2)\sm^{b-a}t^{10(b-a)-3}-(10a-10b-2)\sm^{a-b}t^{10(a-b)-3}\right]\right\}
\\&=\left[-(10p+10q+2)\sm^{p+q}
-(-10p-10q+2)\sm^{-p-q}+(10q-10p-2)\sm^{q-p}+(10p-10q-2)\sm^{p-q}\right]
\\&\;\;\;\;\;
\left(\sm^{a+b}+\sm^{-a-b}-\sm^{b-a}-\sm^{a-b}\right)
+\left(\sm^{p+q}+\sm^{-p-q}-\sm^{q-p}-\sm^{p-q}\right)
\\&\;\;\;\;\;
\left[-(6a+6b+2)\sm^{a+b}-(-6a-6b+2)\sm^{-a-b}+(6b-6a-2)\sm^{b-a}+(6a-6b-2)\sm^{a-b}\right]
\\&\;
-\left[-(6p+6q+2)\sm^{p+q}-(-6p-6q+2)\sm^{-p-q}+(6q-6p-2)\sm^{q-p}+(6p-6q-2)\sm^{p-q}\right]
\\&\;\;\;\;\;\left(\sm^{a+b}+\sm^{-a-b}-\sm^{b-a}-\sm^{a-b})
-(\sm^{p+q}+\sm^{-p-q}-\sm^{q-p}-\sm^{p-q}\right)\\&\;\;\;\;\;\left[-(10a+10b+2)\sm^{a+b}
-(-10a-10b+2)\sm^{-a-b}+(10b-10a-2)\sm^{b-a}+(10a-10b-2)\sm^{a-b}\right]
\\&
=\left(\sm^{a+b}+\sm^{-a-b}-\sm^{b-a}-\sm^{a-b}\right)
\\&\;\;\;\;\left[-(4p+4q)\sm^{p+q}+(4p+4q)\sm^{-p-q}+
(4q-4p)\sm^{q-p}+(4p-4q)\sm^{p-q}\right]
\\&\;+\left(\sm^{p+q}+\sm^{-p-q}-\sm^{q-p}-\sm^{p-q}\right)
\\&\;\;\;\;\;\left[(4a+4b)\sm^{a+b}+(-4a-4b)\sm^{-a-b}-(4b-4a)\sm^{b-a}-(4a-4b)\sm^{a-b}\right]
\\&= \lim_{t\ra -1}\frac{\p\phi(t, \sm)}{\p t}.\end{align*}
Hence Lemma \ref{lem:limit=1} is true. \qed

Now we evaluate
\[\lim_{t\ra -1}\frac{f^{-1}(t, t^4\sm)}{(t^2-t^{-2})g(t,\sm)}.
\]
Applying Lemma \ref{lem:limit=1} and the fact $f(-1,\sm)=0$, we have
\begin{align*}
&\lim_{t\ra -1}f(t,t^4\sm)(t^2-t^{-2})g(t,\sm)\\
&=\lim_{t\ra-1}f(t,t^4\sm)(t^2-t^{-2})
\left[\sm^{ab}t^{-2ab}\d(a,b,n)+f^{-1}(t,t^4\sm)\left(\frac{\sm^{-ab}t^{-6ab}\d(p,q,n+4)\d(a,b,n+2)}{\d(a,b,n+4)}
\right.\right.\\&\;\;+\left.\sm^{ab}t^{10ab}\d(p,q,n+4)\right)
\left.-\sm^{-2ab}t^{-4ab}f^{-1}(t,\sm)\left(\frac{\sm^{-ab}t^{-2ab}\d(p,q,n+2)\d(a,b,n)}{\d(a,b,n+2)}
+\sm^{ab}t^{6ab}\d(p,q,n+2)\right)\right]
\\
&=\lim_{t\ra-1}(t^2-t^{-2})\left[\left(\frac{\sm^{-ab}t^{-6ab}\d(p,q,n+4)\d(a,b,n+2)}{\d(a,b,n+4)}+\sm^{ab}t^{10ab}\d(p,q,n+4)\right)\right.\\
&\;\;\left.-\sm^{-2ab}t^{-4ab}\left(\frac{\sm^{-ab}t^{-2ab}\d(p,q,n+2)\d(a,b,n)}{\d(a,b,n+2)}
+\sm^{ab}t^{6ab}\d(p,q,n+2)\right)\right]
\\
&=(1+\sm^{-2ab})(\sm^{-ab}+\sm^{ab})(\sm^{p+q}+\sm^{-p-q}-\sm^{q-p}-\sm^{p-q}).
\end{align*}
Hence
\[\lim_{t\ra -1}\frac{f^{-1}(t, t^4\sm)}{(t^2-t^{-2})g(t,\sm)}
=\frac{1}{(1+\sm^{-2ab})(\sm^{-ab}+\sm^{ab})(\sm^{p+q}+\sm^{-p-q}-\sm^{q-p}-\sm^{p-q})}.\]
Therefore,
\[\a(-1,\sm,\sl)=a(\sm)(\sl-1)(\sl^2-\sm^{-2ab})(\sl^2-\sm^{-2ab})(\sl^2-\sm^{-4ab})
\]
where \[a(\sm)=\frac{\sm^{3ab}(\sm-\sm^{-1})}
{(1+\sm^{-2ab})(\sm^{-ab}+\sm^{ab})(\sm^{p+q}+\sm^{-p-q}-\sm^{q-p}-\sm^{p-q})
(\sm^{a+b}+\sm^{-a-b}-\sm^{b-a}-\sm^{a-b})}\]
is a well-defined, non-zero function of $\sm$.
Clearly,  $(\sl^2-\sm^{-2ab})$ is a repeated factor of $\a(-1,\sm,\sl)$, and
after deleting it,  $\a(-1,\sm,\sl)$ is equal to the $A$-polynomial of
$T(p,q)\#T(a,b)$ (see Lemma \ref{lem:A-poly of sum} (ii)) up to a non-zero factor depending on $\sm$ only.

We now show that $7$ is the minimum $\sl$-degree among all non-zero annihilators of
$J_{\pq\#\ab}(n)$.
So, suppose $$\begin{array}{l}D_6J_{\pq\#\ab}(n+6) +D_5J_{\pq\#\ab}(n+5)+D_4J_{\pq\#\ab}(n+4)
+D_3J_{\pq\#\ab}(n+3)\\+D_2J_{\pq\#\ab}(n+2)+D_1J_{\pq\#\ab}(n+1)+D_0J_{\pq\#\ab}
=0\end{array}$$
for some $D_i\in\q(t,\sm)$, $i=0,1,\dots,6$. We need to show that all $D_i$ have to be zero.

Note that  in Section  \ref{sec:pqab}, the calculation of
\begin{align*}&D_6J_{\pq\#\ab}(n+6) +D_5J_{\pq\#\ab}(n+5)+D_4J_{\pq\#\ab}(n+4)
+D_3J_{\pq\#\ab}(n+3)\\&+D_2J_{\pq\#\ab}(n+2)+D_1J_{\pq\#\ab}(n+1)+D_0J_{\pq\#\ab}
\\=\;&E_6J_{\pq\#\ab}(n+1)+E_5J_{\pq\#\ab}(n)+E_4J_{\pq}(n+1)
+E_3J_{\pq}(n)+E_2J_{\ab}(n+1)\\&+E_1J_{\ab}(n)+E_0\end{align*}
is still valid when $pq=ab$, and so is  Lemma \ref{lem:Eizero}; i.e.,
we must have $E_i=0$ for all $i=0,1,\dots,6$.
However, when $pq=ab$, the argument
for the assertion that $E_i=0$ for all $i$ implies $D_i=0$ for all  $i$
needs some modification, because specializing coefficients at $t=-1$ makes
some equations linearly dependent.

Plugging $pq=ab$ into the system $E_i=0$, $i=0,1,\dots,6$,  obtained in  Section  \ref{sec:pqab},
and simplifying    $E_1=0$ and $E_0=0$ as in Section  \ref{sec:pqab},
we get
\begin{align*}
& E_6=\left(\frac{D_5}{[n+5]} t^{-8ab(n+4)}+\frac{D_3}{[n+3]}\right)
 t^{-8ab(n+2)}[n+1]+D_1=0,\\
& E_5=\left[\left(\frac{D_6}{[n+6]} t^{-8ab(n+5)}+\frac{D_4}{[n+4]}\right)
 t^{-8ab(n+3)}+\frac{D_2}{[n+2]}\right]t^{-8ab(n+1)}[n]+D_0=0,
\\&
 E_4=\left(\frac{D_5}{[n+5]} t^{-8ab(n+4)}+\frac{D_3}{[n+3]}\right)
t^{-6ab(n+2)}\delta(a,b,n+1)+\frac{D_5}{[n+5]}t^{-6ab(n+4)}\delta(a,b,n+3)
 t^{-4ab(n+2)}
\\&
\;\;\;\;\;=\left(\frac{t^{-14abn-44ab}\d(a,b,n+1)+t^{-10abn-32ab}\d(a,b,n+3)}{[n+5]}\right)D_5+\frac{t^{-6ab(n+2)}\d(a,b,n+1)}{[n+3]}D_3=0,
\\& E_3=\left[\left(\frac{D_6}{[n+6]} t^{-8ab(n+5)}+\frac{D_4}{[n+4]}\right)
 t^{-8ab(n+3)}+\frac{D_2}{[n+2]}\right]
t^{-6ab(n+1)}\delta(a,b,n)\\
& \;\;\;\;\;\;\;\;\;\; +\left[\left(\frac{D_6}{[n+6]} t^{-8ab(n+5)}+\frac{D_4}{[n+4]}\right)
 t^{-6ab(n+3)}\delta(a,b,n+2)\right.\\
&\;\;\;\;\;\;\;\;\;\;\;\;\;\;  +
\left.\frac{D_6}{[n+6]}t^{-6ab(n+5)}\delta(a,b,n+4)
 t^{-4ab(n+3)} \right]
t^{-4ab(n+1)}\\
&\;\;\;\;\;=
\left(\frac{t^{-22abn-70ab}\d(a,b,n)+t^{-18abn-62ab}\d(a,b,n+2)
+t^{-14abn-46ab}\d(a,b,n+4)}{[n+6]}\right)D_6\\&
\;\;\;\;\;\;\;\;+\left(\frac{t^{-14abn-30ab}\d(a,b,n)+t^{-10abn-22ab}\d(a,b,n+2)
}{[n+4]}\right)D_4+\frac{t^{-6ab(n+1)}\d(a,b,n)}{[n+2]}D_2=0,
\\&
E_2=\left(\frac{D_5}{[n+5]} t^{-8ab(n+4)}+\frac{D_3}{[n+3]}\right)
t^{-6ab(n+2)}\delta(p,q,n+1)
+\frac{D_5}{[n+5]}t^{-6ab(n+4)}\delta(p,q,n+3)
 t^{-4ab(n+2)}\\
&\;\;\;\;=\left(\frac{t^{-14abn-44ab}\d(p,q,n+1)+t^{-10abn-32ab}\d(p,q,n+3)}{[n+5]}\right)D_5+\frac{t^{-6ab(n+2)}\d(p,q,n+1)}{[n+3]}D_3=0,
\\&
 E_1=\left[\left(\frac{D_6}{[n+6]} t^{-8ab(n+5)}+\frac{D_4}{[n+4]}\right)
 t^{-8ab(n+3)}+\frac{D_2}{[n+2]}\right]
 t^{-6ab(n+1)}\delta(p,q,n)
\\& \;\;\;\; +\left[\left(\frac{D_6}{[n+6]} t^{-8ab(n+5)}+\frac{D_4}{[n+4]}\right)\right.t^{-6ab(n+3)}\delta(p,q,n+2)
+\left.\frac{D_6}{[n+6]}t^{-6ab(n+5)}\delta(p,q,n+4)
t^{-4ab(n+3)} \right]t^{-4ab(n+1)}=0
\\&
\;\lra\;\; (\sm^{-8ab}+\sm^{-6ab}+\sm^{-4ab})D_6
+(\sm^{-4ab}+\sm^{-2ab})D_4+D_2=0.
\\&  E_0=
\left[\left(\frac{D_6}{[n+6]} t^{-8ab(n+5)}+\frac{D_4}{[n+4]}\right)
 t^{-8ab(n+3)}+\frac{D_2}{[n+2]}\right]
t^{-4ab(n+1)}\delta(p,q,n)\delta(a,b,n)
\\&\;\;\;\;\;\;
+\left[\left(\frac{D_6}{[n+6]} t^{-8ab(n+5)}+\frac{D_4}{[n+4]}\right)
 t^{-6ab(n+3)}\delta(a,b,n+2)\right.\\&\;\;
\;\;\;\;\;\;\;\;\;+\left.\frac{D_6}{[n+6]}t^{-6ab(n+5)}\delta(a,b,n+4)
 t^{-4ab(n+3)} \right]
t^{-2ab(n+1)}\d(p,q,n)
\\&\;\;\;\;\;
+\left[\left(\frac{D_6}{[n+6]} t^{-8ab(n+5)}+\frac{D_4}{[n+4]}\right)\right.t^{-6ab(n+3)}\delta(p,q,n+2)
\\&\;\;\;\;\;\;\;\;\;\;
+\left.\frac{D_6}{[n+6]}t^{-6ab(n+5)}\delta(p,q,n+4)
t^{-4ab(n+3)} \right]
t^{-2ab(n+1)}\d(a,b,n)\\
&\;\;\;\;\;\;+\left(\frac{D_6}{[n+6]} t^{-8ab(n+5)}+\frac{D_4}{[n+4]}\right)t^{-4ab(n+3)}\delta(p,q,n+2)\delta(a,b,n+2)
\\&\;\;\;\;\;\;
+\left(\frac{D_5}{[n+5]} t^{-8ab(n+4)}+\frac{D_3}{[n+3]}\right)
t^{-4ab(n+2)}\delta(p,q,n+1)\delta(a,b,n+1)
\\&\;\;\;\;\;\;
+\frac{D_6}{[n+6]}t^{-6ab(n+5)}\delta(a,b,n+4)
 t^{-2ab(n+3)} \delta(p,q,n+2)
+\frac{D_6}{[n+6]}t^{-6ab(n+5)}\delta(p,q,n+4)
 t^{-2ab(n+3)} \delta(a,b,n+2)
\\&\;\;\;\;\;\;
+\frac{D_6}{[n+6]}t^{-4ab(n+5)}\delta(p,q,n+4)\delta(a,b,n+4)
+\frac{D_5}{[n+5]}t^{-6ab(n+4)}\delta(a,b,n+3)
t^{-2ab(n+2)} \delta(p,q,n+1)
\\&\;\;\;\;\;\;+\frac{D_5}{[n+5]}t^{-6ab(n+4)}\delta(p,q,n+3)
t^{-2ab(n+2)} \delta(a,b,n+1)
+\frac{D_5}{[n+5]}t^{-4ab(n+4)}\delta(p,q,n+3)\delta(a,b,n+3)=0\\
&\;\lra\;\;
(\sm^{-8ab}+2\sm^{-6ab}+3\sm^{-4ab}
+2\sm^{-2ab}
+1)D_6
+(\sm^{-4ab}+2\sm^{-2ab}
+1)D_5\\&\;\;\;\;\;\;\;\;+(\sm^{-4ab}+2\sm^{-2ab}+1)D_4+D_3+D_2=0.
\end{align*}

The determinant of the coefficient matrix for the system $E_4=0$ and $E_2=0$ (in variables $D_5$ and $D_3$)
is
\begin{align*}&\frac{\left(t^{-14abn-44ab}\d(a,b,n+1)+
t^{-10abn-32ab}\d(a,b,n+3)\right)t^{-6ab(n+2)}\d(p,q,n+1)}{[n+5][n+3]}\\&
-\frac{\left(t^{-14abn-44ab}\d(p,q,n+1)+
t^{-10abn-32ab}\d(p,q,n+3)\right)t^{-6ab(n+2)}\d(a,b,n+1)}{[n+5][n+3]}
\\&
=\frac{ t^{-16abn-44ab}\left(\d(a,b,n+3)\d(p,q,n+1)-
\d(p,q,n+3)\d(a,b,n+1)\right)}{[n+5][n+3]},
\end{align*}
which is non-zero by Lemma \ref{lem:f non-zero} (where we replace $n$ with $n+1$).
It follows that  $D_3=D_5=0$. Substituting these into $E_6=0$, we find that  $D_1=0$.

From the simplified $E_1=0$,
we get $D_2=-(\sm^{-8ab}+\sm^{-6ab}+\sm^{-4ab})D_6
-(\sm^{-4ab}+\sm^{-2ab})D_4$.
Plugging it into $E_3=0$ and the simplified $E_0=0$, we get
\begin{align*}&\left(\frac{t^{-22abn-70ab}\d(a,b,n)+t^{-18abn-62ab}\d(a,b,n+2)
+t^{-14abn-46ab}\d(a,b,n+4)}{[n+6]}\right.
\\&\left.\;\;-\frac{t^{-6ab(n+1)}\d(a,b,n)\left(\sm^{-8ab}
+\sm^{-6ab}+\sm^{-4ab}\right)}{[n+2]}\right)D_6\\
&+\left(\frac{t^{-14abn-30ab}\d(a,b,n)+t^{-10abn-22ab}\d(a,b,n+2)
}{[n+4]}-\frac{t^{-6ab(n+1)}\d(a,b,n)\left(\sm^{-4ab}+\sm^{-2ab}\right)}{[n+2]}\right)D_4
=0
\end{align*}
and
\begin{align*}
\left(\sm^{-6ab}+2\sm^{-4ab}
+2\sm^{-2ab}+1\right)D_6+\left(\sm^{-2ab}+1\right)D_4=0.
\end{align*}
The determinant of the coefficient matrix of these two equations is
\begin{align*}
&\;\left(\frac{t^{-22abn-70ab}\d(a,b,n)+t^{-18abn-62ab}\d(a,b,n+2)
+t^{-14abn-46ab}\d(a,b,n+4)}{[n+6]}\right.
\\&\left.\;\;\;-\frac{t^{-6ab(n+1)}\d(a,b,n)\left(\sm^{-8ab}
+\sm^{-6ab}+\sm^{-4ab}\right)}{[n+2]}\right)
\left(\sm^{-2ab}+1\right)\\
&-\left(\frac{t^{-14abn-30ab}\d(a,b,n)+t^{-10abn-22ab}\d(a,b,n+2)
}{[n+4]}-\frac{t^{-6ab(n+1)}\d(a,b,n)\left(\sm^{-4ab}+\sm^{-2ab}\right)}{[n+2]}\right)
\\&\;\;\;\;\;\left(\sm^{-6ab}+2\sm^{-4ab}
+2\sm^{-2ab}+1\right)
\\&
=\;\frac{\left(t^{-22abn-70ab}\d(a,b,n)+t^{-18abn-62ab}\d(a,b,n+2)
+t^{-14abn-46ab}\d(a,b,n+4)\right)\left(\sm^{-2ab}+1\right)}{[n+6]}
\\&\;\;\;\;-\frac{\left(t^{-14abn-30ab}\d(a,b,n)+t^{-10abn-22ab}\d(a,b,n+2)\right)
\left(\sm^{-6ab}+2\sm^{-4ab}
+2\sm^{-2ab}+1\right)}{[n+4]}\\&
\;\;\;\;+\frac{t^{-6ab(n+1)}\d(a,b,n)\left(\sm^{-8ab}+2\sm^{-6ab}+2\sm^{-4ab}+\sm^{-2ab}\right)}{[n+2]}.
\end{align*}

\begin{lemma}\label{lem:nonzero det}The preceding determinant function is non-zero.
\end{lemma}
\pf  The function can be expressed as a fraction $\displaystyle \frac{\phi(t, \sm)}{\psi(t,\sm)}$ where
\begin{align*}
&\phi(t,\sm)\\&=(t^2-t^{-2})^3\Big[[n+4][n+2]\left(t^{-22abn-70ab}\d(a,b,n)+t^{-18abn-62ab}\d(a,b,n+2)
+t^{-14abn-46ab}\d(a,b,n+4)\right).\\
&\;\;\;\;\left(\sm^{-2ab}+1\right)-[n+6][n+2]\left(t^{-14abn-30ab}\d(a,b,n)+t^{-10abn-22ab}\d(a,b,n+2)\right)\\&\;\;\;\;
\left(\sm^{-6ab}+2\sm^{-4ab}
+2\sm^{-2ab}+1\right)+[n+6][n+4]t^{-6ab(n+1)}\d(a,b,n)\left(\sm^{-8ab}+2\sm^{-6ab}+2\sm^{-4ab}+\sm^{-2ab}\right)\Big]
\end{align*}
and
\begin{align*}
\psi(t,\sm)&=(t^2-t^{-2})^3[n+6][n+4][n+2].
\end{align*}
Since  $\psi(t,\sm)$ is non-zero, it suffices to show that $\phi(t,\sm)$ is non-zero
(although $\phi(-1,\sm)=0$).
This can be seen by considering the degree of the Laurent polynomial
$\phi(t, t^{2n})$.

In fact, when $a>0$, the first product  in the  expression of $\phi(t,\sm)$
(after substituting $\sm$ with $t^{2n}$) has a lowest degree of
$-(4+26ab+2a+2b)n-10-70ab-2a-2b$,
the second product has a lowest degree of
$-(4+26ab+2a+2b)n-14-30ab-2a-2b$,
and the third product has a lowest degree of
$-(4+22ab+2a+2b)n-18-2a-2b$.
Hence, the lowest degree of $\phi(t, t^{2n})$  is
$-(4+26ab+2a+2b)n-10-70ab-2a-2b<0$.

Similarly, when $a<0$, the first product of $\phi(t,t^{2n})$
has a highest degree of $-(4+26ab+2a-2b)n-14-70ab-2a+2b$,
the second product has a highest degree of
$-(4+26ab+2a-2b)n-18-30ab-2a+2b$,
and the third product has a highest degree of
$-(4+22ab+2a-2b)n-22-2a+2b$.
Hence,  the highest degree of $\phi(t,t^{2n})$  is
$-(4+26ab+2a-2b)n-14-70ab-2a+2b>0$.
\qed

By Lemma \ref{lem:nonzero det}, we have $D_6=D_4=0$,
which in turn implies  $D_2=0$. Then  from $E_5=0$, we have $D_0=0$.

\section{Case $|p|>q>2$, $|a|>b=2$, $pq\ne 2a$, $p, a$ have the same sign}\label{sec:pqa2}

The procedure follows similarly to that of Section \ref{sec:pqab}; the  only difference is that  a
 different set of formulas is applied.
 Applying (\ref{connected sum for CJ}),
 (\ref{relation of J_T(n+2)}) and (\ref{(p,2,n+2)}), we have
\begin{equation}\label{equa:pq2a}\begin{aligned}
&[n+2]J_{\pq\#\aa}(n+2)
=J_{\pq}(n+2)J_{\aa}(n+2)
\\&
=\left(t^{-4pq(n+1)}J_{\pq}(n)+t^{-2pq(n+1)}\d(p,q,n)\right)\\
&\;\;\;\;\left(t^{-8a(n+1)}
J_{\aa}(n)-t^{-6a(n+1)}[2n+1]
+t^{-2a(n+1)}[2n+3]\right)
\\&=t^{-(4pq+8a)(n+1)}J_{\pq}(n)J_{\aa}(n)\\
&\;\;\;\;+\left(t^{-(4pq+2a)(n+1)}[2n+3]-t^{-(4pq+6a)(n+1)}[2n+1]\right)J_{\pq}(n)
\\&
\;\;\;\;+t^{-(2pq+8a)(n+1)}\d(p,q,n)J_{\aa}(n)
+\left(t^{-(2pq+2a)(n+1)}[2n+3]\right.\\&\;\;\;\;
\left.-t^{-(2pq+6a)(n+1)}[2n+1]\right)\d(p,q,n)
\\&=t^{-(4pq+8a)(n+1)}[n]J_{\pq\#\aa}(n)
+\left(t^{-(4pq+2a)(n+1)}[2n+3]\right.\\&\;\;\;\;\left.-t^{-(4pq+6a)(n+1)}[2n+1]\right)J_{\pq}(n)
+t^{-(2pq+8a)(n+1)}\d(p,q,n)J_{\aa}(n)
\\&\;\;\;\;+\left(t^{-(2pq+2a)(n+1)}[2n+3]-t^{-(2pq+6a)(n+1)}[2n+1]\right)\d(p,q,n)
\end{aligned}\end{equation}
from which we have
\begin{equation}\label{pq2a'}
\begin{aligned}
&\left([n+2]\sl^2-\sm^{-2(pq+2a)}t^{-(4pq+8a)}[n]\right)J_{\pq\#\aa}(n)
\\&
=\left(t^{-(4pq+2a)(n+1)}[2n+3]-t^{-(4pq+6a)(n+1)}[2n+1]\right)J_{\pq}(n)
\\&\;\;\;\;+t^{-(2pq+8a)(n+1)}\d(p,q,n)J_{\aa}(n)
\\&\;\;\;\;+\left(t^{-(2pq+2a)(n+1)}[2n+3]-t^{-(2pq+6a)(n+1)}[2n+1]\right)\d(p,q,n).
\end{aligned}
\end{equation}
The function
\begin{align*}
f(t,\sm)&=t^{-(4pq+2a)(n+1)}[2n+3]-t^{-(4pq+6a)(n+1)}[2n+1]
\\&=\sm^{-2pq-a}t^{-4pq-2a}\frac{\sm^{2}t^{6}-\sm^{-2}t^{-6}}{t^2-t^{-2}}
-\sm^{-2pq-3a}t^{-4pq-6a}\frac{\sm^{2}t^2-\sm^{-2}t^{-2}}{t^2-t^{-2}}
\end{align*}  is non-zero
since
\begin{equation}\label{f(-1,m)}\lim_{t\ra -1}(t^2-t^{-2})f(t,\sm)=
(\sm^2-\sm^{-2})(\sm^{-2pq-a}-\sm^{-2pq-3a})\ne 0.
\end{equation}
Applying the operator $(\sl^2-\sm^{-2pq}t^{-4pq})f^{-1}(t,\sm)$
from the left to both sides of  (\ref{pq2a'}) and then applying (\ref{an of J_T(n+2)}),
 (\ref{equ:operator}) (treating $\sm$ as $t^{2n}$),
(\ref{relation of J_T(n+2)}), and (\ref{(p,2,n+2)}),  we get
{\small
{\allowdisplaybreaks\begin{equation}\label{equ:pqa2}
\begin{aligned}
&\left(\sl^2-\sm^{-2pq}t^{-4pq}\right)f^{-1}(t,\sm)\left([n+2]\sl^2-\sm^{-2(pq+2a)}t^{-(4pq+8a)}[n]\right)J_{\pq\# \aa}(n)
\\&
=\left(\sl^2-\sm^{-2pq}t^{-4pq}\right)J_{\pq}(n)
+\left(\sl^2-\sm^{-2pq}t^{-4pq}\right)f^{-1}(t,\sm)t^{-(2pq+8a)(n+1)}\d(p,q,n)J_{\aa}(n)
\\&\;\;\;\;
+\left(\sl^2-\sm^{-2pq}t^{-4pq}\right)f^{-1}(t,\sm)\left(t^{-(2pq+2a)(n+1)}[2n+3]-t^{-(2pq+6a)(n+1)}[2n+1]\right)\d(p,q,n)
\\&
=\sm^{-pq}t^{-2pq}\d(p,q,n)
+f^{-1}(t,t^4\sm)t^{-(2pq+8a)(n+3)}\d(p,q,n+2)J_{\aa}(n+2)\\&\;\;\;\;
-\sm^{-2pq}t^{-4pq}f^{-1}(t,\sm)t^{-(2pq+8a)(n+1)}\d(p,q,n)J_{\aa}(n)
\\&\;\;\;\;+f^{-1}(t,t^4\sm)\left(t^{-(2pq+2a)(n+3)}[2n+7]-t^{-(2pq+6a)(n+3)}[2n+5]\right)\d(p,q,n+2)\\&\;\;\;\;
-\sm^{-2pq}t^{-4pq}f^{-1}(t,\sm)\left(t^{-(2pq+2a)(n+1)}[2n+3]-t^{-(2pq+6a)(n+1)}[2n+1]\right)\d(p,q,n)
\\&
=\sm^{-pq}t^{-2pq}\d(p,q,n)
+f^{-1}(t,t^4\sm)t^{-(2pq+8a)(n+3)}\d(p,q,n+2)
\\&\;\;\;\;\left(t^{-8a(n+1)}
J_{\aa}(n)-t^{-6a(n+1)}[2n+1]
+t^{-2a(n+1)}[2n+3]\right)\\&\;\;\;\;
-\sm^{-2pq}t^{-4pq}f^{-1}(t,\sm)t^{-(2pq+8a)(n+1)}\d(p,q,n)J_{\aa}(n)
\\&\;\;\;\;+f^{-1}(t,t^4\sm)\left(t^{-(2pq+2a)(n+3)}[2n+7]-t^{-(2pq+6a)(n+3)}[2n+5]\right)\d(p,q,n+2)\\&\;\;\;\;
-\sm^{-2pq}t^{-4pq}f^{-1}(t,\sm)\left(t^{-(2pq+2a)(n+1)}[2n+3]-t^{-(2pq+6a)(n+1)}[2n+1]\right)\d(p,q,n)
\\&
=\left(f^{-1}(t,t^4\sm)\sm^{-pq-8a}t^{-6pq-32a}\d(p,q,n+2)
-\sm^{-3pq-4a}t^{-6pq-8a}f^{-1}(t,\sm)\d(p,q,n)\right)J_{\aa}(n)\\&\;\;\;\;
+\sm^{-pq}t^{-2pq}\d(p,q,n)-f^{-1}(t,t^4\sm)
\sm^{-pq-7a}t^{-6pq-30a}\d(p,q,n+2)[2n+1]\\&\;\;\;\;
+f^{-1}(t,t^4\sm)
\sm^{-pq-5a}t^{-6pq-26a}\d(p,q,n+2)[2n+3]
\\&\;\;\;\;+f^{-1}(t,t^4\sm)\left(\sm^{-pq-a}t^{-6pq-6a}[2n+7]
-\sm^{-pq-3a}t^{-6pq-18a}[2n+5]\right)\d(p,q,n+2)\\&\;\;\;\;
-\sm^{-2pq}t^{-4pq}f^{-1}(t,\sm)\left(\sm^{-pq-a}t^{-2pq-2a}[2n+3]
-\sm^{-pq-3a}t^{-2pq-6a}[2n+1]\right)\d(p,q,n).
\end{aligned}
\end{equation}}}
The function
\[g(t, \sm)=f^{-1}(t,t^4\sm)\sm^{-pq-8a}t^{-6pq-32a}\d(p,q,n+2)
-\sm^{-3pq-4a}t^{-6pq-8a}f^{-1}(t,\sm)\d(p,q,n)\] is non-zero
since  from  (\ref{f(-1,m)}) and (\ref{equ:delta}) we see
\begin{equation}\label{g not zero} g(-1,\sm)=\frac{(\sm^{p+q}+\sm^{-p-q}-\sm^{q-p}-\sm^{p-q})
(\sm^{-pq-8a}-\sm^{-3pq-4a})}{(\sm^2-\sm^{-2})
(\sm^{-2pq-a}-\sm^{-2pq-3a})}\ne 0\end{equation}
because $pq\ne 2a$.

Therefore, we may apply
 the operator $(\sl+\sm^{-2a}t^{-2a})g^{-1}(t,\sm)$
from the left to  both sides of (\ref{equ:pqa2}),  and then apply (\ref{(p,2,n+1)}) and
 (\ref{equ:operator}), and we get
\begin{align*}
&\left(\sl+\sm^{-2a}t^{-2a}\right)g^{-1}(t,\sm)\left(\sl^2-\sm^{-2pq}t^{-4pq}\right)f^{-1}(t,\sm)\left([n+2]\sl^2-\sm^{-2(pq+2a)}t^{-(4pq+8a)}[n]\right)J_{\pq\#\aa}(n)\\
&=\left(\sl+\sm^{-2a}t^{-2a}\right)J_{\aa}(n)+
\left(\sl+\sm^{-2a}t^{-2a}\right)g^{-1}(t,\sm)\Big[\sm^{-pq}t^{-2pq}\d(p,q,n)
\\&\;\;\;\;-f^{-1}(t,t^4\sm)
\sm^{-pq-7a}t^{-6pq-30a}\d(p,q,n+2)[2n+1]
+f^{-1}(t,t^4\sm)
\sm^{-pq-5a}t^{-6pq-26a}\d(p,q,n+2)[2n+3]
\\&\;\;\;\;+f^{-1}(t,t^4\sm)\left(\sm^{-pq-a}t^{-6pq-6a}[2n+7]
-\sm^{-pq-3a}t^{-6pq-18a}[2n+5]\right)\d(p,q,n+2)\\&\;\;\;\;
-\sm^{-2pq}t^{-4pq}f^{-1}(t,\sm)\left(\sm^{-pq-a}t^{-2pq-2a}[2n+3]
-\sm^{-pq-3a}t^{-2pq-6a}[2n+1]\right)\d(p,q,n)\Big]\\
&
=\sm^{-a}[2n+1]+g^{-1}(t,t^2\sm)\Big[\sm^{-pq}t^{-4pq}\d(p,q,n+1)-f^{-1}(t,t^6\sm)
\sm^{-pq-7a}t^{-8pq-44a}\d(p,q,n+3)[2n+3]\\&\;\;\;\;
+f^{-1}(t,t^6\sm)
\sm^{-pq-5a}t^{-8pq-36a}\d(p,q,n+3)[2n+5]
+f^{-1}(t,t^6\sm)\left(\sm^{-pq-a}t^{-8pq-8a}[2n+9]\right.\\&\;\;\;\;
\left.-\sm^{-pq-3a}t^{-8pq-24a}[2n+7]\right)\d(p,q,n+3)
-\sm^{-2pq}t^{-8pq}f^{-1}(t,t^2\sm)\left(\sm^{-pq-a}t^{-4pq-4a}[2n+5]
\right.\\&\;\;\;\;\left.-\sm^{-pq-3a}t^{-4pq-12a}[2n+3]\right)\d(p,q,n+1)\Big]
+\sm^{-2a}t^{-2a}g^{-1}(t,\sm)\Big[\sm^{-pq}t^{-2pq}\d(p,q,n)
\\&\;\;\;\;-f^{-1}(t,t^4\sm)
\sm^{-pq-7a}t^{-6pq-30a}\d(p,q,n+2)[2n+1]
+f^{-1}(t,t^4\sm)
\sm^{-pq-5a}t^{-6pq-26a}\d(p,q,n+2)[2n+3]\\&\;\;\;\;
+f^{-1}(t,t^4\sm)\left(\sm^{-pq-a}t^{-6pq-6a}[2n+7]
-\sm^{-pq-3a}t^{-6pq-18a}[2n+5]\right)\d(p,q,n+2)\\&\;\;\;\;
-\sm^{-2pq}t^{-4pq}f^{-1}(t,\sm)\left(\sm^{-pq-a}t^{-2pq-2a}[2n+3]-
\sm^{-pq-3a}t^{-2pq-6a}[2n+1]\right)\d(p,q,n)\Big].
\end{align*}
Let $h(t,\sm)$ be the function on the right-hand side of the last equality.
Then $h(t, \sm)\ne 0$. Indeed, together with (\ref{g not zero}) and (\ref{f(-1,m)}),
we have
{\small\begin{equation}\label{h not zero}
\begin{aligned}
&\lim_{t\ra -1}(t^2-t^{-2})h(t,\sm)
\\&=\sm^{-a}(\sm^2-\sm^{-2})+\frac{\left(\sm^2-\sm^{-2}\right)
\left(\sm^{-2pq-a}-\sm^{-2pq-3a}\right)}{\left(\sm^{p+q}+\sm^{-p-q}-\sm^{q-p}-\sm^{p-q}\right)
\left(\sm^{-pq-8a}-\sm^{-3pq-4a}\right)}\\&\;\;\;\;
\left(\sm^{-pq}+
\frac{-\sm^{-pq-7a}+\sm^{-pq-5a}+\sm^{-pq-a}-\sm^{-pq-3a}
-\sm^{-3pq-a}+\sm^{-3pq-3a}}{\sm^{-2pq-a}-\sm^{-2pq-3a}}\right)\\
&\;\;\;\;\left(\sm^{p+q}+\sm^{-p-q}-\sm^{q-p}-\sm^{p-q}\right)\left(1+\sm^{-2a}\right)
\\&=\sm^{-a}(\sm^2-\sm^{-2})+\frac{(\sm^2-\sm^{-2})
 (-\sm^{-pq-7a}+\sm^{-pq-5a}+\sm^{-pq-a}-\sm^{-pq-3a}
)(1+\sm^{-2a})}{\sm^{-pq-8a}-\sm^{-3pq-4a}}\\
&=\frac{(\sm^2-\sm^{-2})(-\sm^{-3pq-5a}
+\sm^{-pq-a})}
{\sm^{-pq-8a}-\sm^{-3pq-4a}},
\end{aligned}\end{equation}}
which is a well-defined, non-zero function of $\sm$  since $pq\ne 2a$ and $pq\ne -2a$.
Hence,
\begin{align*}\a(t,\sm,\sl)&=
(\sl-1)\frac{1}{(t^2-t^{-2})h(t,\sm)}
\left(\sl+\sm^{-2a}t^{-2a}\right)g^{-1}(t,\sm)\left(\sl^2-\sm^{-2pq}t^{-4pq}\right)\\&
\;\;\;\;f^{-1}(t,\sm)\left([n+2]\sl^2-\sm^{-2(pq+2a)}t^{-(4pq+8a)}
[n]\right)
\end{align*}
is an annihilator of $J_{\pq\# T(a,2)}(n)$ of $\sl$-degree $6$.
Moreover using  (\ref{g not zero}), (\ref{h not zero}), (\ref{f(-1,m)}), (\ref{[n]}), and Lemma \ref{lem:A-poly of sum} (iv),
we have
\begin{align*}\a(-1,\sm,\sl)\;=&
 a(\sm) (\sl-1)\left(\sl+\sm^{-2a}\right)
\left(\sl^2-\sm^{-2pq}\right)\left(\sl^2-\sm^{-2(pq+2a)}\right)\\
= & a(\sm) \sm^{-4pq-6a}A_{\pq\# \aa}(\sm,\sl)
\end{align*}
where $$a(\sm)=\frac{1}{(-\sm^{-3pq-5a}
+\sm^{-pq-a})(\sm^{p+q}+\sm^{-p-q}-\sm^{q-p}-\sm^{p-q})
(\sm+\sm^{-1})}$$
is a well-defined, non-zero function in $\sm$.

We now proceed to show that $6$ is the minimum $\sl$-degree
among all non-zero annihilators of  $J_{\pq\#T(a,2)}(n)$.
So  suppose
$D_5\sl^5+D_4\sl^4 + D_3\sl^3 + D_2\sl^2 + D_1\sl + D_0$ is an annihilator of $J_{\pq\#T(a,2)}(n)$
with $D_0,\dots,D_5\in \q(t, \sm)$,   that is,
\begin{align*}
D_5J_{\pq\#\aa}(n+5)+D_4 J_{\pq\#\aa}(n+4)
 + D_3J_{\pq\#\aa}(n+3)\\ + D_2 J_{\pq\#\aa}(n+2)
+ D_1 J_{\pq\#\aa}(n+1)
 + D_0 J_{\pq\#\aa}(n) = 0.
\end{align*}
We wish to show that $D_i = 0$ for all $i$.

Applying  (\ref{equa:pq2a}),   (\ref{relation of J_T(n+2)}), (\ref{(p,2,n+2)})  and (\ref{(p,2)}), repeatedly, we have
\begin{align*}&
D_5J_{\pq\#\aa}(n+5)+D_4 J_{\pq\#\aa}(n+4)
 + D_3J_{\pq\#\aa}(n+3)\\& + D_2 J_{\pq\#\aa}(n+2)
+ D_1 J_{\pq\#\aa}(n+1)
 + D_0 J_{\pq\#\aa}(n)
 \\&=\frac{D_5}{[n+5]}
\left[ t^{-(4pq+8a)(n+4)}[n+3]J_{\pq\#\aa}(n+3)
+\left(t^{-(4pq+2a)(n+4)}[2n+9]\right.\right.\\
&\;\;\;\;\;\;\;\;\;\;\;\;\;\;\left.-t^{-(4pq+6a)(n+4)}[2n+7]\right)J_{\pq}(n+3)
+t^{-(2pq+8a)(n+4)}\d(p,q,n+3)J_{\aa}(n+3)
\\&
\;\;\;\;\;\;\;\;\;\;\;\;\;\;+\left.\left(t^{-(2pq+2a)(n+4)}[2n+9]-t^{-(2pq+6a)(n+4)}[2n+7]\right)\d(p,q,n+3)\right]
\\&\;\;+\frac{D_4}{[n+4]}
\left[ t^{-(4pq+8a)(n+3)}[n+2]J_{\pq\#\aa}(n+2)
+\left(t^{-(4pq+2a)(n+3)}[2n+7]\right.\right.\\
&\;\;\;\;\;\;\;\;\;\;\;\;\;\;\;\;\left.-t^{-(4pq+6a)(n+3)}[2n+5]\right)J_{\pq}(n+2)
+t^{-(2pq+8a)(n+3)}\d(p,q,n+2)J_{\aa}(n+2)
\\&
\;\;\;\;\;\;\;\;\;\;\;\;\;\;\;\;+\left.\left(t^{-(2pq+2a)(n+3)}[2n+7]-t^{-(2pq+6a)(n+3)}[2n+5]\right)\d(p,q,n+2)\right]\\
 &\;\;+ D_3J_{\pq\#\aa}(n+3)+ D_2 J_{\pq\#\aa}(n+2) + D_1 J_{\pq\#\aa}(n+1)
 + D_0 J_{\pq\#\aa}(n)\\
&=\left(\frac{D_5}{[n+5]}t^{-(4pq+8a)(n+4)}[n+3]
+D_3\right)J_{\pq\#\aa}(n+3)\\
&\;\;+\left(\frac{D_4}{[n+4]}t^{-(4pq+8a)(n+3)}[n+2]
+D_2\right)J_{\pq\#\aa}(n+2)\\
&\;\;+\frac{D_5}{[n+5]}
\left(t^{-(4pq+2a)(n+4)}[2n+9]
-t^{-(4pq+6a)(n+4)}[2n+7]\right)J_{\pq}(n+3)\\
&\;\;+\frac{D_5}{[n+5]}t^{-(2pq+8a)(n+4)}\d(p,q,n+3)J_{\aa}(n+3)\\&\;\;
+\frac{D_4}{[n+4]}\left(t^{-(4pq+2a)(n+3)}[2n+7]
-t^{-(4pq+6a)(n+3)}[2n+5]\right)J_{\pq}(n+2)\\&
\;\;+\frac{D_4}{[n+4]}t^{-(2pq+8a)(n+3)}\d(p,q,n+2)J_{\aa}(n+2)\\&\;\;
+\frac{D_5}{[n+5]}\left(t^{-(2pq+2a)(n+4)}[2n+9]-t^{-(2pq+6a)(n+4)}[2n+7]\right)\d(p,q,n+3)\\&\;\;
+\frac{D_4}{[n+4]}\left(t^{-(2pq+2a)(n+3)}[2n+7]-t^{-(2pq+6a)(n+3)}[2n+5]\right)\d(p,q,n+2)\\&\;\;
+ D_1 J_{\pq\#\aa}(n+1)
 + D_0 J_{\pq\#\aa}(n)\\&
=\left(\frac{D_5}{[n+5]}t^{-(4pq+8a)(n+4)}
+\frac{D_3}{[n+3]}\right)\left  [t^{-(4pq+8a)(n+2)}[n+1]J_{\pq\#\aa}(n+1)
\right.\\&\;\;\;\;\;+\left(t^{-(4pq+2a)(n+2)}[2n+5]-t^{-(4pq+6a)(n+2)}[2n+3]\right)J_{\pq}(n+1)
+t^{-(2pq+8a)(n+2)}\d(p,q,n+1)J_{\aa}(n+1)\\&\;\;\;\;\;
+\left(t^{-(2pq+2a)(n+2)}[2n+5]\right.- \left.\left.t^{-(2pq+6a)(n+2)}[2n+3]\right)\d(p,q,n+1)\right]
\\&\;\;
+\left(\frac{D_4}{[n+4]}t^{-(4pq+8a)(n+3)}
+\frac{D_2}{[n+2]}\right)\left [
t^{-(4pq+8a)(n+1)}[n]J_{\pq\#\aa}(n)\right.\\&\;\;\;\;\;\;
+
\left(t^{-(4pq+2a)(n+1)}[2n+3]-t^{-(4pq+6a)(n+1)}[2n+1]\right)J_{\pq}(n)\\&\;\;\;\;\;\;
+t^{-(2pq+8a)(n+1)}\d(p,q,n)J_{\aa}(n)
+\left(t^{-(2pq+2a)(n+1)}[2n+3]\right.-\left.\left.t^{-(2pq+6a)(n+1)}[2n+1]\right)\d(p,q,n)\right ]
\\&\;\;
+\frac{D_5}{[n+5]}
\left(t^{-(4pq+2a)(n+4)}[2n+9]
-t^{-(4pq+6a)(n+4)}[2n+7]\right)
 \left(t^{-4pq(n+2)} J_{\pq}(n+1)+t^{-2pq(n+2)} \delta(p,q,n+1)\right)
\\&
\;\;+\frac{D_5}{[n+5]}t^{-(2pq+8a)(n+4)}\d(p,q,n+3)
\left(t^{-8a(n+2)}J_{\aa}(n+1)-t^{-6a(n+2)}[2n+3]+t^{-2a(n+2)}[2n+5]\right)\\&\;\;
+\frac{D_4}{[n+4]}\left(t^{-(4pq+2a)(n+3)}[2n+7]
-t^{-(4pq+6a)(n+3)}[2n+5]\right)\left( t^{-4pq(n+1)} J_{\pq}(n)+t^{-2pq(n+1)} \delta(p,q,n)\right)
\\&\;\;
+\frac{D_4}{[n+4]}t^{-(2pq+8a)(n+3)}\d(p,q,n+2)\left(t^{-8a(n+1)}J_{\aa}(n)-t^{-6a(n+1)}[2n+1]+t^{-2a(n+1)}[2n+3]\right)\\&\;\;
+\frac{D_5}{[n+5]}\left(t^{-(2pq+2a)(n+4)}[2n+9]-t^{-(2pq+6a)(n+4)}[2n+7]\right)\d(p,q,n+3)\\&\;\;
+\frac{D_4}{[n+4]}\left(t^{-(2pq+2a)(n+3)}[2n+7]-t^{-(2pq+6a)(n+3)}[2n+5]\right)\d(p,q,n+2)\\&\;\;
+ D_1 J_{\pq\#\aa}(n+1)
 + D_0 J_{\pq\#\aa}(n)\\&
=\left[\left(\frac{D_5}{[n+5]}t^{-(4pq+8a)(n+4)}
+\frac{D_3}{[n+3]}\right)t^{-(4pq+8a)(n+2)}[n+1]+D_1\right]J_{\pq\#\aa}(n+1)\\&\;\;
+\left[\left(\frac{D_4}{[n+4]}t^{-(4pq+8a)(n+3)}
+\frac{D_2}{[n+2]}\right)
t^{-(4pq+8a)(n+1)}[n]+D_0\right]J_{\pq\#\aa}(n)
\\&\;\;
+\left[\left(\frac{D_5}{[n+5]}t^{-(4pq+8a)(n+4)}
+\frac{D_3}{[n+3]}\right)
\left(t^{-(4pq+2a)(n+2)}[2n+5]-t^{-(4pq+6a)(n+2)}[2n+3]\right)
\right.\\&\;\;\;\;\;\;
+\left.\frac{D_5}{[n+5]}
\left(t^{-(4pq+2a)(n+4)}[2n+9]
-t^{-(4pq+6a)(n+4)}[2n+7]\right)
t^{-4pq(n+2)}\right] J_{\pq}(n+1)\\&\;\;
+\left[\left(\frac{D_4}{[n+4]}t^{-(4pq+8a)(n+3)}
+\frac{D_2}{[n+2]}\right)
\left(t^{-(4pq+2a)(n+1)}[2n+3]-t^{-(4pq+6a)(n+1)}[2n+1]\right)\right.
\\&\;\;\;\;\;\;
+\left.\frac{D_4}{[n+4]}\left(t^{-(4pq+2a)(n+3)}[2n+7]
-t^{-(4pq+6a)(n+3)}[2n+5]\right) t^{-4pq(n+1)}\right]J_{\pq}(n)
\\&\;\;
+\left[\left(\frac{D_5}{[n+5]}t^{-(4pq+8a)(n+4)}
+\frac{D_3}{[n+3]}\right)t^{-(2pq+8a)(n+2)}\d(p,q,n+1)
\right.\\&\;\;\;\;\;\;
+\left.\frac{D_5}{[n+5]}t^{-(2pq+8a)(n+4)}\d(p,q,n+3)
t^{-8a(n+2)}\right]J_{\aa)}(n+1)\\&
\;\;+\left[\left(\frac{D_4}{[n+4]}t^{-(4pq+8a)(n+3)}
+\frac{D_2}{[n+2]}\right)t^{-(2pq+8a)(n+1)}\d(p,q,n)\right.\\&\;\;\;\;\;\;
+\left.\frac{D_4}{[n+4]}t^{-(2pq+8a)(n+3)}\d(p,q,n+2)t^{-8a(n+1)}\right]J_{\aa}(n)
\\&\;\;
+\left(\frac{D_5}{[n+5]}t^{-(4pq+8a)(n+4)}
+\frac{D_3}{[n+3]}\right)\left(t^{-(2pq+2a)(n+2)}[2n+5]-
t^{-(2pq+6a)(n+2)}[2n+3]\right)\d(p,q,n+1)\\&\;\;
+\left(\frac{D_4}{[n+4]}t^{-(4pq+8a)(n+3)}
+\frac{D_2}{[n+2]}\right)\left(t^{-(2pq+2a)(n+1)}[2n+3]-
t^{-(2pq+6a)(n+1)}[2n+1]\right)\d(p,q,n)
\\&\;\;
+\frac{D_5}{[n+5]}
\left(t^{-(4pq+2a)(n+4)}[2n+9]
-t^{-(4pq+6a)(n+4)}[2n+7]\right)
t^{-2pq(n+2)} \delta(p,q,n+1)
\\&
\;\;+\frac{D_5}{[n+5]}t^{-(2pq+8a)(n+4)}\d(p,q,n+3)
\left(-t^{-6a(n+2)}[2n+3]+t^{-2a(n+2)}[2n+5]\right)
\\&\;\;
+\frac{D_4}{[n+4]}\left(t^{-(4pq+2a)(n+3)}[2n+7]
-t^{-(4pq+6a)(n+3)}[2n+5]\right)
t^{-2pq(n+1)} \delta(p,q,n)\\
&\;\;+\frac{D_4}{[n+4]}t^{-(2pq+8a)(n+3)}\d(p,q,n+2)\left(-t^{-6a(n+1)}[2n+1]+t^{-2a(n+1)}[2n+3]\right)
\\&\;\;
+\frac{D_5}{[n+5]}\left(t^{-(2pq+2a)(n+4)}[2n+9]-t^{-(2pq+6a)(n+4)}[2n+7]\right)\d(p,q,n+3)\\&\;\;
+\frac{D_4}{[n+4]}\left(t^{-(2pq+2a)(n+3)}[2n+7]-t^{-(2pq+6a)(n+3)}[2n+5]\right)\d(p,q,n+2)
\\&
=\left[\left(\frac{D_5}{[n+5]}t^{-(4pq+8a)(n+4)}
+\frac{D_3}{[n+3]}\right)t^{-(4pq+8a)(n+2)}[n+1]+D_1\right]J_{\pq\#\aa}(n+1)\\&\;\;
+\left[\left(\frac{D_4}{[n+4]}t^{-(4pq+8a)(n+3)}
+\frac{D_2}{[n+2]}\right)
t^{-(4pq+8a)(n+1)}[n]+D_0\right]J_{\pq\#\aa}(n)
\\&\;\;
+\left[\left(\frac{D_5}{[n+5]}t^{-(4pq+8a)(n+4)}
+\frac{D_3}{[n+3]}\right)
\left(t^{-(4pq+2a)(n+2)}[2n+5]-t^{-(4pq+6a)(n+2)}[2n+3]\right)
\right.\\&\;\;\;\;\;\;
+\left.\frac{D_5}{[n+5]}
\left(t^{-(4pq+2a)(n+4)}[2n+9]
-t^{-(4pq+6a)(n+4)}[2n+7]\right)
t^{-4pq(n+2)}\right] J_{\pq}(n+1)\\&\;\;
+\left[\left(\frac{D_4}{[n+4]}t^{-(4pq+8a)(n+3)}
+\frac{D_2}{[n+2]}\right)
\left(t^{-(4pq+2a)(n+1)}[2n+3]-t^{-(4pq+6a)(n+1)}[2n+1]\right)\right.
\\&\;\;\;\;\;\;
+\left.\frac{D_4}{[n+4]}\left(t^{-(4pq+2a)(n+3)}[2n+7]
-t^{-(4pq+6a)(n+3)}[2n+5]\right) t^{-4pq(n+1)}\right]J_{\pq}(n)
\\&\;\;
+\left\{\left[\left(\frac{D_5}{[n+5]}t^{-(4pq+8a)(n+4)}
+\frac{D_3}{[n+3]}\right)t^{-(2pq+8a)(n+2)}\d(p,q,n+1)
\right.\right.\\&\;\;\;\;\;\;\;\;
+\left.\frac{D_5}{[n+5]}t^{-(2pq+8a)(n+4)}\d(p,q,n+3)
t^{-8a(n+2)}\right]\left(-t^{-2a(2n+1)}\right)\\&
\;\;\;\;\;\;+\left[\left(\frac{D_4}{[n+4]}t^{-(4pq+8a)(n+3)}
+\frac{D_2}{[n+2]}\right)t^{-(2pq+8a)(n+1)}\d(p,q,n)\right.\\&\;\;\;\;\;\;\;\;
+\left.\left.\frac{D_4}{[n+4]}t^{-(2pq+8a)(n+3)}\d(p,q,n+2)t^{-8a(n+1)}\right]\right\}J_{\aa}(n)
\\&\;\;
+\left[\left(\frac{D_5}{[n+5]}t^{-(4pq+8a)(n+4)}
+\frac{D_3}{[n+3]}\right)t^{-(2pq+8a)(n+2)}\d(p,q,n+1)
\right.\\&\;\;\;\;\;\;\;\;
+\left.\frac{D_5}{[n+5]}t^{-(2pq+8a)(n+4)}\d(p,q,n+3)
t^{-8a(n+2)}\right]t^{-2an}[2n+1]\\
&\;\;
+\left(\frac{D_5}{[n+5]}t^{-(4pq+8a)(n+4)}
+\frac{D_3}{[n+3]}\right)\left(t^{-(2pq+2a)(n+2)}[2n+5]-
t^{-(2pq+6a)(n+2)}[2n+3]\right)\d(p,q,n+1)\\&\;\;
+\left(\frac{D_4}{[n+4]}t^{-(4pq+8a)(n+3)}
+\frac{D_2}{[n+2]}\right)\left(t^{-(2pq+2a)(n+1)}[2n+3]-
t^{-(2pq+6a)(n+1)}[2n+1]\right)\d(p,q,n)
\\&\;\;
+\frac{D_5}{[n+5]}
\left(t^{-(4pq+2a)(n+4)}[2n+9]
-t^{-(4pq+6a)(n+4)}[2n+7]\right)
t^{-2pq(n+2)} \delta(p,q,n+1)
\\&
\;\;+\frac{D_5}{[n+5]}t^{-(2pq+8a)(n+4)}\d(p,q,n+3)
\left(-t^{-6a(n+2)}[2n+3]+t^{-2a(n+2)}[2n+5]\right)
\\&\;\;
+\frac{D_4}{[n+4]}\left(t^{-(4pq+2a)(n+3)}[2n+7]
-t^{-(4pq+6a)(n+3)}[2n+5]\right)
t^{-2pq(n+1)} \delta(p,q,n)
\\&\;\;
+\frac{D_4}{[n+4]}t^{-(2pq+8a)(n+3)}\d(p,q,n+2)\left(-t^{-6a(n+1)}[2n+1]+t^{-2a(n+1)}[2n+3]\right)
\\&\;\;
+\frac{D_5}{[n+5]}\left(t^{-(2pq+2a)(n+4)}[2n+9]-t^{-(2pq+6a)(n+4)}[2n+7]\right)\d(p,q,n+3)\\&\;\;
+\frac{D_4}{[n+4]}\left(t^{-(2pq+2a)(n+3)}[2n+7]-t^{-(2pq+6a)(n+3)}[2n+5]\right)\d(p,q,n+2)
\\&
=E_5J_{\pq\#\aa}(n+1)+E_4J_{\pq\#\aa}(n)
+E_3J_{\pq}(n+1)+E_2J_{\pq}(n)+E_1J_\aa(n)+E_0.\end{align*}

\begin{lemma}\label{lem:pqa2}The equation
\[E_5J_{\pq\#\aa}(n+1)+E_4J_{\pq\#\aa}(n)
+E_3J_{\pq}(n+1)+E_2J_{\pq}(n)+E_1J_\aa(n)+E_0=0\] implies
 $E_i=0$  for  all $i$.
\end{lemma}

The proof of this lemma is entirely similar to that of Lemma \ref{Ei are zero}.

 Now similarly as in  Section \ref{sec:pqab},  we show that the system of $6$  homogeneous  linear equations $E_i=0$,
$i=0,1,\dots,5$,
in $6$ variables $D_0, D_1,\dots, D_5$, has only the trivial solution: $D_i=0$  for all $i$.
\begin{align*}
& E_5=\left(\frac{D_5}{[n+5]}t^{-(4pq+8a)(n+4)}
+\frac{D_3}{[n+3]}\right)t^{-(4pq+8a)(n+2)}[n+1]+D_1=0
\\&
\lra\;\;\;\;\sm^{-4pq-8a}D_5+\sm^{-2pq-4a}D_3+D_1=0,
\\&
E_4=\left(\frac{D_4}{[n+4]}t^{-(4pq+8a)(n+3)}
+\frac{D_2}{[n+2]}\right)
t^{-(4pq+8a)(n+1)}[n]+D_0=0
\\&
\lra\;\;\;\;\sm^{-4pq-8a}D_4+\sm^{-2pq-4a}D_2+D_0=0,
\\&
E_3=\left(\frac{D_5}{[n+5]}t^{-(4pq+8a)(n+4)}
+\frac{D_3}{[n+3]}\right)
\left(t^{-(4pq+2a)(n+2)}[2n+5]-t^{-(4pq+6a)(n+2)}[2n+3]\right)
\\&\;\;\;\;\;\;\;\;
+\frac{D_5}{[n+5]}
\left(t^{-(4pq+2a)(n+4)}[2n+9]
-t^{-(4pq+6a)(n+4)}[2n+7]\right)
t^{-4pq(n+2)}=0
\\&
\lra\;\;\;\;\left(\sm^{-2pq-4a}+\sm^{-2pq}\right)D_5
+D_3=0,
\\&
 E_2=\left(\frac{D_4}{[n+4]}t^{-(4pq+8a)(n+3)}
+\frac{D_2}{[n+2]}\right)
\left(t^{-(4pq+2a)(n+1)}[2n+3]-t^{-(4pq+6a)(n+1)}[2n+1]\right)
\\&\;\;\;\;\;\;\;\;
+\frac{D_4}{[n+4]}\left(t^{-(4pq+2a)(n+3)}[2n+7]
-t^{-(4pq+6a)(n+3)}[2n+5]\right) t^{-4pq(n+1)}=0
\\&
\lra\;\;\;\;\left(\sm^{-2pq-4a}+\sm^{-2pq}\right)D_4
+D_2=0,
\\&
E_1=\left[\left(\frac{D_5}{[n+5]}t^{-(4pq+8a)(n+4)}
+\frac{D_3}{[n+3]}\right)t^{-(2pq+8a)(n+2)}\d(p,q,n+1)
\right.\\&\;\;\;\;\;\;\;\;
+\left.\frac{D_5}{[n+5]}t^{-(2pq+8a)(n+4)}\d(p,q,n+3)
t^{-8a(n+2)}\right]\left(-t^{-2a(2n+1)}\right)\\
&\;\;\;\;\;\;+\left[\left(\frac{D_4}{[n+4]}t^{-(4pq+8a)(n+3)}
+\frac{D_2}{[n+2]}\right)t^{-(2pq+8a)(n+1)}\d(p,q,n)\right.\\&\;\;\;\;\;\;\;\;
+\left.\frac{D_4}{[n+4]}t^{-(2pq+8a)(n+3)}\d(p,q,n+2)t^{-8a(n+1)}\right]=0
\\&
\lra\;\;\;\;\left(-\sm^{-2pq-6a}-\sm^{-6a}\right)D_5
+\left(\sm^{-2pq-4a}+\sm^{-4a}\right)D_4-\sm^{-2a}D_3+D_2=0,
\\&
E_0=\left[\left(\frac{D_5}{[n+5]}t^{-(4pq+8a)(n+4)}
+\frac{D_3}{[n+3]}\right)t^{-(2pq+8a)(n+2)}\d(p,q,n+1)
\right.\\&\;\;\;\;\;\;\;\;\;\;
+\left.\frac{D_5}{[n+5]}t^{-(2pq+8a)(n+4)}\d(p,q,n+3)
t^{-8a(n+2)}\right]t^{-2an}[2n+1]\\&
\;\;\;\;\;\;\;\;
+\left(\frac{D_5}{[n+5]}t^{-(4pq+8a)(n+4)}
+\frac{D_3}{[n+3]}\right)\left(t^{-(2pq+2a)(n+2)}[2n+5]-
t^{-(2pq+6a)(n+2)}[2n+3]\right)\d(p,q,n+1)\\&\;\;\;\;\;\;\;\;
+\left(\frac{D_4}{[n+4]}t^{-(4pq+8a)(n+3)}
+\frac{D_2}{[n+2]}\right)\left(t^{-(2pq+2a)(n+1)}[2n+3]-
t^{-(2pq+6a)(n+1)}[2n+1]\right)\d(p,q,n)
\\&\;\;\;\;\;\;\;\;
+\frac{D_5}{[n+5]}
\left(t^{-(4pq+2a)(n+4)}[2n+9]
-t^{-(4pq+6a)(n+4)}[2n+7]\right)
t^{-2pq(n+2)} \delta(p,q,n+1)
\\&
\;\;\;\;\;\;\;\;+\frac{D_5}{[n+5]}t^{-(2pq+8a)(n+4)}\d(p,q,n+3)
\left(-t^{-6a(n+2)}[2n+3]+t^{-2a(n+2)}[2n+5]\right)
\\&\;\;\;\;\;\;\;\;
+\frac{D_4}{[n+4]}\left(t^{-(4pq+2a)(n+3)}[2n+7]
-t^{-(4pq+6a)(n+3)}[2n+5]\right)
t^{-2pq(n+1)} \delta(p,q,n)
\\&\;\;\;\;\;\;\;\;
+\frac{D_4}{[n+4]}t^{-(2pq+8a)(n+3)}\d(p,q,n+2)\left(-t^{-6a(n+1)}[2n+1]+t^{-2a(n+1)}[2n+3]\right)
\\&\;\;\;\;\;\;\;\;
+\frac{D_5}{[n+5]}\left(t^{-(2pq+2a)(n+4)}[2n+9]-t^{-(2pq+6a)(n+4)}[2n+7]\right)\d(p,q,n+3)\\&\;\;\;\;\;\;\;\;
+\frac{D_4}{[n+4]}\left(t^{-(2pq+2a)(n+3)}[2n+7]-t^{-(2pq+6a)(n+3)}[2n+5]\right)\d(p,q,n+2)=0
\\&
\lra\;\;\;\;
 (\sm^{-3pq-9a}+\sm^{-pq-9a}+\sm^{-3pq-5a}-\sm^{-3pq-7a}
+\sm^{-3pq-a}
-\sm^{-3pq-3a}
 -\sm^{-pq-7a}+\sm^{-pq-5a}
\\&\;\;\;\;\;\;\;\;+\sm^{-pq-a}-\sm^{-pq-3a})D_5
+(\sm^{-3pq-5a}-\sm^{-3pq-7a}+\sm^{-3pq-a}
-\sm^{-3pq-3a}
-\sm^{-pq-7a}+\sm^{-pq-5a}\\&\;\;\;\;\;\;\;\;
+\sm^{-pq-a}-\sm^{-pq-3a})D_4
+(\sm^{-pq-5a}
+\sm^{-pq-a}-\sm^{-pq-3a})
D_3+(\sm^{-pq-a}-\sm^{-pq-3a})D_2
=0.
\end{align*}
That is, the system of linear equations is simplified to the
following one over $\z[\sm^{\pm 1}]$:
\begin{align*}&
\sm^{-4pq-8a}D_5+\sm^{-2pq-4a}D_3+D_1=0,
\\&
\sm^{-4pq-8a}D_4+\sm^{-2pq-4a}D_2+D_0=0,
\\&
\left(\sm^{-2pq-4a}+\sm^{-2pq}\right)D_5
+D_3=0,
\\&
\left(\sm^{-2pq-4a}+\sm^{-2pq}\right)D_4
+D_2=0,
\\&
\left(-\sm^{-2pq-6a}-\sm^{-6a}\right)D_5
+\left(\sm^{-2pq-4a}+\sm^{-4a}\right)D_4-\sm^{-2a}D_3+D_2=0,
\\&
(\sm^{-3pq-9a}+\sm^{-pq-9a}+\sm^{-3pq-5a}-\sm^{-3pq-7a}
+\sm^{-3pq-a}
-\sm^{-3pq-3a}
 -\sm^{-pq-7a}+\sm^{-pq-5a}
\\&+\sm^{-pq-a}-\sm^{-pq-3a})D_5
+(\sm^{-3pq-5a}-\sm^{-3pq-7a}+\sm^{-3pq-a}
-\sm^{-3pq-3a}
-\sm^{-pq-7a}+\sm^{-pq-5a}\\
&+\sm^{-pq-a}-\sm^{-pq-3a})D_4
+(\sm^{-pq-5a}
+\sm^{-pq-a}-\sm^{-pq-3a})
D_3+(\sm^{-pq-a}-\sm^{-pq-3a})D_2
=0.
\end{align*}

From the third and fourth equations, we get
$D_3=-\left(\sm^{-2pq-4a}+\sm^{-2pq}\right)D_5$
and $D_2=-\left(\sm^{-2pq-4a}+\sm^{-2pq}\right)D_4$
Plugging them into the fifth and sixth equations, we get
\begin{align*}&
\left(-\sm^{-6a}
+\sm^{-2pq-2a}\right)D_5
+\left(\sm^{-4a}-\sm^{-2pq}\right)D_4
=0,
\\&
 (\sm^{-pq-9a}
 -\sm^{-pq-7a}+\sm^{-pq-5a}
-\sm^{-pq-3a}+\sm^{-pq-a}
-\sm^{-3pq-5a})D_5\\
&+(-\sm^{-pq-7a}+\sm^{-pq-5a}
+\sm^{-pq-a}-\sm^{-pq-3a}
)D_4
=0.
\end{align*}
The determinant of these two equations
is
\begin{align*}&
\left(-\sm^{-6a}
+\sm^{-2pq-2a}\right)(-\sm^{-pq-7a}+\sm^{-pq-5a}
+\sm^{-pq-a}-\sm^{-pq-3a}
)\\&-\left(\sm^{-4a}-\sm^{-2pq}\right) (\sm^{-pq-9a}
 -\sm^{-pq-7a}+\sm^{-pq-5a}
-\sm^{-pq-3a}+\sm^{-pq-a}
-\sm^{-3pq-5a})\\
&
=\sm^{-pq-13a}-\sm^{-pq-11a}
-\sm^{-pq-7a}+\sm^{-pq-9a}-\sm^{-3pq-9a}+\sm^{-3pq-7a}
+\sm^{-3pq-3a}-\sm^{-3pq-5a}
\\&
-\sm^{-pq-13a}
 +\sm^{-pq-11a}-\sm^{-pq-9a}
+\sm^{-pq-7a}-\sm^{-pq-5a}+\sm^{-3pq-9a}\\
&
+\sm^{-3pq-9a}
 -\sm^{-3pq-7a}+\sm^{-3pq-5a}
-\sm^{-3pq-3a}+\sm^{-3pq-a}
-\sm^{-5pq-5a}\\&=
-\sm^{-pq-5a}
+\sm^{-3pq-9a}
+\sm^{-3pq-a}
-\sm^{-5pq-5a}
=(\sm^{-4a}-\sm^{-2pq})(\sm^{-3pq-5a}-\sm^{-pq-a})\end{align*}
which is non-zero since $pq\ne 2a$ and $pq\ne -2a$.
Therefore, $D_5=D_4=0$,  and thus $D_3=D_2=0$.
Then by the first two equations we have $D_1=D_0=0$.

\section{Case $|p|>q>2$, $|a|>b=2$, $pq=2a$}

This section is handled in a similar way to  Section \ref{sec:pq=ab}.
Formula (\ref{equ:pqa2}) in Section \ref{sec:pqa2} remains valid when $pq=2a$, yielding
\begin{equation}\label{equ:pq=2a}
\begin{aligned}
&\left(\sl^2-\sm^{-4a}t^{-8a}\right)f^{-1}(t,\sm)\left([n+2]\sl^2-\sm^{-8a}t^{-16a}[n]\right)J_{\pq\# \aa}(n)
\\&=\left(f^{-1}(t,t^4\sm)\sm^{-10a}t^{-44a}\d(p,q, n+2)
-\sm^{-10a}t^{-20a}f^{-1}(t,\sm)\d(p,q,n)\right)J_{\aa}(n)\\&\;\;\;\;
+\sm^{-2a}t^{-4a}\d(p,q,n)-f^{-1}(t,t^4\sm)
\sm^{-9a}t^{-42a}\d(p,q,n+2)[2n+1]\\&\;\;\;\;
+f^{-1}(t,t^4\sm)
\sm^{-7a}t^{-38a}\d(p,q,n+2)[2n+3]
\\&\;\;\;\;+f^{-1}(t,t^4\sm)\left(\sm^{-3a}t^{-18a}[2n+7]-\sm^{-5a}t^{-30a}[2n+5]\right)\d(p,q,n+2)\\&\;\;\;\;
-\sm^{-4a}t^{-8a}f^{-1}(t,\sm)\left(\sm^{-3a}t^{-6a}[2n+3]-\sm^{-5a}t^{-10a}[2n+1]\right)\d(p,q,n)
\end{aligned}\end{equation}
where
\begin{align*}&f(t,\sm)=\sm^{-5a}t^{-10a}\frac{\sm^2t^6-\sm^{-2}t^{-6}}{t^2-t^{-2}}
-\sm^{-7a}t^{-14a}\frac{\sm^2t^2-\sm^{-2}t^{-2}}{t^2-t^{-2}}
\\&\;\;\;\;\;\;\;\;\;\;\;\;=\frac{\sm^{-5a+2}t^{-10a+6}-\sm^{-5a-2}t^{-10a-6}
-\sm^{-7a+2}t^{-14a+2}+\sm^{-7a-2}t^{-14a-2}}{t^2-t^{-2}}\end{align*}
which is non-zero since
\[ \lim_{t\ra -1} (t^2-t^{-2})f(t,\sm)=
(\sm^2-\sm^{-2})(\sm^{-5a}-\sm^{-7a})\ne 0.\]
The function $\displaystyle g(t, \sm)=
f^{-1}(t,t^4\sm)\sm^{-10a}t^{-44a}\d(p,q, n+2)
-\sm^{-10a}t^{-20a}f^{-1}(t,\sm)\d(p,q,n)$, which is the coefficient
of $J_{\aa}(n)$, is non-zero (but $g(-1,\sm)=0$, cf. (\ref{g not zero})).
Indeed
\begin{align*}&
g(t,\sm)\\&=
\frac{\sm^{-10a}t^{-44a}\left(\sm^{p+q}t^{6(p+q)+2}+\sm^{-p-q}t^{-6(p+q)+2}
-\sm^{q-p}t^{6(q-p)-2}-\sm^{p-q}t^{-6(q-p)-2}\right)}{\sm^{-5a+2}t^{-30a+14}-\sm^{-5a-2}t^{-30a-14}
-\sm^{-7a+2}t^{-42a+10}+\sm^{-7a-2}t^{-42a-10}}\\
&\;\;-\frac{\sm^{-10a}t^{-20a}\left(\sm^{p+q}t^{2(p+q)+2}+\sm^{-p-q}t^{-2(p+q)+2}
-\sm^{q-p}t^{2(q-p)-2}-\sm^{p-q}t^{-2(q-p)-2}\right)}{\sm^{-5a+2}t^{-10a+6}-\sm^{-5a-2}t^{-10a-6}
-\sm^{-7a+2}t^{-14a+2}+\sm^{-7a-2}t^{-14a-2}}
\end{align*}
from which we can express $g(t,\sm)$ as the fraction
\begin{equation}
\label{g(t,m)}
g(t,\sm)=\frac{\phi(t,\sm)}{\psi(t,\sm)},
\end{equation}
where
\begin{align*}
\phi(t,\sm)&
=\sm^{-10a}t^{-44a}\left(\sm^{p+q}t^{6(p+q)+2}+\sm^{-p-q}t^{-6(p+q)+2}
-\sm^{q-p}t^{6(q-p)-2}-\sm^{p-q}t^{-6(q-p)-2}\right)\\
&\;\;\;\;\;\left(\sm^{-5a+2}t^{-10a+6}-\sm^{-5a-2}t^{-10a-6}
-\sm^{-7a+2}t^{-14a+2}+\sm^{-7a-2}t^{-14a-2}\right)\\
&\;\;-\sm^{-10a}t^{-20a}\left(\sm^{p+q}t^{2(p+q)+2}+\sm^{-p-q}t^{-2(p+q)+2}
-\sm^{q-p}t^{2(q-p)-2}-\sm^{p-q}t^{-2(q-p)-2}\right)\\
&\;\;\;\;\;\;\left(\sm^{-5a+2}t^{-30a+14}-\sm^{-5a-2}t^{-30a-14}
-\sm^{-7a+2}t^{-42a+10}+\sm^{-7a-2}t^{-42a-10}\right)
\end{align*}
and
\begin{align*}
\psi(t,\sm)
&=\left(\sm^{-5a+2}t^{-30a+14}-\sm^{-5a-2}t^{-30a-14}
-\sm^{-7a+2}t^{-42a+10}+\sm^{-7a-2}t^{-42a-10}\right)\\
&\;\;\;\;\;\left(\sm^{-5a+2}t^{-10a+6}-\sm^{-5a-2}t^{-10a-6}
-\sm^{-7a+2}t^{-14a+2}+\sm^{-7a-2}t^{-14a-2}\right)
\end{align*}
Note that $\phi(-1,\sm)=0$ and $\psi(-1,\sm)\ne 0$.
So to show that $g(t,\sm)\ne 0$, we just need to show that $\phi(t,\sm)\ne 0$.
This can be demonstrated by considering the degree of the Laurent  polynomial
$\phi(t, t^{2n})$.

When $a>0$ (and thus $a>p>q>2$),  the first  product term of $\phi(t,t^{2n})$
 has the lowest degree:  $-20an-44a-2(p+q)n-6(p+q)+2-14an-4n-14a-2=-34an-2(p+q)n-4n-58a-6(p+q)$.
The second product term has the  lowest degree: $-20an-20a-2(p+q)n-2(p+q)+2-14an-4n-42a-10=-34an-2(p+q)n-4n-62a-2(p+q)-8$.
The difference between these degrees is:
\[-4(p+q)+4a+8=4(a+2-p-q)=4(a+2-p-\frac{2a}{p})=\frac{4(ap+2p-p^2-2a)}{p}
=\frac{4(a-p)(p-2)}{p}>0.\]
Similarly, when $a<0$ (and thus $a<p<-q<-2$),
the first product term of $\phi(t,t^{2n})$
 has the highest degree:  $-20an-44a+2(q-p)n+6(q-p)-2-14an+4n-14a+2=
-34an+2(q-p)n+4n-58a+6(q-p)$. The second product term has the  highest degree: $-20an-20a+2(q-p)n+2(q-p)-2-14an+4n-42a+10=-34an+2(q-p)n+4n-62a+2(q-p)+8$.
The difference is:
\[
4(q-p)+4a-8=4(q-p+a-2)=4(\frac{2a}{p}-p+a-2)=\frac{4(2a-p^2+ap-2p)}{p}
=\frac{4(a-p)(p+2)}{p}<0.\]
Hence, $\phi(t,\sm)\ne 0$, which implies  $g(t,\sm)\ne 0$.

Therefore, we may continue to apply the operator  $\left(\sl+\sm^{-2a}t^{-2a}\right)g^{-1}(t,\sm)$
from the left to (\ref{equ:pq=2a}), yielding
\begin{align*}
&\left(\sl+\sm^{-2a}t^{-2a}\right)g^{-1}(t,\sm)\left(\sl^2-\sm^{-4a}t^{-8a}\right)f^{-1}(t,\sm)\left([n+2]\sl^2-\sm^{-8a}t^{-16a}[n]\right)J_{\pq\#\aa}(n)\\
&=\sm^{-a}[2n+1]+g^{-1}(t,t^2\sm)\left[\sm^{-2a}t^{-8a}\d(p,q,n+1)-f^{-1}(t,t^6\sm)
\sm^{-9a}t^{-60a}\d(p,q,n+3)[2n+3]\right.\\&\;\;\;\;
+f^{-1}(t,t^6\sm)
\sm^{-7a}t^{-52a}\d(p,q,n+3)[2n+5]
+f^{-1}(t,t^6\sm)\left(\sm^{-3a}t^{-24a}[2n+9]\right.\\&\;\;\;\;
\left.-\sm^{-5a}t^{-40a}[2n+7]\right)\d(p,q,n+3)
-\sm^{-4a}t^{-16a}f^{-1}(t,t^2\sm)\left(\sm^{-3a}t^{-12a}[2n+5]
\right.\\&\;\;\;\;\left.\left.-\sm^{-5a}t^{-20a}[2n+3]\right)\d(p,q,n+1)\right]
+\sm^{-2a}t^{-2a}g^{-1}(t,\sm)\left[\sm^{-2a}t^{-4a}\d(p,q,n)\right.
\\&\;\;\;\;-f^{-1}(t,t^4\sm)
\sm^{-9a}t^{-42a}\d(p,q,n+2)[2n+1]
+f^{-1}(t,t^4\sm)
\sm^{-7a}t^{-38a}\d(p,q,n+2)[2n+3]\\&\;\;\;\;
+f^{-1}(t,t^4\sm)\left(\sm^{-3a}t^{-18a}[2n+7]
-\sm^{-5a}t^{-30a}[2n+5]\right)\d(p,q,n+2)\\&\;\;\;\;
-\sm^{-4a}t^{-8a}f^{-1}(t,\sm)\left(\sm^{-3a}t^{-6a}[2n+3]-
\left.\sm^{-5a}t^{-10a}[2n+1]\right)\d(p,q,n)\right].
\end{align*}
Let $h(t,\sm)$ be the function on the right-hand side of the last equality.
Then $h(t,\sm)\ne 0$
since
$\displaystyle\lim_{t\ra -1}(t^2-t^{-2})h(t,\sm)
=\infty$
(cf. (\ref{h not zero})).
Therefore,
\begin{align*}&\a(t,\sm,\sl)\\&=(\sl-1)\frac{h^{-1}(t,\sm)}{t^2-t^{-2}}\left(\sl+\sm^{-2a}t^{-2a}\right)g^{-1}(t,\sm)\left(\sl^2-\sm^{-4a}t^{-8a}\right)f^{-1}(t,\sm)\left([n+2]\sl^2-\sm^{-8a}t^{-16a}[n]\right)
\\&=(\sl-1)\frac{h^{-1}(t,\sm)g^{-1}(t,t^2\sm)}{t^2-t^{-2}}\left(\sl+\frac{g^{-1}(t,\sm)}{g^{-1}(t,t^2\sm)}\sm^{-2a}t^{-2a}\right)\left(\sl^2-\sm^{-4a}t^{-8a}\right)\\
&\;\;\;\; f^{-1}(t,\sm)\left([n+2]\sl^2-\sm^{-8a}t^{-16a}[n]\right)
\end{align*}
is an annihilator of $J_{\pq\#\aa}(n)$. To evaluate $\a(t,\sm,\sl)$ at $t=-1$, we need the following lemma.

\begin{lemma}\label{pphi}$\displaystyle\lim_{t\ra -1}\frac{g^{-1}(t,\sm)}{g^{-1}(t,t^2\sm)}=1.$
\end{lemma}

\pf The proof is similar to that of Lemma \ref{lem:limit=1}.
By  (\ref{g(t,m)}) and the note following it,  we have
\begin{align*}\lim_{t\ra -1}\frac{g^{-1}(t,\sm)}{g^{-1}(t, t^2\sm)}
=\lim_{t\ra -1}\frac{\psi(t,\sm)\phi(t, t^2\sm)}{\phi(t, \sm)\psi(t,t^2\sm)}
=\lim_{t\ra -1}\frac{\psi(t,\sm)}{\psi(t,t^2\sm)}\lim_{t\ra-1}
\frac{\phi(t, t^2\sm)}{\phi(t, \sm)}=\lim_{t\ra-1}
\frac{\phi(t, t^2\sm)}{\phi(t, \sm)}.
\end{align*}
Applying   L'H\^opital's rule (treating $\sm$ as an independent variable), we have
 \[ \lim_{t\ra -1}\frac{\phi(t, t^2\sm)}{\phi(t,\sm)}
= \lim_{t\ra -1}\frac{\dfrac{\p\phi(t, t^2\sm)}{\p t}}{\dfrac{\p \phi(t,\sm)}{\p t}}.\]
Now  one can check that \[\lim_{t\ra -1}\frac{\p \phi(t,\sm)}{\p t}=\lim_{t\ra -1}\frac{\p \phi(t,t^2\sm)}{\p t}\ne 0.\]
Hence the lemma is true. \qed

Applying Lemma \ref{pphi} and the fact that $g(-1,\sm)=0$, we can evaluate $\displaystyle \lim_{t\ra -1}\frac{h^{-1}(t,\sm)g^{-1}(t,t^2\sm)}{t^2-t^{-2}}$.
\begin{align*}
& \lim_{t\ra-1}(t^2-t^{-2})g(t,t^2\sm)h(t,\sm)\\
&=\lim_{t\ra-1}
\left\{(t^2-t^{-2})g(t, t^2\sm)\sm^{-a}[2n+1]+(t^2-t^{-2})\Big[\sm^{-2a}t^{-8a}\d(p,q,n+1)\right.\\&\;\;\;\;
-f^{-1}(t,t^6\sm)
\sm^{-9a}t^{-60a}\d(p,q,n+3)[2n+3]
+f^{-1}(t,t^6\sm)
\sm^{-7a}t^{-52a}\d(p,q,n+3)[2n+5]\\&\;\;\;\;
+f^{-1}(t,t^6\sm)\left(\sm^{-3a}t^{-24a}[2n+9]
-\sm^{-5a}t^{-40a}[2n+7]\right)\d(p,q,n+3)
\\&\;\;\;\;-\sm^{-4a}t^{-16a}f^{-1}(t,t^2\sm)\left(\sm^{-3a}t^{-12a}[2n+5]
-\sm^{-5a}t^{-20a}[2n+3]\right)\d(p,q,n+1)\Big]
\\&\;\;\;\;+\sm^{-2a}t^{-2a}(t^2-t^{-2})\frac{g(t, t^{2}\sm)}{g(t,\sm)}\Big[\sm^{-2a}t^{-4a}\d(p,q,n)
\\&\;\;\;\;-f^{-1}(t,t^4\sm)
\sm^{-9a}t^{-42a}\d(p,q,n+2)[2n+1]
+f^{-1}(t,t^4\sm)
\sm^{-7a}t^{-38a}\d(p,q,n+2)[2n+3]\\&\;\;\;\;
+f^{-1}(t,t^4\sm)\left(\sm^{-3a}t^{-18a}[2n+7]
-\sm^{-5a}t^{-30a}[2n+5]\right)\d(p,q,n+2)\\&\;\;\;\;
-\sm^{-4a}t^{-8a}f^{-1}(t,\sm)\left(\sm^{-3a}t^{-6a}[2n+3]-
\left.\sm^{-5a}t^{-10a}[2n+1]\right)\d(p,q,n)\Big]\right\}
\\
&=\lim_{t\ra-1}
\left\{(t^2-t^{-2})\Big[\sm^{-2a}t^{-8a}\d(p,q,n+1)-f^{-1}(t,t^6\sm)
\sm^{-9a}t^{-60a}\d(p,q,n+3)[2n+3]\right.\\&\;\;\;\;
+f^{-1}(t,t^6\sm)
\sm^{-7a}t^{-52a}\d(p,q,n+3)[2n+5]
+f^{-1}(t,t^6\sm)\left(\sm^{-3a}t^{-24a}[2n+9]\right.\\&\;\;\;\;
\left.-\sm^{-5a}t^{-40a}[2n+7]\right)\d(p,q,n+3)
-\sm^{-4a}t^{-16a}f^{-1}(t,t^2\sm)\left(\sm^{-3a}t^{-12a}[2n+5]
\right.\\&\;\;\;\;\left.-\sm^{-5a}t^{-20a}[2n+3]\right)\d(p,q,n+1)\Big]
+\sm^{-2a}t^{-2a}(t^2-t^{-2})\Big[\sm^{-2a}t^{-4a}\d(p,q,n)
\\&\;\;\;\;-f^{-1}(t,t^4\sm)
\sm^{-9a}t^{-42a}\d(p,q,n+2)[2n+1]
+f^{-1}(t,t^4\sm)
\sm^{-7a}t^{-38a}\d(p,q,n+2)[2n+3]\\&\;\;\;\;
+f^{-1}(t,t^4\sm)\left(\sm^{-3a}t^{-18a}[2n+7]
-\sm^{-5a}t^{-30a}[2n+5]\right)\d(p,q,n+2)\\&\;\;\;\;
-\sm^{-4a}t^{-8a}f^{-1}(t,\sm)\left(\sm^{-3a}t^{-6a}[2n+3]-
\left.\sm^{-5a}t^{-10a}[2n+1]\right)\d(p,q,n)\Big]\right\}
\\
&=
(1+\sm^{-2a})\left(\sm^{-2a}+\frac{\sm^{-3a}-\sm^{-5a}}
{\sm^{-5a}-\sm^{-7a}}\right)(\sm^{p+q}+\sm^{-p-q}-\sm^{q-p}-\sm^{p-q})\\
&=\frac{(1+\sm^{-2a})(\sm^{-7a}-\sm^{-9a} +\sm^{-3a}-\sm^{-5a})(\sm^{p+q}+\sm^{-p-q}-\sm^{q-p}-\sm^{p-q})
}{\sm^{-5a}-\sm^{-7a}}\end{align*}
Therefore, \[\a(-1,\sm, \sl)=a(\sm)(\sl-1)(\sl+\sm^{-2a})(\sl^2-\sm^{-4a})(\sl^2-\sm^{-8a})\]
where \[a(\sm)=\frac{1}{(1+\sm^{-2a})(\sm^{-7a}-\sm^{-9a} +\sm^{-3a}-\sm^{-5a})(\sm^{p+q}+\sm^{-p-q}-\sm^{q-p}-\sm^{p-q})(\sm+\sm^{-1})}\] is a well-defined, non-zero function.
Visibly $\a(-1,\sm,\sl)$ contains repeated factor  $\sl+\sm^{-2a}$, and after removing this factor,
 is equal to the $A$-polynomial of $T(p,q)\# T(a,2)$ (where $pq=2a$) up to a non-zero factor
only depending  on $\sm$.

We now show that $6$ is the least $\sl$-degree among all non-zero annihilators of
$J_{\pq\#\aa}(n)$.
So suppose \begin{align*}&D_5J_{T(p,q)\#T(a,2)}(n+5)+D_4J_{T(p,q)\#T(a,2)}(n+4)
+D_3J_{T(p,q)\#T(a,2)}(n+3)\\&+D_2J_{T(p,q)\#T(a,2)}(n+2)+D_1J_{T(p,q)\#T(a,2)}(n+1)+D_0J_{T(p,q)\#T(a,2)}(n)
=0\end{align*}
for some $D_i\in\q(t,\sm)$, $i=0,1, \dots, 5$, we want to show that all $D_i$ have to be zero.

Note that  in Section  \ref{sec:pqa2} the calculation of
\begin{align*}&D_5J_{T(p,q)\#T(a,2)}(n+5)+D_4J_{T(p,q)\#T(a,2)}(n+4)
+D_3J_{T(p,q)\#T(a,2)}(n+3)\\&+D_2J_{T(p,q)\#T(a,2)}(n+2)+D_1J_{T(p,q)\#T(a,2)}(n+1)+D_0J_{T(p,q)\#T(a,2)}(n)
\\&=E_5J_{T(p,q)\#T(a,2)}(n+1)+E_4J_{T(p,q)\#T(a,2)}(n)+E_3J_{(p,q)}(n+1)
+E_2J_{T(p,q)}(n)+E_1J_{T(a,2)}(n)+E_0\end{align*}
is still valid when $pq=2a$. By  Lemma \ref{lem:pqa2},
$E_i=0$ for all $i=0,1,\dots, 5$.
But when $pq=2a$, the argument
for the statement that $E_i=0$ for all $i$ implies $D_i=0$ for all $i$
needs some modification.

Plugging $pq=2a$ into the system $E_i=0$, $i=0,1,\dots, 5$,  obtained in  Section  \ref{sec:pqa2},
and simplifying    $E_2=0$, $E_1=0$ and $E_0=0$ as in Section  \ref{sec:pqa2},
we get
\begin{align*}
&E_5=(\frac{D_5}{[n+5]}t^{-16a(n+4)}+\frac{D_3}{[n+3]})t^{-16(n+2)}[n+1]+D_1=0,
\\&
 E_4=\left(\frac{D_4}{[n+4]}t^{-16a(n+3)}
+\frac{D_2}{[n+2]}\right)
t^{-16a(n+1)}[n]+D_0=0,
\\&
E_3=\left(\frac{D_5}{[n+5]}t^{-16a(n+4)}
+\frac{D_3}{[n+3]}\right)
\left(t^{-10a(n+2)}[2n+5]-t^{-14a(n+2)}[2n+3]\right)
\\&\;\;\;\;\;\;\;\;
+\frac{D_5}{[n+5]}
\left(t^{-10a(n+4)}[2n+9]
-t^{-14a(n+4)}[2n+7]\right)
t^{-8a(n+2)}\\
&\;\;\;\;=
\frac{D_5}{[n+5]}
\left(t^{-18an-56a}[2n+9]
-t^{-22an-72a}[2n+7]+t^{-26an-84a}[2n+5]
-t^{-30an-92a}[2n+3]\right)
\\&\;\;\;\;\;\;\;\;+\frac{D_3}{[n+3]}
\left(t^{-10an-20a}[2n+5]-t^{-14an-28a}[2n+3]\right)=0,
\\&
E_2=0\lra
\left(\sm^{-8a}+\sm^{-4a}\right)D_4
+D_2=0,
\\&
E_1=0\lra \left(-\sm^{-10a}-\sm^{-6a}\right)D_5
+\left(\sm^{-8a}+\sm^{-4a}\right)D_4-\sm^{-2a}D_3+D_2=0,
\\&
E_0=0
\lra
(\sm^{-15a}+2\sm^{-11a}-\sm^{-13a}
+2\sm^{-7a}
-2\sm^{-9a}
+\sm^{-3a}-\sm^{-5a})D_5
+(\sm^{-11a}-\sm^{-13a}
+2\sm^{-7a}\\&\;\;\;\;\;\;\;\;\;\;\;\;\;\;\;\;-2\sm^{-9a}
+\sm^{-3a}-\sm^{-5a})D_4
+(\sm^{-7a}
+\sm^{-3a}-\sm^{-5a})
D_3+(\sm^{-3a}-\sm^{-5a})D_2
=0.
\end{align*}
 From $E_3=0$, we get
\[
 D_3=
-\frac{[n+3]\left(t^{-8an-36a}[2n+9]
-t^{-12an-52a}[2n+7]+t^{-16an-64a}[2n+5]
-t^{-20an-72a}[2n+3]\right)D_5}{[n+5]\left([2n+5]-t^{-4an-8a}[2n+3]\right)}.
\]
From the simplified version of $E_2=0$, we have
\[
D_2
=-\left(\sm^{-8a}+\sm^{-4a}\right)D_4.
\]
Plugging $D_3$ and $D_2$ into the simplified $E_1=0$ and $E_0=0$,
we have
\begin{align*}&
\left(-\sm^{-10a}-\sm^{-6a}\right)D_5
+\left(\sm^{-8a}+\sm^{-4a}\right)D_4
\\&
+\sm^{-2a}\frac{[n+3]\left(t^{-8an-36a}[2n+9]
-t^{-12an-52a}[2n+7]+t^{-16an-64a}[2n+5]
-t^{-20an-72a}[2n+3]\right)D_5}{[n+5]\left([2n+5]-t^{-4an-8a}[2n+3]\right)}\\&
-\left(\sm^{-8a}+\sm^{-4a}\right)D_4=0,
\\\text{and}
\\&
 (\sm^{-15a}+2\sm^{-11a}-\sm^{-13a}
+2\sm^{-7a}
-2\sm^{-9a}
+\sm^{-3a}-\sm^{-5a})D_5
+(\sm^{-11a}-\sm^{-13a}
+2\sm^{-7a}\\&-2\sm^{-9a}
+\sm^{-3a}-\sm^{-5a})D_4
-(\sm^{-7a}
+\sm^{-3a}-\sm^{-5a})
\\&
\frac{[n+3]\left(t^{-8an-36a}[2n+9]
-t^{-12an-52a}[2n+7]+t^{-16an-64a}[2n+5]
-t^{-20an-72a}[2n+3]\right)D_5}{[n+5]\left([2n+5]-t^{-4an-8a}[2n+3]\right)}
\\&-(\sm^{-3a}-\sm^{-5a})\left(\sm^{-8a}+\sm^{-4a}\right)D_4
=0.
\end{align*}
That is, we have the following two equations:
{\small
\begin{align*}&
\left(-\sm^{-10a}-\sm^{-6a}
+\frac{[n+3]\left(t^{-12an-36a}[2n+9]
-t^{-16an-52a}[2n+7]+t^{-20an-64a}[2n+5]
-t^{-24an-72a}[2n+3]\right)}{[n+5]\left([2n+5]-t^{-4an-8a}[2n+3]\right)}\right)D_5=0,
\\\\&
\Big[\sm^{-15a}+2\sm^{-11a}-\sm^{-13a}
+2\sm^{-7a}
-2\sm^{-9a}
+\sm^{-3a}-\sm^{-5a}\\
&-(\sm^{-7a}
+\sm^{-3a}-\sm^{-5a})
\frac{[n+3]\left(t^{-8an-36a}[2n+9]
-t^{-12an-52a}[2n+7]+t^{-16an-64a}[2n+5]
-t^{-20an-72a}[2n+3]\right)}{[n+5]\left([2n+5]-t^{-4an-8a}[2n+3]\right)}\Big]D_5\\
&+\left(-\sm^{-9a}
+\sm^{-7a}
+\sm^{-3a}-\sm^{-5a}\right)D_4
=0.
\end{align*}}
Just need to show that the coefficient of $D_5$ in the first equation  is non-zero.
For then $D_5=0$ and  the second equation shows $D_4=0$ since the coefficient of $D_4$ is obviously non-zero.
In turn $D_ 3=D_2=0$,
and then  from $E_1=0$ and $E_2=0$, we have $D_1=D_0=0$.

Indeed,
\begin{align*}&
-\sm^{-10a}-\sm^{-6a}
+\frac{[n+3]\left(t^{-12an-36a}[2n+9]
-t^{-16an-52a}[2n+7]+t^{-20an-64a}[2n+5]
-t^{-24an-72a}[2n+3]\right)}{[n+5]\left([2n+5]-t^{-4an-8a}[2n+3]\right)}
\\&=\Big[(-\sm^{-10a}-\sm^{-6a})[n+5]\left([2n+5]-t^{-4an-8a}[2n+3]\right)
\\&\;\;\;\;\;\;+[n+3]\left(t^{-12an-36a}[2n+9]
-t^{-16an-52a}[2n+7]+t^{-20an-64a}[2n+5]
-t^{-24an-72a}[2n+3]\right)\Big]\\
&\;\;\;\;\;\;\Big/\Big[[n+5]\left([2n+5]-t^{-4an-8a}[2n+3]\right)\Big]\ne 0
\end{align*}
as it is easy to see that the lowest degree in $t$ (treating $\sm$ as $t^{2n}$)
of (the numerator multiplied by $(t^2-t^{-2})^2$) is $-24an-6n-72a-12<0$.

\section{Case  $|p|>q=2$, $|a|>b=2$, $p\ne a$ and $p,a$ have the same sign}\label{sec:2p2a}

The procedure is similar to that   of Section \ref{sec:pqab}; the only difference is that a different set of formulas  is applied. Applying (\ref{connected sum for CJ}) and (\ref{(p,2)}), we have
\begin{equation}\label{equ:p2a2}
\begin{aligned}&
[n+1]J_{T(p,2)\#T(a,2)}(n+1)
=J_{T(p,2)}(n+1)J_{T(a,2)}(n+1)
\\&
=\left( -t^{-4pn-2p}J_{T(p,2)}(n)+t^{-2pn}[2n+1]\right)
\left( -t^{-4an-2a}J_{T(a,2)}(n)+t^{-2an}[2n+1]\right)\\
&=t^{-4(p+a)n-2(p+a)}J_{T(p,2)}(n)J_{T(a,2)}(n)
-t^{-(4p+2a)n-2p}[2n+1]J_{T(p,2)}(n)
\\&
\;\;\;-t^{-(4a+2p)n-2a}[2n+1]J_{T(a,2)}(n)+t^{-2(p+a)n}[2n+1]^2\\&
=t^{-4(p+a)n-2(p+a)}[n]J_{T(p,2)\#T(a,2)}(n)
-t^{-(4p+2a)n-2p}[2n+1]J_{T(p,2)}(n)
\\&
\;\;\;-t^{-(4a+2p)n-2a}[2n+1]J_{T(a,2)}(n)+t^{-2(p+a)n}[2n+1]^2
\end{aligned}
\end{equation}from which  we get
\begin{align*}
&\left([n+1]\sl- \sm^{-2(p+a)}t^{-2(p+a)}[n]\right)J_{T(p,2)\#T(a,2)}(n)\\
&=-\sm^{-2p-a}t^{-2p}[2n+1]J_{T(p,2)}(n)-\sm^{-2a-p}t^{-2a}[2n+1]J_{T(a,2)}(n)+\sm^{-(p+a)}[2n+1]^2.
\end{align*}
Applying the operator $-(\sl+\sm^{-2p}t^{-2p})\frac{\sm^{2p+a}t^{2p}}{[2n+1]}$ from the left to
both sides of this equality and then applying (\ref{(p,2,n+1)}), (\ref{equ:operator}),  and (\ref{(p,2)}), we have
\begin{equation}\label{2p2a}
\begin{aligned}
&-(\sl+\sm^{-2p}t^{-2p})\frac{\sm^{2p+a}t^{2p}}{[2n+1]}\left([n+1]\sl- \sm^{-2(p+a)}t^{-2(p+a)}[n]\right)J_{T(p,2)\#T(a,2)}( n)\\
&=(\sl+\sm^{-2p}t^{-2p})
\left(J_{T(p,2)}(n)
+\sm^{p-a}t^{2p-2a}J_{T(a,2)}(n)-\sm^pt^{2p}[2n+1]\right)
\\&=\sm^{-p}[2n+1]+\sm^{p-a}t^{4p-4a}J_{T(a,2)}(n+1)-\sm^pt^{4p}[2n+3]\\
&\;\;\;\;
+\sm^{-p-a}t^{-2a}J_{T(a,2)}(n)-\sm^{-p}[2n+1]\\
&
=\sm^{p-a}t^{4p-4a}\left( -t^{-4an-2a}J_{T(a,2)}(n)+t^{-2an}[2n+1]\right)
-\sm^pt^{4p}[2n+3]+\sm^{-p-a}t^{-2a}J_{T(a,2)}(n)\\
&=(\sm^{-p-a}t^{-2a}-\sm^{p-3a}t^{4p-6a})J_{T(a,2)}(n)
+\sm^{p-2a}t^{4p-4a}[2n+1]
-\sm^{p} t^{4p}[2n+3].
\end{aligned}
\end{equation}
The function $\sm^{-p-a}t^{-2a}-\sm^{p-3a}t^{4p-6a}$ is non-zero since $p\ne a$.
So we may apply from the left the operator
 \[(\sl+\sm^{-2a}t^{-2a})\frac{1}{(\sm^{-p-a}t^{-2a}-\sm^{p-3a}t^{4p-6a})}\]
  to the preceding equation, and obtain
\begin{align*}
&-(\sl+\sm^{-2a}t^{-2a})\frac{1}{(\sm^{-p-a}t^{-2a}-\sm^{p-3a}t^{4p-6a})}(\sl+\sm^{-2p}t^{-2p})
\frac{\sm^{2p+a}t^{2p}}{[2n+1]}\\
&\left([n+1]\sl- \sm^{-2(p+a)}t^{-2p-2a}[n]\right)J_{T(p,2)\#T(a,2)}(n)
\\
&=(\sl+\sm^{-2a}t^{-2a})\left(J_{T(a,2)}(n)
+\frac{\sm^{p-2a}t^{4p-4a}[2n+1]
-\sm^{p} t^{4p}[2n+3]}{\sm^{-p-a}t^{-2a}-\sm^{p-3a}t^{4p-6a}}
\right)\\
&=\sm^{-a}[2n+1]+\frac{\sm^{p-2a}t^{6p-8a}[2n+3]
-\sm^{p} t^{6p}[2n+5]}{\sm^{-p-a}t^{-2p-4a}
-\sm^{p-3a}t^{6p-12a}}\\&
\;\;\;+\frac{\sm^{p-4a}t^{4p-6a}[2n+1]
-\sm^{p-2a} t^{4p-2a}[2n+3]}{\sm^{-p-a}t^{-2a}-\sm^{p-3a}t^{4p-6a}}
\\&
=\sm^{-a}\frac{\sm^2t^2-\sm^{-2}t^{-2}}
{t^2-t^{-2}}+
\frac{\sm^{p-2a}t^{6p-8a}(\sm^2t^6-\sm^{-2}t^{-6})
-\sm^{p} t^{6p}(\sm^2t^{10}-\sm^{-2}t^{-10})}{(t^2-t^{-2})(\sm^{-p-a}t^{-2p-4a}
-\sm^{p-3a}t^{6p-12a})}\\&
\;\;\;+\frac{\sm^{p-4a}t^{4p-6a}(\sm^2t^2-\sm^{-2}t^{-2})
-\sm^{p-2a} t^{4p-2a}(\sm^2t^6-\sm^{-2}t^{-6})}
{(t^2-t^{-2})(\sm^{-p-a}t^{-2a}-\sm^{p-3a}t^{4p-6a})}.
\end{align*}
Let $g(t, M)$ be the function on the right-hand side  of the last equality.
Then
\begin{align*}\lim_{t\ra-1}(t^2-t^{-1})g(t, \sm)&=\sm^{-a}(\sm^2-\sm^{-2})
+
\frac{\sm^{p-2a}(\sm^2-\sm^{-2})
-\sm^{p}(\sm^2-\sm^{-2})}{\sm^{-p-a}
-\sm^{p-3a}}\\&
\;\;\;+\frac{\sm^{p-4a}(\sm^2-\sm^{-2})
-\sm^{p-2a} (\sm^2-\sm^{-2})}
{\sm^{-p-a}-\sm^{p-3a}}
\\&=\frac{(\sm^2-\sm^{-2})(\sm^{-p-2a}-\sm^{p})}
{\sm^{-p-a}-\sm^{p-3a}}\end{align*}
is a well-defined, non-zero function of $\sm$.
In particular  $(t^2-t^{-2})g(t, \sm)$ is a non-zero function.
Hence
\begin{align*}
\a(t,\sm, \sl)&=(\sl-1)\frac{1}{(t^2-t^{-2})g(t,\sm)}(\sl+\sm^{-2a}t^{-2a})\frac{1}{(\sm^{-p-a}t^{-2a}-\sm^{p-3a}t^{4p-6a})}\\&\;\;\;\;(\sl+\sm^{-2p}t^{-2p})
\frac{\sm^{2p+a}t^{2p}}{[2n+1]}
\left([n+1]\sl- \sm^{-2(p+a)}t^{-2p-2a}[n]\right)
\end{align*}
is an annihilator of $J_{T(p,2)\#T(a,2)}(n)$.
Also $$\a(-1, \sm,\sl)=a(\sm)(\sl-1)(\sl+\sm^{-2a})(\sl+\sm^{-2p})(\sl-\sm^{-2p-2a})
=a(\sm)\sm^{-4p-4a}A_{T(p,2)\#T(a,2)}(\sm,\sl)$$
where \[a(\sm)=\frac{\sm^{2p+a}}{(\sm+\sm^{-1})(\sm^2-\sm^{-2})(\sm^{-p-2a}-\sm^p)}
=\frac{\sm^{p+a}}{(\sm+\sm^{-1})(\sm^2-\sm^{-2})(\sm^{-2p-2a}-1)}\]
is a well-defined, non-zero function of $\sm$.

We now verify that $4$ is  the minimal $\sl$-degree of all non-zero annihilators of $J_{T(p,2)\#T(p,2)}(n)$.
So  suppose
\[D_3J_{T(p,2)\# T(a,2)}(n+3)+D_2J_{T(p,2)\# T(a,2)}(n+2) + D_1 J_{T(p,2)\# T(a,2)}(n+1)
 + D_0 J_{T(p,2)\# T(a,2)}(n)= 0\]
with $D_i\in \q(t,\sm)$, $i=0,1,2$,  we need to  show that $D_i = 0$ for $i = 0,1,2$.

Applying (\ref{equ:p2a2}), (\ref{(p,2,n+2)}) and (\ref{(p,2)}), we have
\begin{align*}&
D_3J_{T(p,2)\# T(a,2)}(n+3)+D_2J_{T(p,2)\# T(a,2)}(n+2) + D_1 J_{T(p,2)\# T(a,2)}(n+1)
 + D_0 J_{T(p,2)\# T(a,2)}(n)\\&=\frac{D_3}{[n+3]}\left(t^{-4(p+a)(n+2)-2(p+a)}[n+2]J_{T(p,2)\#T(a,2)}(n+2)
-t^{-(4p+2a)(n+2)-2p}[2n+5]J_{T(p,2)}(n+2)\right.\\&\;\;\;\;\;\;\;\;\;\;\;\;\;\;
\left.-t^{-(4a+2p)(n+2)-2a}[2n+5]J_{T(a,2)}(n+2)+t^{-2(p+a)(n+2)}[2n+5]^2\right)
\\&\;\;\;\;+D_2J_{T(p,2)\# T(a,2)}(n+2) + D_1 J_{T(p,2)\# T(a,2)}(n+1)
 + D_0 J_{T(p,2)\# T(a,2)}(n)
\\&=\left(\frac{D_3}{[n+3]}t^{-4(p+a)n-10(p+a)}+\frac{D_2}{[n+2]}\right)
\left(t^{-4(p+a)(n+1)-2(p+a)}[n+1]J_{T(p,2)\#T(a,2)}(n+1)\right.\\&\;\;\;\;\;
\left.-t^{-(4p+2a)(n+1)-2p}[2n+3]J_{T(p,2)}(n+1)-t^{-(4a+2p)(n+1)-2a}[2n+3]J_{T(a,2)}(n+1)+t^{-2(p+a)(n+1)}[2n+3]^2
\right)\\
&\;\;\;\;-\frac{D_3}{[n+3]}t^{-(4p+2a)(n+2)-2p}[2n+5]\left(-t^{-8p(n+1)}J_{T(p,2)}(n)-t^{-6p(n+1)}[2n+1]+t^{-2p(n+1)}[2n+3]\right)
\\&\;\;\;\;
-\frac{D_3}{[n+3]}t^{-(4a+2p)(n+2)-2a}[2n+5]\left(-t^{-8a(n+1)}J_{T(a,2)}(n)-t^{-6a(n+1)}[2n+1]+t^{-2a(n+1)}[2n+3]\right)\\&\;\;\;\;
+\frac{D_3}{[n+3]}t^{-2(p+a)(n+2)}[2n+5]^2
+D_1 J_{T(p,2)\# T(a,2)}(n+1)
 + D_0 J_{T(p,2)\# T(a,2)}(n)
\\&
=\left(\frac{D_3}{[n+3]}t^{-8(p+a)n-16(p+a)}+\frac{D_2}{[n+2]}t^{-4(p+a)n-6(p+a)}+\frac{D_1}{[n+1]}\right)
\left(t^{-4(p+a)n-2(p+a)}[n]J_{T(p,2)\#T(a,2)}(n)\right.
\\&
\;\;\;\;\;\;\;\;\left.-t^{-(4p+2a)n-2p}[2n+1]J_{T(p,2)}(n)
-t^{-(4a+2p)n-2a}[2n+1]J_{T(a,2)}(n)+t^{-2(p+a)n}[2n+1]^2\right)\\
&\;\;\;\;-\left(\frac{D_3}{[n+3]}t^{-4(p+a)n-10(p+a)}+\frac{D_2}{[n+2]}\right)
t^{-(4p+2a)(n+1)-2p}[2n+3]\left(-t^{-4pn-2p}J_{(p,2)}(n)+t^{-2pn}[2n+1]\right)
\\&
\;\;\;\;-\left(\frac{D_3}{[n+3]}t^{-4(p+a)n-10(p+a)}+\frac{D_2}{[n+2]}\right)
t^{-(4a+2p)(n+1)-2a}[2n+3]\left(-t^{-4an-2a}J_{(a,2)}(n)+t^{-2an}[2n+1]\right)\\
&\;\;\;\;
+\left(\frac{D_3}{[n+3]}t^{-4(p+a)n-10(p+a)}+\frac{D_2}{[n+2]}\right)t^{-2(p+a)(n+1)}[2n+3]^2
\\&\;\;\;\;-\frac{D_3}{[n+3]}t^{-(4p+2a)(n+2)-2p}[2n+5]\left(-t^{-8p(n+1)}J_{T(p,2)}(n)-t^{-6p(n+1)}[2n+1]+t^{-2p(n+1)}[2n+3]\right)
\\&\;\;\;\;
-\frac{D_3}{[n+3]}t^{-(4a+2p)(n+2)-2a}[2n+5]\left(-t^{-8a(n+1)}J_{T(a,2)}(n)-t^{-6a(n+1)}[2n+1]+t^{-2a(n+1)}[2n+3]\right)\\&\;\;\;\;
+\frac{D_3}{[n+3]}t^{-2(p+a)(n+2)}[2n+5]^2
 + D_0 J_{T(p,2)\# T(a,2)}(n)\\&
=\left[\left(\frac{D_3}{[n+3]}t^{-8(p+a)n-16(p+a)}+\frac{D_2}{[n+2]}t^{-4(p+a)n-6(p+a)}+\frac{D_1}{[n+1]}\right)
t^{-4(p+a)n-2(p+a)}[n]+D_0\right] J_{T(p,2)\# T(a,2)}(n)
\\&\;\;\;\;
+\left[-\left(\frac{D_3}{[n+3]}t^{-8(p+a)n-16(p+a)}+\frac{D_2}{[n+2]}t^{-4(p+a)n-6(p+a)}+\frac{D_1}{[n+1]}\right)
t^{-(4p+2a)n-2p}[2n+1]\right.\\&\;\;\;\;\;\;\;\;\;\;
+\left(\frac{D_3}{[n+3]}t^{-4(p+a)n-10(p+a)}+\frac{D_2}{[n+2]}\right)
t^{-(4p+2a)(n+1)-2p}[2n+3]t^{-4pn-2p}\\&
\;\;\;\;\;\;\;\;\;\;\left.+\frac{D_3}{[n+3]}t^{-(4p+2a)(n+2)-2p}[2n+5]t^{-8p(n+1)}\right]J_{T(p,2)}(n)
\\&\;\;\;\;
+\left[-\left(\frac{D_3}{[n+3]}t^{-8(p+a)n-16(p+a)}+\frac{D_2}{[n+2]}t^{-4(p+a)n-6(p+a)}+\frac{D_1}{[n+1]}\right)
t^{-(4a+2p)n-2a}[2n+1]\right.\\&\;\;\;\;\;\;\;\;\;\;
+\left(\frac{D_3}{[n+3]}t^{-4(p+a)n-10(p+a)}+\frac{D_2}{[n+2]}\right)
t^{-(4a+2p)(n+1)-2a}[2n+3]t^{-4an-2a}\\&
\;\;\;\;\;\;\;\;\;\;\left.+\frac{D_3}{[n+3]}t^{-(4a+2p)(n+2)-2a}[2n+5]t^{-8a(n+1)}\right]J_{T(a,2)}(n)
\\&\;\;\;\;
+\left(\frac{D_3}{[n+3]}t^{-8(p+a)n-16(p+a)}+\frac{D_2}{[n+2]}t^{-4(p+a)n-6(p+a)}+\frac{D_1}{[n+1]}\right)
t^{-2(p+a)n}[2n+1]^2
\\&\;\;\;\;
-\left(\frac{D_3}{[n+3]}t^{-4(p+a)n-10(p+a)}+\frac{D_2}{[n+2]}\right)t^{-(4p+2a)(n+1)-2p}[2n+3]t^{-2pn}[2n+1]
\\&\;\;\;\;
-\left(\frac{D_3}{[n+3]}t^{-4(p+a)n-10(p+a)}+\frac{D_2}{[n+2]}\right)t^{-(4a+2p)(n+1)-2a}[2n+3]t^{-2an}[2n+1]
\\&
\;\;\;\;
+\left(\frac{D_3}{[n+3]}t^{-4(p+a)n-10(p+a)}+\frac{D_2}{[n+2]}\right)t^{-2(p+a)(n+1)}[2n+3]^2
\\&\;\;\;\;-\frac{D_3}{[n+3]}t^{-(4p+2a)(n+2)-2p}[2n+5]\left(-t^{-6p(n+1)}[2n+1]+t^{-2p(n+1)}[2n+3]\right)
\\&\;\;\;\;-\frac{D_3}{[n+3]}t^{-(4a+2p)(n+2)-2a}[2n+5]\left(-t^{-6a(n+1)}[2n+1]+t^{-2a(n+1)}[2n+3]\right)
\\&\;\;\;\;+\frac{D_3}{[n+3]}t^{-2(p+a)(n+2)}[2n+5]^2\\&
=E_3J_{T(p,2)\# T(a,2)}(n)+E_2J_{T(p,2)}(n)+E_1J_{T(a,2)}(n)+E_0.\end{align*}
\begin{lemma}\label{lem:Ei=0p2a2}The equation
$E_3J_{T(p,2)\# T(a,2)}(n)+E_2J_{T(p,2)}(n)+E_1J_{T(a,2)}(n)+E_0=0$ implies
$E_i=0$ for all $i\in\{ 0,1,2,3\}.$
\end{lemma}
This lemma can be proved similarly to Lemma \ref{lem:Eizero}.
Now we show that $E_i=0$ for all $\in \{0,1,2,3\}$  implies $D_i=0$ for all $i\in\{0,1,2,3\}$.
The method is similar to that given in Section \ref{sec:pqab}.
\begin{align*}
&E_3=\left(\frac{D_3}{[n+3]}t^{-8(p+a)n-16(p+a)}+\frac{D_2}{[n+2]}t^{-4(p+a)n-6(p+a)}+\frac{D_1}{[n+1]}\right)
t^{-4(p+a)n-2(p+a)}[n]+D_0=0\\&
\lra\;\;\sm^{-6(p+q)}D_3+\sm^{-4(p+a)}D_2+\sm^{-2(p+a)}D_1+D_0=0,\\&
E_2=-\left(\frac{D_3}{[n+3]}t^{-8(p+a)n-16(p+a)}+\frac{D_2}{[n+2]}t^{-4(p+a)n-6(p+a)}+\frac{D_1}{[n+1]}\right)
t^{-(4p+2a)n-2p}[2n+1]\\&\;\;\;\;\;\;\;\;\;\;
+\left(\frac{D_3}{[n+3]}t^{-4(p+a)n-10(p+a)}+\frac{D_2}{[n+2]}\right)
t^{-(4p+2a)(n+1)-2p}[2n+3]t^{-4pn-2p}\\&
\;\;\;\;\;\;\;\;\;\;+\frac{D_3}{[n+3]}t^{-(4p+2a)(n+2)-2p}[2n+5]t^{-8p(n+1)}
\\&
\lra \;\;\left(-\sm^{-4p-4a}+\sm^{-4p-2a}+\sm^{-4p}\right)D_3+\left(-\sm^{-2p-2a}+\sm^{-2p}\right)D_2+D_1=0,\\&
E_1=-\left(\frac{D_3}{[n+3]}t^{-8(p+a)n-16(p+a)}+\frac{D_2}{[n+2]}t^{-4(p+a)n-6(p+a)}+\frac{D_1}{[n+1]}\right)
t^{-(4a+2p)n-2a}[2n+1]\\&\;\;\;\;\;\;\;\;\;\;
+\left(\frac{D_3}{[n+3]}t^{-4(p+a)n-10(p+a)}+\frac{D_2}{[n+2]}\right)
t^{-(4a+2p)(n+1)-2a}[2n+3]t^{-4an-2a}\\&
\;\;\;\;\;\;\;\;\;\;+\frac{D_3}{[n+3]}t^{-(4a+2p)(n+2)-2a}[2n+5]t^{-8a(n+1)}
\\&
\lra \;\;\left(-\sm^{-4a-4p}+\sm^{-4a-2p}+\sm^{-4a}\right)D_3+\left(-\sm^{-2a-2p}+\sm^{-2a}\right)D_2+D_1=0,
\\&
E_0=\left(\frac{D_3}{[n+3]}t^{-8(p+a)n-16(p+a)}+\frac{D_2}{[n+2]}t^{-4(p+a)n-6(p+a)}+\frac{D_1}{[n+1]}\right)
t^{-2(p+a)n}[2n+1]^2
\\&\;\;\;\;
-\left(\frac{D_3}{[n+3]}t^{-4(p+a)n-10(p+a)}+\frac{D_2}{[n+2]}\right)t^{-(4p+2a)(n+1)-2p}[2n+3]t^{-2pn}[2n+1]
\\&\;\;\;\;
-\left(\frac{D_3}{[n+3]}t^{-4(p+a)n-10(p+a)}+\frac{D_2}{[n+2]}\right)t^{-(4a+2p)(n+1)-2a}[2n+3]t^{-2an}[2n+1]
\\&
\;\;\;\;
+\left(\frac{D_3}{[n+3]}t^{-4(p+a)n-10(p+a)}+\frac{D_2}{[n+2]}\right)t^{-2(p+a)(n+1)}[2n+3]^2
\\&\;\;\;\;-\frac{D_3}{[n+3]}t^{-(4p+2a)(n+2)-2p}[2n+5]\left(-t^{-6p(n+1)}[2n+1]+t^{-2p(n+1)}[2n+3]\right)
\\&\;\;\;\;-\frac{D_3}{[n+3]}t^{-(4a+2p)(n+2)-2a}[2n+5]\left(-t^{-6a(n+1)}[2n+1]+t^{-2a(n+1)}[2n+3]\right)
\\&\;\;\;\;+\frac{D_3}{[n+3]}t^{-2(p+a)(n+2)}[2n+5]^2
\\&\lra\;\;
\left(\sm^{-4p-4a}-\sm^{-4p-2a}-\sm^{-4a-2p}+\sm^{-2p-2a}+\sm^{-4p}-\sm^{-2p}
+\sm^{-4a}-\sm^{-2a}+1\right)D_3\\&
\;\;\;\;\;\;\;\;+\left(\sm^{-2p-2a}-\sm^{-2p}-\sm^{-2a}+1\right)D_2+D_1=0\end{align*}
The determinant of the last three simplified linear equations is
\begin{align*}
-2\sm^{-6p-2a}+\sm^{-4p-2a}+2\sm^{-6p}-2\sm^{-4p}
+2\sm^{-2p-6a}-2\sm^{-6a}+2\sm^{-4a}
+\sm^{-2p}-\sm^{-2p-4a}
-\sm^{-2a}\end{align*}
which is non-zero since $p\ne a$, Hence $D_3=D_2=D_1=0$, and then from $E_3=0$, we get $D_0=0$.

 \section{Case   $|p|>q=2$, $|a|>b=2$, $p= a$}
This  case
 quickly follows from what has been done in Section \ref{sec:2p2a}.
Equation (\ref{2p2a}) remains true  when $p=a$, yielding
\begin{align*}
&-(\sl+\sm^{-2p}t^{-2p})
\frac{\sm^{3p}t^{2p}}{[2n+1]}\left([n+1]\sl- \sm^{-4p}t^{-4p}[n]\right)J_{T(p,2)\#T(p,2)}(n)\\
&=\sm^{-p}[2n+1]-\sm^pt^{4p}[2n+3].
\end{align*}
Let $g(t,\sm)$ be the function on the right-hand side  of the last equality, which is non-zero.
Hence
\[\a(t, \sm, \sl)=(\sl-1)\frac{1}{(t^2-t^{-2})g(t,\sm,\sl)}(\sl+\sm^{-2p}t^{-2p})
\frac{\sm^{3p}t^{2p}}{[2n+1]}\left([n+1]\sl- \sm^{-4p}t^{-4p}[n]\right)\]
is an annihilator of $J_{T(p,2)\#T(p,2)}(n)$.
Also\[\a(-1,\sm,\sl)=a(\sm)(\sl-1)(\sl+\sm^{-2p})(\sl-\sm^{-4p})=a(\sm)\sm^{-6p}A_{T(p,2)\#T(p,2)}(\sm,\sl)\]
where \[a(\sm)=\frac{\sm^{3p}}{(\sm^2-\sm^{-2})(\sm^{-p}-\sm^p)(\sm+\sm^{-1})}.\]
What remaining to show is that if
\[D_2J_{T(p,2)\# T(p,2)}(n+2) + D_1 J_{T(p,2)\# T(p,2)}(n+1)
 + D_0 J_{T(p,2)\# T(p,2)}(n)= 0.\]
for some $D_0,D_1,D_2\in \q(t,\sm)$, then we must have  $D_i = 0$ for all $i$.

Plugging $D_3=0$ and $a=p$ into the corresponding calculation in Section \ref{sec:2p2a}, we obtain
\begin{align*}&
D_2J_{T(p,2)\# T(p,2)}(n+2) + D_1 J_{T(p,2)\# T(p,2)}(n+1)
 + D_0 J_{T(p,2)\# T(p,2)}(n)\\&
=\left[\left(\frac{D_2}{[n+2]}
t^{-8pn-12p}+\frac{D_1}{[n+1]}\right)
t^{-8pn-4p}[n]+D_0\right]
J_{T(p,2)\# T(p,2)}(n)\\&\;\;\;\;\;\;
+\left[-\left(\frac{D_2}{[n+2]}
t^{-8pn-12p}+\frac{D_1}{[n+1]}\right)2t^{-6pn-2p}[2n+1]+\frac{D_2}{[n+2]}2t^{-6pn-8p}[2n+3]t^{-4pn-2p}
\right]J_{T(p,2)}(n)\\&\;\;\;\;\;\;
+\left(\frac{D_2}{[n+2]}
t^{-8pn-12p}+\frac{D_1}{[n+1]}\right)t^{-4pn}[2n+1]^2\\&\;\;\;\;\;\;
-\frac{D_2}{[n+2]}2t^{-6pn-8p}[2n+3]t^{-2pn}[2n+1]
+\frac{D_2}{[n+2]}t^{-4pn-4p}[2n+3]^2
\\&\;\;=E_2J_{T(p,2)\# T(p,2)}(n)+E_1J_{T(p,2)}(n)+E_0.\end{align*}
\begin{lemma}The equation  $E_2J_{T(p,2)\# T(p,2)}(n)+E_1J_{T(p,2)}(n)+E_0=0$
implies $E_i=0$ for all $i=0,1,2$.
\end{lemma}
This lemma can be treated as a subcase of Lemma \ref{lem:Ei=0p2a2}.

 Now we show that the system of linear equations $E_i=0$, $i=0,1,2$,
in variables $D_2, D_1, D_0$ has only the trivial solution.
Again plugging $D_3=0$ and $a=p$ into the relevant calculation in Section \ref{sec:2p2a},
 we obtain the following simplified linear system in variables $D_2$, $D_1, D_0$:
\begin{align*}
&\sm^{-8p}D_2+\sm^{-4p}D_1+D_0=0,\\&
\left(-\sm^{-4p}+\sm^{-2p}\right)D_2+D_1=0,\\&
\left(\sm^{-4p}-2\sm^{-2p}+1\right)D_2+D_1=0.\end{align*}
The determinant of this system is $-2\sm^{-4p}+3\sm^{-2p}-1\ne 0$.
Hence $D_2=D_1=D_0=0$.

\end{document}